\documentclass[10pt,a4paper,reqno]{book}

\usepackage[english]{babel}

\usepackage[latin1]{inputenc}
\RequirePackage{slashed}
\usepackage{makeidx}
\usepackage[T1]{fontenc}
\usepackage{aeguill}
\usepackage{amsmath,amsthm, amssymb}
\usepackage{amsfonts}
\usepackage{color}
\usepackage{graphics}
\usepackage[all]{xy}
\usepackage{epsfig}
\usepackage[left=2.5cm,top=2.5cm,right=2.5cm,bottom=2.5cm,bindingoffset=0cm]{geometry}
\usepackage{hyperref}
\hypersetup{
pdfstartview=FitH,
backref=true, 
pagebackref=true,
hyperindex=true, 
colorlinks=false, 
breaklinks=true, 
urlcolor=black, 
linkcolor=black, 
bookmarks=true, 
pdfpagemode=none,
bookmarksopen=true, 
pdftitle={Vertex Operator Algebra}, 
pdfauthor={Christophe NOZARADAN}, 
pdfsubject={Mathématique} 
}

\newcommand{\ket}[1]{|#1\rangle}
\newcommand{\braket}[2]{\langle #1|#2\rangle}

\theoremstyle{plain}
\newtheorem{theo}{Theorem}

\newtheorem{coroll}[theo]{Corollary}
\theoremstyle{definition}\newtheorem{defin}[theo]{Definition}
\newtheorem{exe}[theo]{Example}

\theoremstyle{plain}\newtheorem{lemme}[theo]{Lemma}

\newtheorem{prop}[theo]{Proposition}
\theoremstyle{remark}\newtheorem{remarque}[theo]{Remark}

\numberwithin{theo}{section}

\newcommand{\F}[3]{F_{#1,#2}^{#3}}
\newcommand{\Fzd}[2]{F_{#1}^{#2}}
\newcommand{\Fz}{F_{z,w}^{\lambda}}
\newcommand{\Fzz}{F_{z}^{\lambda}}

\newcommand{\Res}[1]{\text{Res}_{#1}}
\newcommand{\beqn}{\begin{eqnarray}}
\newcommand{\eeqn}{\end{eqnarray}}
\newcommand{\beqns}{\begin{eqnarray*}}
\newcommand{\eeqns}{\end{eqnarray*}}

\newcommand{\rmq}[1]{}
\newcommand{\corps}[1]{\mathbb{#1}}

\title{\textsf{\textbf{Introduction to Vertex Algebras}}}
\author{\textsf{\textbf{Christophe NOZARADAN$^\dagger$}}}
\date{
\vskip 1truecm
        {\it\small   {Mathematics department, \\
        Université Catholique de Louvain, \\ 
        Louvain-la-Neuve, Belgium} \\
        [3mm]\small e-mail:} {\small\tt christophe.nozaradan@uclouvain.be}
        \vskip 1truecm
        {\small\textbf{Abstract}}
        \vskip 0.5truecm
        \begin{changemargin}{1.5cm}{1.5cm}
\noindent These lecture notes are intended to give a modest impulse to anyone willing to start or pursue a journey into the theory of \textit{vertex algebras} by reading one of the books \cite{kac1} or \cite{lepow}.  Therefore, the primary goal is to provide required tools and help being acquainted with the \textit{machinery} which the whole theory is based on.  
The exposition follows Kac's approach~\cite{kac1}. Fundamental examples relevant in Theoretical Physics are also discussed.  No particular prerequisites are assumed.
\vskip 1truecm
\small Last edited on \today.
\vskip 2cm
$^\dagger$ Research Fellow F.R.S. - FNRS (Aspirant)
\end{changemargin}
} 

\makeindex

\begin{document}
\numberwithin{equation}{chapter}
\frontmatter

\maketitle

\chapter*{Prologue}
\begin{changemargin}{1cm}{1cm} 
\baselineskip=15pt

The notion of \textit{vertex algebra}\footnote{Although similar, \textit{vertex algebras} should be distinguished from \textit{vertex operator algebras}, the latter being a slightly more restricted notion than the former (see the main text).} was first axiomatized by Richard Borcherds in 1986.  In retrospect, this algebraic structure has been proved to be equivalent to Conformal Theory of chiral fields in two dimensions~\cite{BPol}, for which it provides a rigorous mathematical formulation.  Thereby, a new area was born, stressing beautifully the interplay between physics and mathematics.   

Originally, vertex operators arose in String Theory. They are used to describe certain types of interactions, between different particles or strings, localized at vertices -- hence the name ``vertex'' -- of the corresponding Feynman diagrams.  From this point of view, the notion of vertex operator thus appeared \textit{before} the underlying concept of vertex algebra.     

The concept of \textit{vertex algebra} happens to be a powerful tool -- and initially introduced by Borcherds to this purpose -- in the proof~\cite{Bo} of the so-called ``Moonshine Monstrous'' conjectures, formulated by Conway and Norton.  This notion is used to build representations of sporadic groups (the Monster group in particular), but has also a prominent role in the infinite-dimensional case (Kac-Moody algebras).  In addition, it is not surprising to see the notion of Vertex Algebra emerging from String Theory, as the worldsheet of a string exhibits the symmetry of a two-dimensional conformal field theory.  Seiberg and Witten~\cite{SW} then made the connection between Vertex Algebras and Connes's Noncommutative Geometry, which can already be shown to appear when studying strings ending in $Dp$-branes.  

Those few motivations illustrate this interdisciplinary dialog, based on a relatively abstract notion, source of a fruitful complicity : the one enhances itself and in turn provides new perspectives to the others.

These present notes do not pretend to contain original material. It is to be considered as an account of the author's personal journey in the theory of \textit{vertex algebras}, that began his last undergraduate year in Physics. The goal intented is to bring this active field of research to the attention of a wider audience, by offering a more readable exposition than it is believed to be available at present (in the author's own physicist point of view), and to give interested readers a modest impulse in their own journey into the theory. They are invited to pursue by reading references such as \cite{kac1}, which is the main bibliography of these notes, or \cite{lepow}. The reference~\cite{kac4} has also been of great help, in addition to stimulating discussions by e-mail with R. Heluani, which I thank for his patience. Even though the best effort has been made to avoid errors or mistakes, the readers are encouraged to be vigilant and report any left ones, from minor typos to any more significant errors.  

\end{changemargin}

\chapter*{Introduction}
\begin{changemargin}{1cm}{1cm} 
\baselineskip=15pt

\textit{This part has been used as an introductory material in a four hours lecture given by the author at the Fourth International Modave Summer School on Mathematical Physics, held in Modave, Belgium, September 2008.}\\

We introduce the subject by reviewing basic aspects of the \textit{Monstrous moonshine conjectures}\index{moonshine conjectures} and in particular the way the concept of vertex algebras arises in this context. We then give a first overlook on its algebraic structure to motivate the organization of the main text.

\section*{Moonshine conjectures}
Historically, the notion of vertex algebra first appeared in the proof (Borcherds, 1992) of the \textit{monstrous moonshine conjectures}, stated by Conway and Norton in 1979.  Those conjectures reveal the surprising (hence the name ``moonshine'') connection between \textit{sporadic simple groups}\index{sporadic groups} on the one hand and \textit{modular functions}\index{modular functions} on the other.  We mention here only some needed facts to state the conjectures (see \cite{gannon} for more details).  

A great deal of work has been necessary to classify all finite dimensional simple groups. It took over thirty years and involved contributions of about one hundred authors. In the end, simple groups have been shown to fall into four classes~: the cyclic groups, the alternating groups of degree $n\geq5$, the Lie-type groups (over a finite field) and finally those not entering any of the other classes~: the \textit{sporadic} groups.   The largest sporadic group is called the \textit{Monster}\index{Monster group}.  It has about $8.10^{53}$ elements.  For further reference, let's note that its smallest non trivial irreducible representations are $196883$, $21296876$ and $842609326$-dimensional.

Among modular functions, we will focus on the so-called $j$-function\index{$j$-function}.  This function is holomorphic on the upper-half complex plane $\cal H$ and invariant under the action of the group $SL_2(\corps{Z})$ on $\cal H$ :
$$ j(z)=j\left(\frac{az+b}{cz+d}\right) \ \ \forall \begin{pmatrix} a & b \\ c & d \end{pmatrix} \in SL_2(\corps{Z}) \ , \ z\in {\cal H} \ ,$$
where 
$SL_2(\corps{Z})$ is the group of two-by-two matrices with entries in $\corps{Z}$ and determinant $1$.
Since $\begin{pmatrix} 1 & 1 \\ 0 & 1 \end{pmatrix} \in SL_2(\corps{Z})$, we see that $j(z+1)=j(z)$ for $z\in{\cal H}$.  This allows us to write the so-called \textit{$q$-expansion} of $j$, its Laurent series (since bounded below) in the origin in which it admits a pole and outside of which it is holomorphic (in the region restricted to the unit disk).  Basically, it looks as follows~:
$$ j(q) = \frac{1}{q} + 744 + 196884 q + 21493760 q^2 + 8642909970 q^3 + \cdots \ ,$$
where $q=e^{2\pi i z}$. Looking at the numerical values of the coefficients of $q$, McKay first noticed that
$$ 196884 = 1 + 196883 \ .$$
Now recall that $196883$ is the dimension of the smallest non trivial irreducible representation of the monster group $\corps{M}$. Considering $1$ as the dimension of its trivial representation, it is tempting to say that the coefficient of $j$ appearing on the LHS coincides with the dimension of a reducible representation (the direct sum of the first two irreducible representations, known as the \textit{Griess algebra}\index{Griess algebra}) of $\corps{M}$ appearing on the RHS.  This is a nice result, but it may just be a lucky coincidence.  Now considering the coefficients of $q^2$ and $q^3$ in the $q$-expansion of $j$, McKay and Thompson then showed that
\beqns
 21493760 &=& 1 + 196883 + 21296876 \\
 864299970 &=& 2 \cdot 1 + 2 \cdot 196883 + 21296876 + 842609326 \ .
\eeqns
There again, a linear combination of dimensions of irreducible representations of $\corps{M}$ appears on the RHS of each equation, so that two other coefficients of $j$ are shown to coincide with the dimension of two other reducible representations of $\corps{M}$. This seems to be more than a coincidence.  Part of the conjecture states that this is true for each coefficient of $j$, the constant term $744$ being kindly put aside.  Since the dimension of a reducible representation $V_n$ coincides with the trace of the matrix representing the unit element $e$ of $\corps{M}$ on $V_n$, we can write
$$ j(z) - 744 = \frac{1}{q}\sum_{n=0}^{\infty} q^n \text{Tr}(e|V_n) \ ,$$
where $V_n$ denote reducible representations of $\corps{M}$, with $V_0$ the trivial representation space, $V_1$ the empty space, $V_2$ the underlying space of the Griess algebra, and so forth.  

The conjecture even goes further. It generalizes the latter equality by replacing $e$ by any element $g$ of $\corps{M}$, for which other modular functions occur, belonging to the \textit{same} family as $j$ (all associated to Riemann surfaces of genus $0$). The generalized equalities then obtained are called the \textit{McKay-Thompson series}\index{McKay-Thompson series}.

The space $$V=\bigoplus_{n=0} V_n$$ constitutes an infinite-dimensional graded module of $\corps{M}$. Acted on by the Monster, it admits in fact a structure of vertex algebra, known as the \textit{Monster vertex algebra}.  The latter notion thus turns out to constitute a ``bridge'' between the study of modular functions and the representations of sporadic simple groups. As one could have argued at that time, it would have been naive to think that this fully illuminates the connection between both fields, since the original question had then been replaced by the -- equally complex -- question of understanding vertex algebraic structures. Nevertheless, it is undeniable that the introduction of the concept of vertex algebra shed some \textit{heaven-sent} light on the contours of this complex matter, as Pandora's box, whose secrets were yet to be discovered. 

\section*{What is a vertex algebra ?}
For the main part of the lecture, we will try to answer the legitimate question~: \textit{What is a vertex algebra?}  The latter is a title borrowed from one of Borcherds' talks. 
There are several definitions -- five to my knowledge -- of the notion of vertex algebra, all equivalent to one another.  Each one is convenient in a particular domain, which are remarkably numerous~: representation theory and modular functions, as seen so far, but also number theory, algebraic geometry, integrable systems and of course mathematical physics, conformal theories as string theory, among other areas. 
Borcherds first defined the notion of vertex algebra in the following way.

\noindent \textbf{Definition 1.} A \textit{vertex algebra} is a vector space $V$ endowed with an \textit{infinite} number of bilinear products.  We index those products by an integer $n$ and write them as bilinear maps from $V\otimes V$ to $V$ such that
$a \otimes b \to a_n b \ \ \forall a,b \in V$, where $a_n b$ is to be understood as the endomorphism $a_n\in End(V)$ acting on $b$ ($n\in\corps{Z}$) :
$$a_n b := a_n(b) \ \ \text{where } a_n\in End(V) \ , \ b\in V \ .$$
Those products are subject to the following conditions :
\begin{enumerate}
\renewcommand{\labelenumi}{(\theenumi)}
\item For all $a,b\in V$, $a_n b = 0$ for large enough $n\in\corps{Z}$;
\item There exists an element $1\in V$ such that $a_n 1 = 0$ for $n\geq 0$ and $a_{-1} 1 = a$;
\item Those products are related to one another by the so-called \textit{Borcherds' identity}\index{Borcherds' identity} :
\begin{eqnarray*}
\sum_{i\in \corps{Z}} \binom{m}{i} \left(a_{q+i}b \right)_{m+n-i}c \quad\quad\quad\quad\quad\quad\quad\quad\quad\quad\quad\quad\quad\quad\quad\quad\quad\quad\\ 
=  \sum_{i\in \corps{Z}} (-1)^i \binom{q}{i} \Big( a_{m+q-i}\left(b_{n+i}c\right) - (-1)^q b_{n+q-i}\left(a_{m+i}c\right)\Big) 
\end{eqnarray*}
\end{enumerate}

This is the most concise way to define a vertex algebra.
The latter -- rather impressive -- identity contains in fact a lot of information and is thus considered as the most important condition for vertex algebras.  A parallel is often made with the Jacobi identity for Lie algebras.

We will put aside Borcherds' identity for now though and focus on the bilinear products.
As seen above, the bilinear product $a_n b$ is to be understood as resulting from an endomorphism $a_n$ of $V$ acting on $b\in V$.  The chosen notation $a_n$ is not an accident, but the exact relation between $a_n \in End(V)$ and $a\in V$ remains unclear at this stage and is yet to be specified.
Now, in order to do so and to deal with an infinite amount of products, it is convenient to gather all of them into a single ``formal'' expansion.  Let's introduce a ``formal'' variable $z$ (in the formal setting, the \textit{indeterminate}\index{indeterminate} $z$ never takes any value of any kind) and gather all the products in one big expansion
$$ \sum_{n\in\corps{Z}} a_{n} b \ z^{n}  =: a(z)b \ ,$$
where 
$$a(z):=\sum_{n\in\corps{Z}} a_{n} z^{n}$$ is called an $End(V)$-valued \textit{formal distribution}\index{formal distribution}.
As a convention, we make the substitution $n\to -n-1$ and define $a_{(n)}:=a_{-n-1}$, so that positive integers index the singular terms of the expansion.  The resulting expansion 
$$a(z) = \sum_{n\in\corps{Z}} a_{(n)} z^{-n-1} $$ is called the \textit{Fourier expansion}\index{Fourier expansion} of $a(z)$ and its coefficients $a_{(n)}\in End(V)$ ($n\in\corps{Z}$) its \textit{Fourier modes}\index{Fourier modes}.
The bilinear products can be recovered from $a(z)b$ since we have
$$a_{(n)} b = \text{Res}_z z^n a(z)b \ . $$
Therefore, $a(z)b$ contains the same information as $\{ a_{(n)} b \}_{n\in\corps{Z}}$ and we can say that a vertex algebra \textit{resembles} an ordinary algebra, except that \textit{its} product law is \textit{parameter-dependant}. As we will see, in a vertex algebra, any element $a\in V$ is uniquely associated to an $End(V)$-valued formal distribution $a(z)$, denoted from now on by $Y(a,z)$. As a consequence, from the algebraic point of view, it is really the \textit{parameter-dependant} product $Y(a,z)b$ that can be seen as the analogue of an ordinary product between two elements $a,b\in V$.  For this reason, this product is often written 
$$a_z b := Y(a,z)b \ .$$

This analogy goes further, as we recall the definition of an algebra showing up as the ``classical limit'' of a vertex algebra.  This algebra is nothing else but a standard \textit{associative, commutative and unital} (i.e. it has a unit element) \textit{algebra}. To stress the analogy with its ``quantum'' counterpart, we present its definition in the following way. Such an algebra can be defined as a linear space $V$ (over $\corps{C}$) endowed with a linear map
$$ Y : V \to End(V) $$
satisfying 
$$ [Y(v),Y(w)] = 0 \ \ \forall v,w\in V \ , $$       
with a unit element $1$ being such that 
$$ Y(1)=I_V \quad\text{and}\quad Y(v)1=v \ \forall v\in V \ . $$
This formulation is completely equivalent to the standard one.  Indeed, defining the product law by 
$$v\cdot w := Y(v)w \ ,$$ for all $v,w \in V$, it is not difficult to verify that the properties of commutativity and associativity are satisfied, and to check that $1$ is the unit element of $V$ with respect to this product, as expected.

Similarly, the structure of vertex algebra can be seen as a linear space $V$ (over $\corps{C}$), except that this time the linear map $Y$ goes from $V$ to the space of $End(V)$-valued formal distributions, which we denote by $\text{End}(V)[[z,z^{-1}]]$~: 
$$ Y(\cdot, z) \,:\, V \to \text{End}(V)[[z,z^{-1}]] \,:\, a \to Y(a,z) \ ,$$
where $Y(a,z)$ is called a \textit{vertex operator}, and in general, for $a,b\in V$,
$$ [Y(a,z),Y(b,z)] \neq 0 \ ,$$
where $[\ ,\ ]$ is the usual Lie superbracket defined for an associative superalgebra.
Moreover, there exists an element $1\in V$ such that 
$$Y(1,z) = \text{I}_V \quad \text{and}\quad Y(a,z)1=a + \cdots \ ,$$
where $\cdots$ denotes an expression involving only strictly positive powers of $z$.
This element $1$ is called the \textit{vacuum vector}\index{vacuum vector} of $V$ and thus often written $\ket{0}$. Notice that, at the level of the \textit{Fourier modes}\index{Fourier modes}, the second property of the vacuum vector translates as follows~:
$$ a_{(n)}1 = 0 \ \text{for}\ n\geq 0 \quad \text{and} \quad a_{(-1)}1=a \ . $$
This actually coincides with item (2) in Borcherds' definition, which should then be written using the notation $a_{(n)}b$ instead of $a_n b$, to agree with our convention. One can thus get a glimpse of how a vertex algebra can be defined by a different but equivalent set of axioms.  This will allow us to better investigate the structure of vertex algebra and, in particular, to realize how much information is contained in the single Borcherds' identity.
The above considerations invite us to start a more serious study of vertex algebraic structures in the framework of \textit{formal distributions} (following Kac's book).  The plan of the lecture notes is exposed in the next section.

\section*{At a glance}
As already mentioned, the concept of \textit{vertex algebra} and its underlying structures are studied following Victor Kac's book \cite{kac1}.

Although it might be considered as a detour, it is believed that the concept of \textit{vertex algebra} is best introduced by first studying the related notion of Lie Conformal Algebra, also called Vertex Lie Algebra.  The properties of the latter are thought to be easier to grasp and, moreover, help being acquainted with the \textit{machinery} which the whole theory is based on. 

The first chapter is organized as follows.  After setting some terminology, we define formal distributions in analogy with Complex Analysis.  We then study a particular formal distribution, the formal Dirac's $\delta$ distribution and derive its main properties.  Among them, one permits us to introduce the notion of \textit{locality}, which constitutes the core of the theory.  Playing a central role in the second part, the notion of \textit{Fourier transform} closes this chapter.

The subsequent chapter deals with algebraic structures formed with the help of formal distributions.  After a brief recall of \textit{Lie superalgebras}, notions of \textit{$j$-products} and \textit{$\lambda$-bracket} are introduced.  They permit us to define the structures of \textit{formal distribution Lie superalgebras} and \textit{Lie conformal superalgebras}, canonically related to one another (\textit{canonical bijection}).  Examples relevant in Theoretical Physics are then discussed. 

In the third chapter, we start by reviewing common aspects of Conformal Field Theory (CFT) and show how they emerge naturally in the framework of local formal distributions. We start with the notions of \textit{Operator Product Expansions} and \textit{Normal Products}.  The latter, together with the $j$-products, are shown to be special instances of a \textit{generalized  $j$-product}. Entering the definition of \textit{vertex algebra} and useful in applications, \textit{Wick's formulas}, as well as the \textit{quasi-commutativity} of the normal product, are proved.  We close the chapter with the notion of \textit{Eigendistribution} of a certain \textit{weight}.

The last chapter finally focuses on the notion of \textit{vertex algebra}.  We introduce its main definition by analogy with its ``classical limit''~: the structure of \textit{associative, commutative and unital algebra}.  We then study its main properties in light of the algebraic structures from the previous chapters and discuss examples. We also survey fundamental theorems and results, such as the \textit{unicity} and the \textit{$n$-products theorems}, including the property of \textit{quasi-symmetry}.  Written in terms of Fourier modes, the latter gives rise to \textit{Borcherds' $n$-products}, which actually implies the property of quasi-commutativity and \textit{quasi-associativity} of the normal product. We then show that every \textit{vertex algebra} is actually a Lie conformal algebra.  We end by illustrating the interplay of all those notions previously introduced, as they enter an equivalent definition of \textit{vertex algebra}~\cite{DsK2}.
\end{changemargin}

\setcounter{tocdepth}{2}
\tableofcontents
\baselineskip=15pt

\mainmatter

%
%

\chapter{Formal calculus}

\section{Formal distributions}

Throughout the chapters, we will adopt the following terminology.
\begin{defin} Let $\cal A$ be any algebraic structure. 
\begin{eqnarray*}
&&{\cal A}[z,w,\ldots] = \Big\lbrace \, a(z,w,\ldots)=\sum_{n,m,\ldots =0}^{N,M,\ldots}  a_{nm\ldots} z^n w^m \ldots \,\, | \,\, a_{nm\ldots} \in {\cal A} \, ; N,M,\ldots \in \mathbb{Z}^+ \Big\rbrace \\
&&{\cal A}[[z,w,\ldots]] = \Big\lbrace \, a(z,w,\ldots)=\sum_{n,m,\ldots \in \mathbb{Z}^+} a_{nm\ldots} z^n w^m \ldots \,\, | \,\, a_{nm\ldots} \in {\cal A} \,  \Big\rbrace \\
&&{\cal A}[[z,z^{-1},w,w^{-1},\ldots ]] = \Big\lbrace \, a(z,w,\ldots)=\sum_{n,m,\ldots \in\mathbb{Z}} a_{nm\ldots} z^n w^m \ldots \,\, | \,\, a_{nm\ldots} \in {\cal A} \,  \Big\rbrace \ ,
\end{eqnarray*}
where ${\cal A}[z,w,\ldots]$, ${\cal A}[[z,w,\ldots]]$ and ${\cal A}[[z,z^{-1},w,w^{-1},\ldots ]]$ denote the spaces of \textbf{$\cal A$-valued formal polynomials}\index{formal polynomials}, \textbf{formal Taylor series}\index{formal Taylor series} and \textbf{formal distributions}\index{formal distributions} in the formal variables $z,w,\ldots$ respectively.
\end{defin}
\begin{remarque}
Other spaces may appear, the definition of which can easily be deduced from the above cases (e.g. ${\cal A}[z,z^{-1}]$, the space of \textbf{formal Laurent polynomials}\index{formal Laurent polynomials}; ${\cal A}[[z]][z^{-1}]$, the space of \textbf{formal Laurent series}\index{formal Laurent series}, in which all but finitely many singular terms are null).  We will mostly be interested in formal expansions in one or two indeterminates.
\end{remarque}

The \textit{formal} formalism is purely algebraic.  In particular, no notion of convergence with respect to a given topology is needed.  The intern laws of addition and multiplication alone are to be well defined.  We say that the algebraic structure of ${\cal A}$ is transferred to the space ${\cal A}[z]$, because the addition and multiplication of any two of its elements lead to a well defined element of ${\cal A}[z]$.  The same is true for the space ${\cal A}[[z]]$.  

On the contrary, this fails to be the case for the space ${\cal A}[[z,z^{-1}]]$.  The addition of two formal distributions is a well defined operation, because each coefficient of the resulting series is an element of ${\cal A}$ by linearity.  Thus, the problem comes from the multiplication law, which is more delicate.  In this case, each coefficient of the resulting series is actually an \textit{infinite} sum of elements of ${\cal A}$.  Indeed, let $a(z)$ and $b(z) \in {\cal A}[[z,z^{-1}]]$.  Their product leads to :
\begin{eqnarray}
a(z) \, \cdot \, b(z) &=& \, \sum_{n \in \corps{Z}} a_n z^n  \, \cdot \, \, \sum_{m \in \corps{Z}} b_m z^m  \\
&=& \sum_{n,m\in \corps{Z}}  a_n \, b_m  z^{n+m} \\
&=& \sum_{k \in \corps{Z}} \left( \sum_{n+m=k} a_n \, b_m \right) z^{k} \ \label{prod}
\end{eqnarray} 
Let's consider the coefficient of $z^k$.  Its expression contains as many terms as couples $(n,m)$ being solutions of the equation $n+m=k$, where $k \in \corps{Z}$ is fixed.  As $n$ and $m \in \corps{Z}$, the equation $n+m=k$ admits an infinite amount of solutions, so that the second sum in \eqref{prod} is infinite.  But in general, nothing guaranties that it converges \textit{in} ${\cal A}$.  Therefore, the result cannot be considered as a formal distribution of ${\cal A}[[z,z^{-1}]]$.  

Note that this type of behavior is analogous to distributions in the sense of Schwartz's theory.  One cannot multiply distributions with one another either.

\begin{defin}
The \textbf{Fourier expansion}\index{Fourier expansion} of any formal distribution $$a(z)=\sum_{n\in\corps{Z}} a_n z^n \ ,$$with $a_n \in {\cal A}$ for all $n\in\corps{Z}$, is conventionally written
$$a(z)=\sum_{n\in\corps{Z}} a_{(n)} z^{-n-1} \ ,$$ 
where the coefficients
$$a_{(n)}=\Res{z}z^n a(z)=a_{-n-1}\in {\cal A} \ ,$$ with $n\in\corps{Z}$, are called the \textbf{Fourier modes}\index{Fourier modes} of $a(z)$.
\end{defin}
The notion of \textit{residue}\index{formal residue} is taken formally in analogy with Complex Analysis, \textit{i.e.} $$\Res{z}a(z)\doteq a_{-1} \ ,$$ or equivalently, 
$$\Res{z}a(z)\doteq a_{(0)} \ .$$  In this convention, the singular terms of a Fourier expansion are labeled by non-negative integers.  The generalization in several indeterminates is straightforward.  The case of more than two variables is useless for our concern.

\begin{remarque}
We can extend the expression of the Fourier modes by linearity.  Replacing $z^n$ by any element $\varphi(z) \in \corps{C}[z,z^{-1}]$, we can define a non degenerate pairing
$$ {\cal A}[[z^{\pm 1}]] \otimes \corps{C}[z,z^{-1}] \to {\cal A} \, : \, a\otimes \varphi \to \braket{a}{\varphi} := \Res{z}\varphi(z) a(z) \ .$$
Any formal distribution $a(z)\in{\cal A}[[z^{\pm 1}]]$ can then be seen as a distribution acting on the space of test functions $\corps{C}[z,z^{-1}]$, hence the name ``distribution''.
\end{remarque}

\begin{remarque} 
If $a(z)\in L^2(S^1)$, we get the standard Fourier expansion of a square integrable function on the circle.  Consider the orthonormal basis 
$$\big\lbrace \frac{e^{in\theta}}{\sqrt{2\pi}}\, ;\ n\in\corps{Z},\, \theta\in[0,2\pi) \big\rbrace$$ of the Hilbert space $L^2([0,2\pi),d\theta)$, endowed with the hermitian scalar product
$$\braket{\phi}{\psi}\doteq \int_{0}^{2\pi} \phi^*(\theta) \psi(\theta) d\theta \ .$$
Straight from the definition, we just have to pose $z=e^{i\theta}$ (so $d\theta=\frac{dz}{iz}$) to have $z^{-1}=z^*$ (where $z^*$ denotes the complex conjugate of $z$).  Explicitly,
\begin{eqnarray*}
a(z)&=&\sum_{n\in\corps{Z}} z^{n} Res_z z^{-n-1} a(z) \\
&=&\sum_{n\in\corps{Z}} z^{n} \frac{1}{2\pi i} \int_{|z|=1} z^{-n-1} a(z) dz \\
&=&\sum_{n\in\corps{Z}} \frac{e^{in\theta}}{\sqrt{2\pi}} \int_{0}^{2\pi} \frac{e^{-in\theta}}{\sqrt{2\pi}} \tilde{a}(\theta) d\theta 
\end{eqnarray*}
\end{remarque}

\section{Notion of expansion}
In Complex Analysis, rational functions can be expanded in several ways, depending on the domain of convergence considered.  Similarly, in the formal setting, a rational expression admits different expansions.

Even if the notion of convergence is of no matter formally, we still have to distinguish the different expansions that might occur for a given rational function, so we introduce the notion of \textbf{formal expansion}\index{formal expansion}. We will only consider rational functions $F(z,w)$ with \textit{poles} at $z=0$, $w=0$ and $z=w$ at worst.  These are finite $\corps{C}$-linear combinations of elements of the form $z^n w^m (z-w)^k$, where $n,m\in\corps{Z}$ and $k\in\corps{Z}^-$, i.e. $F(z,w)\in \corps{C}[z^{\pm 1},w^{\pm 1},(z-w)^{-1}]$. Considering $w\in\corps{C}$ as a given parameter, $F(z,w)$ admits two different expansions.  These are denoted by $i_{z,w}F(z,w)$ and $i_{w,z}F(z,w)$.  This notation is related to its original meaning~: if $z,w\in \corps{C}$, then the map $i_{z,w}$ (resp. $i_{w,z}$), applied on a rational function in $\corps{C}[z^{\pm 1},w^{\pm 1},(z-w)^{-1}]$, gives its converging expansion in the domain $|z|<|w|$ (resp. $|w|<|z|$).  Those expansions will be derive in proposition \ref{propexp} for the only non trivial case, but let's state the final result first, taken as a \textit{definition} in the formal setting.
\begin{defin}\label{exppp} Let $k\in \corps{Z}$.
\beqns
  i_{z,w} \, (z-w)^k &=& \sum_{j=0}^\infty \binom{k}{j} (-w)^j z^{k-j} \quad \in \quad \corps{C}[z][[z^{-1}]][[w]] \\
  i_{w,z} \, (z-w)^k &=& \sum_{j=0}^\infty \binom{k}{j} z^j (-w)^{k-j} \quad \in \quad \corps{C}[w][[w^{-1}]][[z]] \ , 
\eeqns
where the extended binomial coefficient\index{extended binomial coefficient} is defined below.
\end{defin}
\begin{defin}
Let $j\in\corps{Z}$ and $n\in\corps{Z}^+$. The extended binomial coefficient is defined by
\begin{eqnarray}
\binom{j}{n} = 
\begin{cases} \prod_{k=1}^n \frac{j-k+1}{k} = \frac{j(j-1)\cdots(j-n+1)}{n!} &  \text{if $n\neq 0$}   \\  
1 & \text{if $n=0$}  
\end{cases} 
\label{bin}
\end{eqnarray}
\end{defin}
The definition \ref{exppp} relies on the following proposition.

\begin{prop}\label{propexp}
Let $z$ and $w \, \in \corps{C}$.  Then, the formulas in definition \ref{exppp} are the series expansions of meromorphic functions on $\corps{C}$, converging in their respective domain.  They can naturally be generalized formally.
\end{prop}

\begin{proof}
Let $z$ and $w \, \in \corps{C}$ and $k<0$. We set $x=w/z$.\\

\noindent\underline{First case} : $|z|>|w|$ \\

We will first show that
$$i_{z,w} \, (z-w)^k = z^{-1} w^{k+1}  \sum_{j=0}^\infty \binom{j}{-1-k}  \left(\frac{w}{z}\right)^j \ .$$
The latter expression will then be proved to be equivalent to the expected result.

Let $|z|>|w|$, so that $|x|<1$.
	\begin{eqnarray}
	(z-w)^k &=&  z^k (1-x)^k \\
														 &=& z^k \left(\frac{1}{1-x}\right)^{-k} \\
														 &=& z^k (1+x+x^2 + \ldots)^{-k} \,\,\, \text{for } |x|<1 \label{temp}
	\end{eqnarray}
	But
	\begin{eqnarray}
	(1+x+x^2 + \ldots)^{-k} = \sum_{j=0}^\infty \binom{j}{-1-k}x^{j+k+1} \label{sum}
	\end{eqnarray}
	Indeed, we have 
	\begin{eqnarray}
	 \frac{1}{1-x} &=& \sum_{j=0}^\infty x^j \,\,\, \text{for } |x|<1 \label{de}
	\end{eqnarray}
	Thus, after successive derivations with respect to $x$ and some permitted relabeling of indices,
	\begin{eqnarray} 
	 \frac{1}{(1-x)^{-k}} 	&=& \sum_{j=0}^\infty \frac{j(j-1)\ldots(j-(-k-1)+1)}{(-k-1)!} \, x^{j+k+1} \\
                          &=& \sum_{j=0}^\infty
\frac{j!}{(-k-1)!\,(j-(-k-1))!} \, x^{j+k+1} \\
                          &=& \sum_{j=0}^\infty \label{fin}
\binom{j}{-k-1} \, x^{j+k+1} 
	                         \ ,
	\end{eqnarray}
which proves (\ref{sum}). Here, the binomial coefficient is still the standard one (recall $k<0$).
Inserting (\ref{sum}) in (\ref{temp}), we get the final result by the following computation :
\beqns
  i_{z,w} (z-w)^k &=& \sum_{j=0}^\infty \binom{j}{-1-k} z^{-1-j} w^{k+j+1} \\
                  &=& \sum_{j=0}^\infty \binom{j-1-k}{-1-k} z^{-j+k} w^{j} \\
                  &=& \sum_{j=0}^\infty (-1)^{j} \binom{j-1-k}{-1-k} z^{-j+k} (-w)^{j}
\eeqns
Now, making use of the extended binomial coefficient,
\begin{eqnarray*}
(-1)^{j} \binom{j-1-k}{-k-1} 
&=& (-1)^{j} \frac{(j-k-1)!}{(-k-1)! \, j!} \\
&=& (-1)^{j} \frac{(-k)(-k+1)\cdots(j-k-1)}{j!} \\
&=& \frac{k(k-1)\cdots (k-j+1)}{j!} \,\doteq\, \binom{k}{j}
\end{eqnarray*}

Hence
\beqns
  i_{z,w} (z-w)^k &=& \sum_{j=0}^\infty (-1)^{j} \binom{j-1-k}{-1-k} z^{-j+k} (-w)^{j} \\
                  &=& \sum_{j=0}^\infty \binom{k}{j} (-w)^{j} z^{-j+k} 
\eeqns

\noindent\underline{Second case} : $|z|<|w|$\\

We will now show that
$$i_{w,z} \, (z-w)^k = (-1)^k z^{k+1} w^{-1} \sum_{j=0}^\infty \binom{j}{-1-k} \left(\frac{z}{w}\right)^j \ .$$
Let $|z|<|w|$.  As $(z-w)^k=(-1)^k (w-z)^k$, we get the previous case, with $z$ and $w$ interchanged and the additional factor of $(-1)^k$. The final result is obtained in a similar way to the first case considered.

Those expansions are then proved for $z,w \in \corps{C}$.  But they make also sense formally, because they rely on well defined algebraic operations.
\end{proof}

\begin{remarque}
Notice in the definition \ref{exppp} that, for $k>0$, the sums are finite and the expansions are equal, as expected.  
\end{remarque}

\begin{remarque}
In other references, such as \cite{lepow}, the maps $i_{z,w}$ and $i_{w,z}$ are omitted but implicitly understood in the chosen notation.  For example, when considered as expansions, $(z-w)^k$ and $(-w+z)^k$ then differ in \cite{lepow}, as they respectively mean $i_{z,w} (z-w)^k$ and $i_{w,z} (z-w)^k$ here (as in \cite{kac1}). 
\end{remarque}

Finally, the maps $i_{z,w}$ and $i_{w,z}$ commute with the derivation operators $\partial_z$ and $\partial_w$. For further applications in computations, it is convenient to prove it at this stage.
\begin{prop}
The maps $i_{z,w}$ and $i_{w,z}$ commute with the derivation operators $\partial_z$ and $\partial_w$ on $\corps{C}[z^{\pm 1},w^{\pm 1},(z-w)^{-1}]$.
\end{prop}
\begin{proof}
We will prove that $[i_{z,w},\partial_z]=0$.  By linearity, it suffices to prove it on elements of the type $(z-w)^{-j-1}$, where $j\geq0$ (the case of monomials of the form $z^m$ and $w^m$ being obvious).
We have
\begin{eqnarray}
\partial_z i_{z,w}(z-w)^{-j-1} &=& \partial_z \sum_{m=0}^\infty \binom{m}{j} z^{-m-1} w^{m-j} \label{e1} \\
&=&  - \sum_{m=0}^\infty \binom{m}{j} (m+1) z^{-m-2} w^{m-j} \label{e2} \\
&=&  - \sum_{m=0}^\infty \frac{(m+1)!}{(j+1)!(m+1-j-1)!} (j+1) z^{-(m+1)-1} w^{m+1-j-1} \label{e3} \\
&=&  - (j+1) \sum_{m=1}^\infty \binom{m}{j+1} z^{-m-1} w^{m-j-1}  \label{e4} \\ 
&=&  - (j+1) i_{z,w}(z-w)^{-j-2}  \label{e5} \\
&=&   i_{z,w}\partial_z (z-w)^{-j-1} \label{e6}\ , 
\end{eqnarray}
The proof for $\partial_w$ is analogous. In a similar way, we prove it for $i_{w,z}$.
\end{proof}

\section{Dirac's $\delta$ distribution}
Among all the formal distributions in two variables, the formal \textbf{Dirac's $\delta$ distribution}\index{Dirac's $\delta$ distribution}, denoted by $\delta(z,w)$, is particularly important.  It can be introduced by noticing that, as we will prove it, $\delta(z,w)$ is the unique formal distribution satisfying the relation $\Res{z} a(z)\delta(z,w)=a(w)$, for any $a(z)\in {\cal A}[[z,z^{-1}]]$, where ${\cal A}$ is any algebraic structure over $\corps{C}$. The latter equation stresses the analogy with the usual Dirac's $\delta$ distribution (in Schwarz's theory).
\begin{defin}
The formal Dirac's $\delta$ distribution $\delta(z,w)$ is defined by
\begin{eqnarray}
\delta(z,w)=\sum_{n\in \corps{Z}} z^{-1} \left( \frac{w}{z} \right)^n = \sum_{n\in \corps{Z}} w^{-1} \left( \frac{z}{w} \right)^n\ \, \in \, \corps{C}[[z^{\pm 1},w^{\pm 1}]] \ . \label{deltadef2}
\end{eqnarray}
\end{defin}
\begin{remarque}
The Fourier expansion of $\delta(z,w)$ can easily be shown to be  
\begin{eqnarray}
\delta(z,w)=\sum_{n,m \in \corps{Z}} \delta_{n,-1-m} \, z^n w^m \, \in \, \corps{C}[[z^{\pm 1},w^{\pm 1}]] \ . \label{deltadef1}
\end{eqnarray}
\end{remarque}

Using the notion of expansion, we can give a new insight to the Dirac's $\delta$ distribution.
\begin{defin} An equivalent definition  of the Dirac's $\delta$ distribution $\delta(z,w)$ is given by
\begin{eqnarray}
\delta(z,w)=i_{z,w} \frac{1}{z-w} - i_{w,z} \frac{1}{z-w} \ ,\label{deltadef3}
\end{eqnarray}
where both members have to be considered as elements of $\corps{C}[[z,z^{-1},w,w^{-1}]]$. 
\end{defin}
This permits us to interpret $\delta(z,w)$ as the difference of the expansions of the same rational function, computed in two different domains.  We stress that this interpretation, as the difference in \eqref{deltadef3}, is only valid in the space $\corps{C}[[z,z^{-1},w,w^{-1}]]$.  This property, among others to be discussed below, happens to be shared by distributions in the sense of Schwartz's theory.
\begin{prop}\label{equivvv}
The two definitions of $\delta(z,w)$ given above are indeed equivalent.
\end{prop}
\begin{proof}
The definition given above rely on the following computation.
\begin{eqnarray*}
\delta(z,w)&=& z^{-1} \sum_{n\in \corps{Z}}  \left( \frac{w}{z} \right)^n \\
&=& z^{-1} \sum_{n=0}^\infty  \left( \frac{w}{z} \right)^n + z^{-1} \sum_{n=-1}^{-\infty}  \left( \frac{w}{z} \right)^n \\
&=& z^{-1} \sum_{n=0}^\infty  \left( \frac{w}{z} \right)^n + z^{-1} \sum_{n=1}^{\infty}  \left( \frac{z}{w} \right)^n \\
&=& z^{-1} \sum_{n=0}^\infty  \left( \frac{w}{z} \right)^n + z^{-1} \sum_{n=0}^{\infty} \left( \frac{z}{w} \right)^{n+1} \\
&=& z^{-1} \sum_{n=0}^\infty  \left( \frac{w}{z} \right)^n + w^{-1} \sum_{n=0}^{\infty}  \left( \frac{z}{w} \right)^n \ .
\end{eqnarray*}
We notice that 
\begin{eqnarray*}
z^{-1} \sum_{n=0}^\infty   \left( \frac{w}{z} \right)^n &=& i_{z,w}\frac{1}{z-w}  \in \,\,  \corps{C}[[z^{-1}]][[w]] \,\,  \subset \,\,  \corps{C}[[z^{\pm 1},w^{\pm 1}]] \\
- w^{-1} \sum_{n=0}^{\infty}  \left( \frac{z}{w} \right)^n &=& i_{w,z}\frac{1}{z-w} \in \,\,  \corps{C}[[w^{-1}]][[z]] \,\,   \subset \,\,  \corps{C}[[z^{\pm 1},w^{\pm 1}]] \ ,
\end{eqnarray*}
from which the result is clearly deduced.
\end{proof}

The Dirac's $\delta(z,w)$ is obviously derivable and we enumerate its main properties\index{Dirac's $\delta$ distribution!properties}.
\begin{prop}\label{deltaprop} Let $n,m\in\corps{N}$.  The formal distribution $\delta(z,w)$ satisfies :
\begin{enumerate}
\renewcommand{\labelenumi}{\emph{(\theenumi)}}
\item \label{p1} $(z-w)^m \partial_w^n \delta(z,w) = 0 \quad \forall m > n$. 
\item \label{p2} $(z-w)\frac{1}{n!}\partial_w^n \delta(z,w) = \frac{1}{(n-1)!}\partial_w^{n-1} \delta(z,w)$ if $n\geq 1$. 
\item \label{p3} $\delta(z,w)=\delta(w,z)$. 
\item \label{p4} $\partial_z \delta(z,w)=-\partial_w \delta(w,z)$.
\item \label{p5} $a(z)\delta(z,w)=a(w)\delta(z,w)$ where $a(z)$ is any formal distribution. 
\item \label{p6} $\text{Res}_z a(z)\delta(z,w)=a(w)$, which implies that $\delta(z,w)$ is unique. 
\item \label{p7} $\exp(\lambda(z-w))\partial_w^n \delta(z,w)=(\lambda + \delta_w)^n \delta(z,w)$. 
\end{enumerate}
\end{prop}
\begin{proof}

Let's demonstrate (\ref{p1}).  Clearly, $[z^n,i_{z,w}]=0$ and $[w^n,i_{z,w}]=0$, so that by linearity :
$$\lbrack \sum_{j=0}^m \binom{m}{j} z^{m-j} w^j\, , \, i_{z,w}\rbrack = 0 \ .$$
Therefore $\lbrack (z-w)^m\, , \, i_{z,w}\rbrack = 0$,
and by symmetry, $\lbrack (z-w)^m\, , \, i_{w,z}\rbrack = 0$.
Combining those results to the definition, we obtain
\begin{eqnarray*}
\frac{1}{n!}(z-w)^m  \partial_w^n \delta(z,w)  = \,\, i_{z,w} \frac{(z-w)^m}{(z-w)^{n+1}} - i_{w,z} \frac{(z-w)^m}{(z-w)^{n+1}} \ .
\end{eqnarray*}
But, $\lbrace i_{z,w} - i_{w,z} \rbrace (z-w)^{m-n-1}$ is zero, since $m>n$. This proves (\ref{p1}).

To prove (\ref{p2}), using the appropriate definition, we have
\begin{eqnarray*}
(z-w) \frac{1}{n!}  \partial_w^n \delta(z,w) &=& (z-w) \sum_{j\in \corps{Z}} \binom{j}{n} z^{-j-1} w^{j-n} \\
&=& \sum_{j\in \corps{Z}} \binom{j}{n} z^{-j} w^{j-n} - \sum_{j\in \corps{Z}} \binom{j}{n} z^{-j-1} w^{j-n+1} \\
&=& \sum_{\tilde{j}\in \corps{Z}} \binom{\tilde{j}+1}{n} z^{-\tilde{j}-1} w^{\tilde{j}-n+1} - \sum_{j\in \corps{Z}} \binom{j}{n} z^{-j-1} w^{j-n+1} \ ,
\end{eqnarray*}
where $j=\tilde{j}+1$.
It is easy to show that
$\binom{j+1}{n}-\binom{j}{n} = \binom{j}{n-1}$, 
which permits us to combine the two sums in the preceding equality and obtain
\begin{eqnarray*}
(z-w) \frac{1}{n!} \partial_w^n \delta(z,w) = \sum_{j\in \corps{Z}} \binom{j}{n-1} z^{-j-1} w^{j-n+1} &=& \frac{1}{(n-1)!} \partial_w^{n-1} \delta(z,w) \ ,
\end{eqnarray*}
by definition.  This proves (\ref{p2}).  

The property (\ref{p3}) : straight from the definition.

The following computation 
\begin{eqnarray*}
\partial_z \delta(z,w) &=& \partial_z \sum_{n\in\corps{Z}} w^n z^{-n-1} \,=\, \sum_{n\in\corps{Z}} w^n (-1) (n+1) z^{-n-2} \\
&=& - \sum_{n\in\corps{Z}} (n+1) w^n z^{-n-2} \,=\, - \sum_{m\in\corps{Z}} m\,w^{m-1} z^{-m-1} \quad \text{where } m=n+1 \\
&=& - \partial_w \sum_{m\in\corps{Z}} w^{m} z^{-m-1} \,=\, - \partial_w \delta(z,w) \\
&=& - \partial_w \delta(w,z) \quad \text{by (\ref{p3})}
\end{eqnarray*}
prove (\ref{p4}).

For (\ref{p5}), we start from (\ref{p1}), according to which $(z-w)\delta(z,w)=0$, so that $z^n \delta(z,w) = w^n \delta(z,w)$.  By linearity, we get (\ref{p5}). 

Applying the residue to (\ref{p5}), we obtain
$\text{Res}_z a(z) \delta(z,w) = a(w) \text{Res}_z \delta(z,w) = a(w)$,  
proving (\ref{p6}). The unicity is proved by contradiction. By linearity, it suffices to consider $a(z)=z^k$. Suppose there exist two formal distributions $\delta=\sum \delta_{n,m} z^n w^m$ and $\delta'= \sum \delta'_{n,m} z^n w^m$, with $\delta_{n,m}\neq\delta'_{n,m}\,\,\forall k,m\in\corps{Z}$, such that 
$\text{Res}_z \delta(z,w) z^k=w^k=\text{Res}_z \delta'(z,w) z^k$.
Those equalities imply
$$ \text{Res}_z \delta(z,w) z^k=\sum \delta_{-1-k,m}w^m=\text{Res}_z \delta'(z,w) z^k=\sum \delta'_{-1-k,m}w^m \ ,$$
so $\delta_{-1-k,m}=\delta'_{-1-k,m}\,\,\forall k,m\in\corps{Z}$, and hence the contradiction. Therefore, (\ref{p6}) imply the unicity of $\delta(z,w)$.

At last (\ref{p7}), for which we expand the exponential
\begin{eqnarray*}
e^{\lambda (z-w)} \partial_w^n \delta(z,w) = \sum_{k=0}^\infty \frac{\lambda^k}{k!}n!(z-w)^k \frac{\partial_w^n \delta(z,w)}{n!} \ ,
\end{eqnarray*}
in which
\begin{eqnarray*}
(z-w)^k \frac{\partial_w^n \delta(z,w)}{n!} &=& (z-w)^{k-1} \frac{(z-w)\partial_w^n \delta(z,w)}{n!}\\
&=& (z-w)^{k-1} \frac{\partial_w^{n-1} \delta(z,w)}{(n-1)!} \quad \text{by (\ref{p2})} \\
\text{($k$ times)}&\vdots&\\
&=&\frac{\partial_w^{n-k} \delta(z,w)}{(n-k)!}\ ,
\end{eqnarray*}
where $k\leq n$.
So, 
\begin{eqnarray*}
e^{\lambda (z-w)} \partial_w^n \delta(z,w) &=& \sum_{k=0}^n \frac{\lambda^k n!}{k!} 
\frac{\partial_w^{n-k} \delta(z,w)}{(n-k)!} \\
&=& \left( \sum_{k=0}^n \binom{n}{k} \lambda^k \partial_w^{n-k} \right) \delta(z,w) \\
&=& (\lambda + \partial_w)^n \delta(z,w) \ ,
\end{eqnarray*}
which proves (\ref{p7}).

\end{proof}

\section{Locality}
Operator product expansions are widely used in CFT.  The fundamental notion beneath it is the locality\index{locality} of formal distributions in two variables. We just happened to meet our first local formal distribution, namely the formal Dirac's $\delta$ distribution.
\begin{defin}
Let $a(z,w)\in {\cal A}[[z^{\pm 1},w^{\pm 1}]]$. 
The formal distribution $a(z,w)$ is called \textbf{local} if there exists $n >> 0$ in $\corps{N}$ such that $(z-w)^n a(z,w)=0$.
\end{defin}
\begin{exe}
The Dirac's $\delta(z,w)$ and its derivatives are local.  In general, if $a(z,w)$ is local, then its derivatives (with respect to $z$ or $w$) are local as well, as is $a(w,z)$.
\end{exe}
\begin{proof}
The property \ref{deltaprop}(\ref{p1}) proves the first statement.
For the second one, clearly, if $(z-w)^n a(z,w)=0$, then $(w-z)^n a(w,z) = 0$.
Let $\corps{N} \ni n>>0$ such that $(z-w)^{n-1} a(z,w)=0$.  This implies that $(z-w)^n a(z,w)=0$ trivially. Now, the following computation
\begin{eqnarray*}
(z-w)^n \partial_z a(z,w) &=& \partial_z (z-w)^n a(z,w) - a(z,w) \partial_z (z-w)^n  \\
&=& -n a(z,w) (z-w)^{n-1} \\
&=& 0 \ ,
\end{eqnarray*}
proves that $\partial_z a(z,w)$ is a local distribution.  By symmetry, $\partial_w a(z,w)$ is local too.
\end{proof}

Locality is a strong assumption.  By the following theorem\index{locality!decomposition theorem}, any local non vanishing formal distribution $a(z,w)$ is necessarily singular.  Moreover, the form of this singularity can be determined.
\begin{theo}\label{decomp}
Let $a(z,w)$ be a local $\cal A$-valued formal distribution.  Then $a(z,w)$ can be written as a finite sum of $\delta(z,w)$ and its derivatives :
\begin{eqnarray}
a(z,w)=\sum_{j\in \corps{Z}^+} c^j (w) \frac{\partial_w^j \delta(z,w)}{j!} \label{loc}
\end{eqnarray}
where $c^j(w) \,\in\, {\cal A}[[w,w^{-1}]]$ are formal distributions given by
\begin{eqnarray}
c^j (w)= \text{Res}_z (z-w)^j a(z,w) \ .
\end{eqnarray}
In addition, the converse is true.
\end{theo}
\begin{proof}

Let an integer $N>0$.  We have to prove the following equivalence.
\begin{eqnarray}
(z-w)^N a(z,w) = 0 \,\Leftrightarrow \, a(z,w)=\sum_{j=0}^{N-1} c^j (w) \frac{\partial_w^j \delta(z,w)}{j!} \ ,\label{exuni}
\end{eqnarray}
and the unicity of such a decomposition. 

\noindent $\Leftarrow$ : this is a direct consequence of the first property of the proposition \ref{deltaprop}.\\
\noindent $\Rightarrow$ : this implication is less trivial. Let the formal distribution
$$a(z,w)=\sum_{n,m\in\corps{Z}} a_{n,m} z^{-1-n} w^{-1-m} \ .$$
For $N=1$, we have
\begin{eqnarray}
(z-w) a(z,w) = 0 \, \Rightarrow \, a_{n+1,m}=a_{n,m+1} \, \Rightarrow \, a(z,w)=c(w)\delta(z,w) \label{a01}
\end{eqnarray}
with $c(w)=\text{Res}_z a(z,w)$.

Indeed, the first implication is trivial. For the second one, we note that by successive shifts $n\to n-1$ and $m\to m+1$, the recurrence relation can also be written $a_{n,m}=a_{0,n+m}$. Therefore,
\begin{eqnarray*}
a(z,w) &=& \sum_{n,m\in\corps{Z}} a_{n,m} z^{-1-n} w^{-1-m} \\
&=& \sum_{n,m\in\corps{Z}} a_{0,n+m} z^{-1-n} w^{-1-m} \\
&=& \sum_{k,n\in\corps{Z}} a_{0,k} z^{-1-n} w^{-1-k+n} \\
&=& \sum_{k\in\corps{Z}} a_{0,k} w^{-1-k} \sum_{n\in\corps{Z}} z^{-1-n} w^{n} \\
&=& \text{Res}_z a(z,w) \cdot \delta(z,w) \ ,
\end{eqnarray*}
which proves the implication for $N=1$.

The general case can be deduced from the latter.  Indeed,
\begin{eqnarray*}
(z-w)^N a(z,w) = 0 \Leftrightarrow (z-w) \lbrack (z-w)^{N-1} a(z,w) \rbrack = 0 \ ,
\end{eqnarray*}
thus, by (\ref{a01}),  
\begin{eqnarray*}
(z-w)^{N-1} a(z,w) = c^{N-1} (w) \delta(z,w) \ , 
\end{eqnarray*}
with $c^{N-1} (w) = \text{Res}_z (z-w)^{N-1} a(z,w)$.
Let's use the second property of the proposition \ref{deltaprop}, as many times as necessary in order to make appear the factor $(z-w)^{N-1}$ in the right-hand side, before putting everything in the left-hand side to obtain
\begin{eqnarray*}
(z-w)^{N-1}\, \lbrack a(z,w) -  c^{N-1} (w) \frac{1}{(N-1)!} \partial_w^{N-1} \delta(z,w) \rbrack = 0 \ . 
\end{eqnarray*}
Using (\ref{a01}) again, we get
\begin{eqnarray*}
(z-w)^{N-2}\, \lbrack a(z,w) -  c^{N-1} (w) \frac{1}{(N-1)!} \partial_w^{N-1} \delta(z,w) \rbrack = c^{N-2} (w) \delta(z,w) \ , 
\end{eqnarray*}
where
\begin{eqnarray*}
c^{N-2} (w) &=& \text{Res}_z (z-w)^{N-2} \lbrack a(z,w) -  c^{N-1} (w) \frac{1}{(N-1)!} \partial_w^{N-1} \delta(z,w) \rbrack \\
&=& \text{Res}_z (z-w)^{N-2} a(z,w) \ , 
\end{eqnarray*}
for, using $(N-2)$ times the second property of proposition \ref{deltaprop}, 
 
{\setlength\arraycolsep{0pt}
\begin{eqnarray*}
\text{Res}_z c^{N-1}(w) && (z-w)^{N-2} \frac{1}{(N-1)!} \partial_w^{N-1} \delta(z,w) \\ 
&&= \,\, c^{N-1} (w) \,  \text{Res}_z (z-w)^{N-3} \frac{1}{(N-1)!} (z-w) \partial_w^{N-1} \delta(z,w) \\
&&= \,\,  c^{N-1} (w) \,  \text{Res}_z (z-w)^{N-3} \frac{1}{(N-2)!} \partial_w^{N-2} \delta(z,w) \\
&&\, \vdots \\
&&= \,\, c^{N-1} (w) \,  \text{Res}_z  \frac{1}{(N-N)!} \partial_w^{N-N+1} \delta(z,w) \\
&&= \,\, c^{N-1} (w) \,  \partial_w \text{Res}_z \delta(z,w) \\
&&= \,\,  0 \ ,
\end{eqnarray*}}
because $\text{Res}_z \delta(z,w)=1$.

As above, we use once again the second property of proposition \ref{deltaprop}, as many times as necessary to make appear the factor $(z-w)^{N-2}$ at the RHS, before putting all at the LHS to obtain

{\setlength\arraycolsep{0pt}
\begin{eqnarray*}
(z-w)^{N-2}\,\, \lbrack\,\, && a(z,w) \\
&&\,\, -  c^{N-1} (w) \frac{1}{(N-1)!} \partial_w^{N-1} \delta(z,w) \\
&&\,\, - c^{N-2}(w) \frac{1}{(N-2)!}\partial_w^{N-2} \delta(z,w) \,\,\rbrack = 0 \ . 
\end{eqnarray*}}
We repeat this procedure till we get
{\setlength\arraycolsep{0pt}
\begin{eqnarray*}
(z-w)\,\, \lbrack\,\, && a(z,w) \\
&&\,\, - \,\, c^{N-1} (w) \frac{1}{(N-1)!} \partial_w^{N-1} \delta(z,w) \\
&&\,\, - \,\,\cdots\,\, -\,\, c^1 (w) \partial_w \delta(z,w) \,\,\rbrack = 0 \ . 
\end{eqnarray*}}
At last, using (\ref{a01}) one last time, we get
{\setlength\arraycolsep{0pt}
\begin{eqnarray*}
a(z,w) = c^0 (w)&& \delta(z,w) + c^1 (w) \partial_w \delta(z,w) \\
&&\,\,+ \cdots +  c^{N-1} (w) \frac{1}{(N-1)!} \partial_w^{N-1} \delta(z,w) \ , 
\end{eqnarray*}}
with $c^j (w) = \text{Res}_z (z-w)^j a(z,w)$.

In addition, this decomposition is clearly unique, which concludes the proof.
\end{proof}
\begin{remarque}
Notice that the condition 
\begin{eqnarray} \text{Res}_z a(z,w) (z-w)^j = 0 \text{ for an integer } j\geq N \ \label{weak}  \end{eqnarray}
is weaker than locality.  Indeed, we can show that any $a(z,w)$ satisfying this condition can be written as
$$ a(z,w)= \sum_{j=0}^{N-1} c^j(w) \frac{\partial_w^j \delta(z,w)}{j!} + b(z,w) \ ,$$
with $c^j (w) = \text{Res}_z a(z,w) (z-w)^j$ and $b(z,w)$ holomorphic in $z$, \textit{i.e.} $b(z,w) \in {\cal A}[[z,w^{\pm 1}]]$.  For this reason, $a(z,w)$ is said \textbf{weakly local}\index{locality!weak locality}.  It is precisely since $b(z,w)$ is holomorphic that it vanishes when $a(z,w)$ is local.
\end{remarque}
\begin{remarque}
As will be shown later, \eqref{loc} is in fact a disguised OPE.
\end{remarque}

\section{Fourier transform}\label{fourtrrr}
We now define the notion of \textbf{Fourier transform}\index{Fourier transform} in two cases : in one and two indeterminates.  The first case will be used to show that every vertex algebra is a Lie conformal algebra.  The second one will allow us to precisely define the latter algebraic structure.
\begin{defin}\label{transfour1ind}
Let $a(z) \in {\cal A}[[z^{\pm 1}]]$.  We define the Fourier transform in one indeterminate\index{Fourier transform!in one indeterminate} of $a(z)$ by
$\Fzz a(z)=\text{Res}_z e^{\lambda z}a(z)$.
Consequently, $\Fzz$ is a linear map from ${\cal A}[[z^{\pm 1}]]$ to ${\cal A}[[\lambda]]$.
\end{defin}
\begin{prop}\label{trfour1d} $\Fzz$ satisfies the following properties :
\begin{enumerate}
\renewcommand{\labelenumi}{\emph{(\theenumi)}}
\item $\Fzz \partial_z a(z) = - \lambda \Fzz a(z)$;
\item $\Fzz \left( e^{z T} a(z) \right) = \Fzd{z}{z+T} a(z)$, where $T\in End{\cal U}$ and $a(z)\in {\cal U}((z))$;
\item $\Fzz a(-z) = - \Fzd{z}{-\lambda} a(z)$;
\item $\Fzz \partial_w^n \delta(z,w) = e^{\lambda w} \lambda^n$.
\end{enumerate}
\end{prop}
\begin{proof}
(1) :
\begin{eqnarray*}
\partial_z e^{\lambda z} a(z) = \lambda e^{\lambda z} a(z) + e^{\lambda z} \partial_z a(z) \ .
\end{eqnarray*}
Now, $\Res{z}\partial_z (\cdot) = 0$, so applying the residue to both members :
$$0 = \lambda \Fzz a(z) + \Fzz \partial_z a(z) \ ,$$
which proves (1).\\
(2) :
\begin{eqnarray*}
\Fzz \left( e^{z T}a(z) \right) &=& \Res{z} e^{\lambda z} \left( e^{z T} a(z) \right) \\
&=& \Res{z} \left( e^{(\lambda + T)z} a(z) \right) \\
&=& \Fzd{z}{\lambda + T} a(z)
\end{eqnarray*}
To prove (3), we define the following formal distribution :
$$  b(z)=e^{\lambda z} a(-z) + e^{-\lambda z} a(z) \ ,$$
which is even : $b(z)= b(-z)$.  So all even powers in $z$ must vanish and in particular the power $-1$, i.e. $\Res{z} b(z)=0$.  Now this equality can be written
$\Res{z} e^{\lambda z} a(-z) + \Res{z} e^{-\lambda z} a(z) = 0$.
So that $\Fzz a(-z) = - \Fzd{z}{-\lambda} a(z)$.\\
(4) :
\begin{eqnarray*}
\Fzz \partial_w^n \delta(z,w) &=& \Res{z} e^{\lambda z} \partial_w^n \sum_{m\in \corps{Z}} w^m z^{-1-m} \\
&=& \sum_{m\in \corps{Z}}\sum_{j\in \corps{Z}^+} \partial_w^n w^m \frac{\lambda^j}{j!} \Res{z} z^{j-m-1} \\
&=& \partial_w^n \sum_{j\in \corps{Z}^+} \frac{(\lambda w)^j}{j!} \\
&=& \partial_w^n e^{\lambda w} \,=\, \lambda^n e^{\lambda w} \ .
\end{eqnarray*}

\end{proof}
\begin{defin}\label{transfour}
Let $a(z,w) \in {\cal A}[[z^{\pm 1}, w^{\pm 1}]]$.  We define the Fourier transform in two indeterminates\index{Fourier transform!in two indeterminates} of $a(z,w)$ by
$$\Fz a(z,w)=\text{Res}_z e^{\lambda(z-w)}a(z,w) \ .$$
\end{defin}
Consequently, $\Fz$ is a linear map from ${\cal A}[[z^{\pm 1}, w^{\pm 1}]]$ to ${\cal A}[[w^{\pm 1}]][[\lambda]]$.
 Indeed it can be easily shown that, for $a(z,w) \in {\cal A}[[z^{\pm 1}, w^{\pm 1}]]$, we have 
 $$\Fz a(z,w) = \sum_{j\in \corps{Z}^+} \frac{\lambda^j}{j!}c^j(w) \ ,$$ 
 where $c^j(w) = \text{Res}_z (z-w)^j a(z,w)$.
\begin{prop}\label{Four}
Let $a(z,w)\in {\cal A}[[z^{\pm 1}, w^{\pm 1}]]$. Then,
\begin{enumerate}
\renewcommand{\labelenumi}{\emph{(\theenumi)}}
\item if $a(z,w)$ is local, then $\Fz a(z,w)\in{\cal A}[[w^{\pm 1}]][\lambda]$; \label{p001}
\item \label{p002} $\Fz \partial_z a(z,w) = -\lambda \Fz a(z,w)=\lbrack \partial_w, \Fz \rbrack$; 
\item \label{p003} if $a(z,w)$ is local, then $\Fz a(w,z) = \F{z}{w}{-\lambda-\partial_w} a(z,w)$, where $\F{z}{w}{-\lambda-\partial_w} a(z,w) \doteq \F{z}{w}{\mu} a(z,w)|_{\mu=-\lambda-\partial_w}$;
\item \label{p004} $\Fz \F{z}{w}{\mu} \, : \, {\cal U}[[z^{\pm 1},w^{\pm 1},x^{\pm 1}]] \to {\cal U}[[w^{\pm 1}]][[\lambda,\mu]]$
and $\Fz \F{x}{w}{\mu} a(z,w,x)= \F{x}{w}{\lambda+\mu}\F{z}{x}{\lambda} a(z,w,x)$.
\end{enumerate}
\end{prop}

\begin{proof}
(\ref{p001}) is clearly satisfied if $a(z,w)$ is local.

The first equality of (\ref{p002}) results from the following computation.
\begin{eqnarray*} 
\Fz \partial_z a(z,w) &=& \Res{z} e^{\lambda(z-w)} \partial_z a(z,w) \\
&=& -\Res{z} \partial_z \left( e^{\lambda(z-w)} \right) a(z,w) \quad \text{because $\Res{z} \partial_z (\cdot) = 0$} \\
&=& - \lambda\, \Res{z} e^{\lambda(z-w)} a(z,w) \\
&=& - \lambda\, \Fz a(z,w)
\end{eqnarray*}
For the second equality of (\ref{p002}), we start from
\begin{eqnarray*}
\partial_w \Fz a(z,w) = \partial_w \Res{z} e^{\lambda(z-w)} a(z,w) = \left( - \lambda \Fz + \Fz \partial_w \right) a(z,w) \ ,
\end{eqnarray*}
which implies that 
$$- \lambda \Fz = \partial_w \Fz - \Fz \partial_w = \lbrack \partial_w , \Fz \rbrack \ , $$
which concludes the proof.

For (\ref{p003}), assuming $a(z,w)$ is local and $\Fz$ being linear, we can consider the case where 
$$a(z,w)=c(w) \partial_w^j \delta(z,w) \ ,$$ thanks to the decomposition theorem.  Then, 
\begin{eqnarray*}
\Fz a(w,z) &=& \Fz c(z) \partial_z^j \delta(w,z) \,=\, \Res{z} e^{\lambda(z-w)} c(z) (-\partial_w)^j \delta(z,w) \quad \text{by \ref{deltaprop}(\ref{p4})}\\
&=& \Res{z} c(z) (-\lambda-\partial_w)^j \delta(z,w) \quad \text{by \ref{deltaprop}(\ref{p7})} \\
&=& (-\lambda-\partial_w)^j \Res{z} c(z) \delta(z,w) \,=\, (-\lambda-\partial_w)^j c(w)\quad \text{by \ref{deltaprop}(\ref{p6})} \\
&=& \F{z}{w}{-\lambda-\partial_w} a(z,w)\ ,
\end{eqnarray*}
where, for the last equality,
{\setlength\arraycolsep{0pt}\begin{eqnarray*}
\F{z}{w}{-\lambda-\partial_w}&& a(z,w) := \F{z}{w}{\mu} a(z,w)|_{\mu=-\lambda-\partial_w} \\
&&=\,\, \Res{z} e^{\mu(z-w)} c(w) \partial_w^j \delta(z,w)|_{\mu=-\lambda-\partial_w} \\
&&=\,\, \Res{z} \sum_{k=0}^\infty \frac{\mu^k}{k!} c(w) (z-w)^k \partial_w^j \delta(z,w)|_{\mu=-\lambda-\partial_w} \\
&&=\,\, \sum_{k=0}^j \frac{\mu^k}{k!} c(w) j! \Res{z} (z-w)^k \frac{\partial_w^j \delta(z,w)}{j!}|_{\mu=-\lambda-\partial_w} \quad \text{by \ref{deltaprop}(\ref{p1})}\\
&&=\,\, \sum_{k=0}^j \frac{\mu^k}{k!} c(w) j! \Res{z} \frac{\partial_w^{j-k} \delta(z,w)}{(j-k)!}|_{\mu=-\lambda-\partial_w} \quad \text{by \ref{deltaprop}(\ref{p2})}\\
&&=\,\, \sum_{k=0}^j \frac{\mu^k}{k!} \frac{j!}{(j-k)!} c(w) \partial_w^{j-k}  \Res{z} \delta(z,w)|_{\mu=-\lambda-\partial_w} \\
&&=\,\, \sum_{k=0}^j \frac{\mu^k}{k!} \frac{j!}{(j-k)!} c(w) \delta_{k,j}|_{\mu=-\lambda-\partial_w} \\
&&=\,\, \mu^j c(w)|_{\mu=-\lambda-\partial_w} \,=\, (-\lambda-\partial_w)^j c(w) \ ,
\end{eqnarray*}}
which concludes the proof.

At last, (\ref{p004}) results from the following computation.
\begin{eqnarray*}
\Res{z} e^{\lambda(z-w)} \Res{x} e^{\mu(x-w)} a(z,w,x) = \Res{x} e^{(\lambda+\mu)(x-w)} \Res{z} e^{\lambda(z-x)} a(z,w,x) \ ,
\end{eqnarray*}
which is true, as
\begin{eqnarray*}
e^{\lambda(z-w)+\mu(x-w)}= e^{(\lambda+\mu)(x-w) + \lambda (z-x)} \ .
\end{eqnarray*}
\end{proof}
%

\chapter{Lie conformal algebras}
\section{Lie superalgebra}
In this section, we briefly recall main definitions and needed properties.  See \cite{kac2} for a detailed exposition.

A Lie superalgebra\index{Lie superalgebra} can be seen as a generalization of the notion of \textit{Lie algebra}, made compatible with \textit{supersymmetry}.
This means that we first introduce a decomposition by $\corps{Z}_2$ on a linear space $\mathfrak{g}$, i.e. $\mathfrak{g}=\mathfrak{g}_{\bar{0}} \oplus \mathfrak{g}_{\bar{1}}$, which is then said \textit{$\corps{Z}_2$-graded}\index{graded algebra}.  We then make the antisymmetry property and the Jacobi identity (from a standard Lie algebra) compatible with this gradation. More precisely, we define $p(a,b)\doteq(-1)^{p(a) p(b)}$, with $p(a),\, p(b) \in \corps{Z}_2=\{\bar{0},\bar{1}\}$, as the respective parities of homogeneous elements $a$ and $b$, and we pose $(-1)^{\bar{0}}\doteq1$ and $(-1)^{\bar{1}}\doteq-1$.  Now, we can state the following definition.
\begin{defin} \label{supalg}
A \textbf{Lie superalgebra} (or $\corps{Z}_2$-graded Lie algebra) $\mathfrak{g}$ on $\corps{C}$ is a $\corps{Z}_2$-graded complex linear space $$\mathfrak{g}=\mathfrak{g}_{\bar{0}} \oplus \mathfrak{g}_{\bar{1}} \ ,$$
endowed with the so-called \textit{Lie superbracket}\index{Lie superbracket}
$$[\;,\;]\; :\; \mathfrak{g} \otimes \mathfrak{g} \longrightarrow \mathfrak{g} \ , $$ 
satisfying the following properties :
\begin{enumerate}
\renewcommand{\labelenumi}{(\theenumi)}
\item $\lbrack \; , \; \rbrack$ is $\mathbb{C}$-bilinear;
\item $[\;,\; ]$ is compatible with the gradation, \textit{i.e.} 
$$[\mathfrak{g}_{z_k},\mathfrak{g}_{z_l}] \; \subset \; \mathfrak{g}_{z_k + z_l} \ ,$$ 
where $z_k,z_l \in \corps{Z}_2$;
\item $\lbrack \; , \; \rbrack$ is graded antisymmetric, i.e. 
$$\lbrack a , b \rbrack = - p(a,b) \lbrack b , a \rbrack \ ;$$
\item $\lbrack \; , \; \rbrack$ satisfies the graded Jacobi identity, 
\begin{eqnarray}
p(a,c) \; \lbrack a , \lbrack b, c \rbrack \rbrack \; + \;
p(b,a) \; \lbrack b , \lbrack c, a \rbrack \rbrack \; + \;
p(c,b) \; \lbrack c , \lbrack a, b \rbrack \rbrack \; = \; 0 \label{jaaacob}
\end{eqnarray}
\end{enumerate}
The set $\mathfrak{g}_{\bar{0}}$ (resp. $\mathfrak{g}_{\bar{1}}$) contains the so-called \textit{even} (resp. \textit{odd}) elements of $\mathfrak{g}$.
\end{defin}

\begin{remarque}
Note that the parity can be extended by linearity to non-homogeneous elements.
\end{remarque}

In an associative algebra, we can define the superbracket of homogeneous elements by $[a,b]=ab-p(a,b) ba$. It can then be extended by linearity to non-homogeneous elements, as the parity.
\begin{prop}\label{scp}
The map $p(\cdot, \cdot)$ and the superbracket defined above satisfy the following properties
\begin{enumerate}
\renewcommand{\labelenumi}{\emph{(\theenumi)}}
\item \label{scp1} $p(a,b)^2=1$, $p(a,[b,c])=p(a,b)p(a,c)$ and $p([a,b],c)=p(a,c)p(b,c)$
\item \label{scp2} $[a,b]=-p(a,b) [b,a]$
\item \label{scp3} $[a,[b,c]]=[[a,b],c] + p(a,b)[b,[a,c]]$
\end{enumerate}  
\end{prop}
\begin{proof}
The first equality of (\ref{scp1}) is clear.  For the second one, since the superbracket is compatible with the gradation, we have $[b,c]\in \mathfrak{g}_{p(b)+p(c)}$, so that $p([b,c])=p(b)+p(c)$. Then,
$$p(a,[b,c])=(-1)^{p(a)p([b,c])}=(-1)^{p(a)(p(b)+p(c))}=p(a,b)p(a,c) \ .$$  The third equality is analogous.
Let's add that we also have $p(a,bc)=p(a,b)p(a,c)$, since $p(bc)=p([b,c])$.

The following computation
\begin{eqnarray*}
p(a,b)[a,b] &=& p(a,b)ab - p(a,b)^2 ba \\
&=& -\left( ba - p(a,b)ab \right)\\
&=& -[b,a]
\end{eqnarray*}
proves (\ref{scp2}).
To check (\ref{scp3}), we make use of (\ref{scp1}) to obtain
\begin{eqnarray*}
[a,[b,c]]&=& a[b,c]-p(a,[b,c])[b,c]a \\
         &=& abc - p(b,c)acb - p(a,b)p(a,c)bca  + p(a,b)p(a,c)p(b,c)cba\\
\lbrack\lbrack a,b \rbrack ,c\rbrack &=& abc - p(a,b)bac - p(a,c)p(b,c)cab + p(a,c)p(b,c)p(a,b)cba\\
p(a,b)[b,[a,c]]&=& p(a,b) bac - p(a,b)p(a,c)bca - p(b,c) acb  + p(b,c)p(a,c)cab \ .
\end{eqnarray*}
By subtracting the two last equations from the first one, we prove (\ref{scp3}).
In addition, by (\ref{scp1}) and (\ref{scp2}), the property (\ref{scp3}) can be seen as equivalent to the Jacobi identity :
\begin{eqnarray*}
p(a,c)[a,[b,c]] &=& p(a,c) [[a,b],c] + p(a,c)p(a,b)[b,[a,c]] \quad \text{by (\ref{scp3})} \\
p(b,a) [b,[c,a]] &=& -p(b,a)p(a,c)[b,[a,c]] \quad \text{by (\ref{scp2})}\\
p(c,b)[c,[a,b]] &=& -p(c,b)p(c,[a,b])[[a,b],c] \quad \text{by (\ref{scp2})}\\
&=& -p(c,b)p(c,a)p(c,b)[[a,b],c] \quad \text{by (\ref{scp1})}\\
&=& -p(a,c) [[a,b],c] \quad \text{by (\ref{scp1})} \ .
\end{eqnarray*}
By summing the three equalities, we get Jacobi, as expected.
\end{proof}

Any superbracket satisfying the above properties gives any associative algebra the structure of Lie superalgebra.

\begin{remarque}
The even part $\mathfrak{g}_{\bar{0}}$ of a Lie superalgebra is just a standard Lie algebra.  Indeed, for two elements, either even or with opposite parities, the superbracket acts as a (antisymmetric) commutator.  On the contrary, for two odd elements, the superbracket acts as a (symmetric) commutator.
\end{remarque}

\section{Local family}
Let $\mathfrak{g}$ be a Lie superalgebra.  We first extend the Lie superbracket on $\mathfrak{g}$ to the commutator between two $\mathfrak{g}$-valued formal distributions in one indeterminate.  Starting from $a(z)=\sum_m a_{(m)} z^{-1-m} \in \mathfrak{g}[[z,z^{-1}]]$ and $b(w)=\sum_n b_{(n)} w^{-1-n} \in \mathfrak{g}[[w,w^{-1}]]$, we define a new formal distribution in two indeterminates by posing 
\begin{eqnarray}
[a(z),b(w)] \doteq \sum_{m,n} [a_{(m)},b_{(n)}] z^{-1-m} w^{-1-n} \,\in\mathfrak{g}[[z^{\pm 1},w^{\pm 1}]] \label{blablaa}
\end{eqnarray}
We can now define the notion of \textbf{locality of a pair} of formal distributions.
\begin{defin}
Let $\mathfrak{g}$ be a Lie superalgebra.  A pair $(a(z)$,$b(w))$ of $\mathfrak{g}$-valued formal distributions is said local\index{locality!local pair} if $[a(z),b(w)]$
is local.
\end{defin}

By the decomposition theorem, it means that
\begin{eqnarray}
[a(z), b(w)]=\sum_{j\in\corps{Z}^+} c^j(w) \frac{\partial_w^j \delta(z,w)}{j!} \ , \label{mutlocc}
\end{eqnarray}
with $c^j(w)=\Res{z} (z-w)^j [a(z), b(w)] \, \in \mathfrak{g}[[w^{\pm 1}]]$ for all $j\in \corps{Z}^+$.

As shown before, locality of a formal distribution (in two indeterminates) is a strong assumption, which lead to the theorem \ref{decomp}.  Likewise, locality of a pair of formal distributions (each in one indeterminate) implies strong constraints on the commutator between their Fourier modes, as indicated in the following proposition.
\begin{prop}\label{pairwise}
Let $\mathfrak{g}$ be a Lie superalgebra.  Let $(a(z)$,$b(z))$ be a local pair of $\mathfrak{g}$-valued formal distributions.  Then, the Fourier modes satisfy the following commutation relation on $\mathfrak{g}$ :
\begin{eqnarray}
[a_{(m)},b_{(n)}]=\sum_{j\in\corps{Z}^+} \binom{m}{j} c^j_{(m+n-j)} \ , \label{commexpfour}
\end{eqnarray}
where $c^j(w)=\Res{z} (z-w)^j [a(z),b(w)]$ for $j\in\corps{Z}^+$.
\end{prop}
\begin{proof}
Inserting
\begin{eqnarray*}
\frac{1}{j!} \partial_w^j \delta(z,w) = \sum_{k\in\corps{Z}} \binom{k}{j} z^{-1-k} w^{k-j}
\end{eqnarray*}
in (\ref{mutlocc}), we obtain
\begin{eqnarray}
[a(z), b(w)] = \sum_{k\in\corps{Z}} \sum_{j\in\corps{Z}^+} c^j(w) \binom{k}{j} z^{-1-k} w^{k-j} \label{democom}\ .
\end{eqnarray}
Replacing $c^j(w)$ by its Fourier expansion, i.e.
$$ c^j(w)=\sum_{m\in\corps{Z}} c^j_{(m)} w^{-1-m} \ ,$$ the equality (\ref{democom}) reads
\begin{eqnarray}
[a(z), b(w)] &=& \sum_{k\in\corps{Z}}\sum_{m\in\corps{Z}}\sum_{j\in\corps{Z}^+} c^j_{(m)} \binom{k}{j} z^{-1-k} w^{-1-m+k-j} \\
&=& \sum_{k,\tilde{m}}\sum_{j\in\corps{Z}^+}
c^j_{(\tilde{m}+k-j)} \binom{k}{j} z^{-1-k} w^{-1-\tilde{m}} \ ,
\end{eqnarray}
where $\tilde{m}=m-k+j$.
Comparing with (\ref{blablaa}), we get the expected result.
\end{proof}
\begin{remarque}
It is easy to see that, since $\mathfrak{g}$ is a Lie superalgebra, if $(a,b)$ is a local pair, so is $(b,a)$.  To stress this symmetry, the formal distributions $a$ and $b$ are said \textbf{mutually local}.  This is generally not true for an algebra of any other type.
\end{remarque}
The notion of locality of a pair of formal distributions leads us to the concept of local family\index{locality!local family} of formal distributions.
\begin{defin}
A subset $F\subset \mathfrak{g}[[z,z^{-1}]]$ is called a \emph{local family} of $\mathfrak{g}$-valued formal distributions if all pairs of its constituents are local.
\end{defin}

\section{$j$-products and $\lambda$-bracket}
According to \eqref{mutlocc}, the coefficients $c^j (w)$, with $j\in\corps{Z}^+$, measure of far the mutually local formal distributions $a(z)$ and $b(z)$ are from commuting with each other. These coefficients can be regarded as resulting from a $\corps{C}$-bilinear map, called a \textbf{$j$-product}\index{$j$-product} and defined below.  This will allow us to study the algebraic structure they generate.
\begin{defin}
Let $a(w)$ and $b(w)$ be two $\mathfrak{g}$-valued formal distributions, where $\mathfrak{g}$ is a Lie superalgebra. Their $j$-\emph{product} is the $\corps{C}$-bilinear map 
$$\mathfrak{g}[[w^{\pm 1}]] \otimes \mathfrak{g}[[w^{\pm 1}]] \to  \mathfrak{g}[[w^{\pm 1}]]$$ such that
$$a(w) \otimes b(w) \to a(w)_{(j)} b(w) \ , $$
where $j\in \corps{Z}^+$ and 
$$a(w)_{(j)} b(w) \doteq \Res{z} (z-w)^j [a(z),b(w)] \ .$$
This product will sometimes be simply denoted by $a_{(j)} b$.
\end{defin}
These $j$-products will guide us through the study of a new algebraic structure that will be seen to encode the singular part of an operator product expansion\index{operator product expansion}.

\begin{remarque}
Other types of algebras can be considered. We can then define other types of $j$-products and $\lambda$-product, by replacing $[a(z),b(w)]$ using the corresponding product law.
\end{remarque}

Taking the formal Fourier transform on both parts of \eqref{mutlocc}, we define a new product between $\mathfrak{g}$-valued formal distributions, called a \textbf{$\lambda$-bracket}\index{$\lambda$-bracket} when $\mathfrak{g}$ is a Lie superalgebra.
\begin{defin}
Let $\mathfrak{g}$ be a Lie superalgebra.  The $\lambda$-\emph{bracket} of two $\mathfrak{g}$-valued formal distributions is defined by the $\corps{C}$-bilinear map
$$[\;_\lambda \;] \,:\, \mathfrak{g}[[w^{\pm 1}]] \otimes \mathfrak{g}[[w^{\pm 1}]]  \to \mathfrak{g}[[w^{\pm 1}]][[\lambda]]$$
with $$[a(w)_\lambda b(w)] \doteq \Fz [a(z),b(w)] \ ,$$
where $$\Fz =\Res{z} e^{\lambda (z-w)} \ .$$  
\end{defin}
It can easily be shown that the $\lambda$-bracket is related to the $j$-products by
$$[a_\lambda b] = \sum_{j\in \corps{Z}^+} \frac{\lambda^j}{j!} \left( a_{(j)}b \right) \ .$$
This suggests to see the $\lambda$-bracket as the generating function of the $j$-products.  It allows us to gather all the $j$-products in one product alone, the price to pay being the additional formal variable $\lambda$.
\begin{remarque}
Note that for a local pair, the sum in the expansion of $[a_\lambda b]$ in terms of $j$-products is finite.
\end{remarque}

By defining by $(\partial a)(z)=\partial_z a(z)$ the action of the map $\partial$ on a formal distribution $a(z)\in\mathfrak{g}[[z^{\pm 1}]]$, we obtain the following three propositions.
\begin{prop}\label{o3z2}
The $\lambda$-bracket satisfies the following equalities.
\begin{enumerate}
\renewcommand{\labelenumi}{\emph{(\theenumi)}}
\item \label{o3z2a} $[\partial a_\lambda b]=-\lambda \,\, [a_\lambda b]$
\item \label{o3z2b} $[a_\lambda \partial b]=(\partial + \lambda)\, [a_\lambda b]$
\end{enumerate}
\end{prop}
\begin{proof}
As $[(\partial a)(z),b(w)]\in \mathfrak{g}[[z^{\pm 1},w^{\pm 1}]]$, we can make use of proposition \ref{Four}. For the first equality, we compute
\begin{eqnarray*}
[\partial a_\lambda b] &=& \Fz \, [(\partial a)(z),b(w)]  \\
&=& \Fz \partial_z [a(z),b(w)] \\
&=& -\lambda \Fz ([a(z),b(w)]) \text{ by \ref{Four}(\ref{p002})} \\
&=& -\lambda \,\, [a_\lambda b] \ ,
\end{eqnarray*}

For the second equality, we have 
\begin{eqnarray*}
[a_\lambda \partial b] &=& \Fz [a(z),(\partial b)(w)] \\
&=& \Fz [a(z), \partial_w b(w)]\\
&=& \left( \Fz \partial_w \right)  [a(z),b(w)]\\
&=& \left( -[\partial_w,\Fz] + \partial_w \Fz \right) [a(z),b(w)]\\
&=& \left( \lambda \Fz + \partial_w \Fz \right) [a(z),b(w)] \text{ by \ref{Four}(\ref{p002})} \\
&=& \left( \lambda + \partial_w \right) \Fz [a(z),b(w)] \\
&=& \left( \lambda + \partial_w \right) [a_\lambda b] 
\end{eqnarray*}
\end{proof}
In terms of the $j$-products, the previous proposition is translated as follows.
\begin{prop}\label{o3z1}
The $j$-product satisfies the following equalities.
\begin{enumerate}
\renewcommand{\labelenumi}{\emph{(\theenumi)}}
\item $(\partial a)_{(j)}b=-j a_{(j-1)}b$
\item $a_{(j)}\partial b = \partial (a_{(j)}b)+ja_{(j-1)}b$
\end{enumerate}
\end{prop}
\begin{proof}
The first property is equivalent to \ref{o3z2}(\ref{o3z2a}).  Indeed, starting from
\begin{eqnarray*}
[\partial a_\lambda b]=-\lambda \,\, [a_\lambda b] \ , 
\end{eqnarray*}
the LHS
\begin{eqnarray}
[\partial a_\lambda b] = \sum_{j=0}^\infty \frac{\lambda^j}{j!}\left( \partial a_{(j)}b \right) \ , \label{calculoz1}
\end{eqnarray}
and the RHS,
\begin{eqnarray}
-\lambda \,\, [a_\lambda b] &=& - \sum_{k=0}^\infty \frac{\lambda^{k+1}}{k!}\left( a_{(k)}b \right) \nonumber\\
&=& - \sum_{j=1}^\infty \frac{\lambda^{j}}{(j-1)!}\left( a_{(j-1)}b \right) \quad \text{where }j=k+1 \nonumber\\
&=& \sum_{j=1}^\infty \frac{\lambda^{j}}{(j)!}\left(-j\,\, a_{(j-1)}b \right) \nonumber\\
&=& \sum_{j=0}^\infty \frac{\lambda^{j}}{(j)!}\left(-j\,\, a_{(j-1)}b \right) \label{calculoz2} \ .
\end{eqnarray}
Equating (\ref{calculoz1}) and (\ref{calculoz2}), we obtain the first property.\\
The second property is equivalent to \ref{o3z2}(\ref{o3z2b}).  Thus, we start from 
\begin{eqnarray*}
[a_\lambda \partial b]=(\partial + \lambda)\, [a_\lambda b] \ .
\end{eqnarray*}
The LHS reads
\begin{eqnarray}
[a_\lambda \partial b] = \sum_{j=0}^\infty \frac{\lambda^j}{j!}\left( a_{(j)}\partial b \right) \ , \label{calculoz3}
\end{eqnarray}
and the RHS
\begin{eqnarray}
(\partial + \lambda)\, [a_\lambda b] &=& (\partial + \lambda)\, \sum_{k=0}^\infty \frac{\lambda^k}{k!}\left( a_{(k)}b \right) \nonumber\\
&=& \sum_{k=0}^\infty \frac{\lambda^k}{k!}\partial\left( a_{(k)}b \right) + \sum_{k=0}^\infty \frac{\lambda^{k+1}}{k!} \left( a_{(k)}b \right) \nonumber\\
&=& \sum_{j=0}^\infty \frac{\lambda^j}{j!}\partial\left( a_{(j)}b \right) + \sum_{j=1}^\infty \frac{\lambda^{j}}{j!} \left(j \, a_{(j-1)}b \right) \quad \text{where }j=k+1 \nonumber \\
&=& \sum_{j=0}^\infty \frac{\lambda^j}{j!}\partial\left( a_{(j)}b \right) + \sum_{j=0}^\infty \frac{\lambda^{j}}{j!} \left(j \, a_{(j-1)}b \right) \label{calculoz4}\ .
\end{eqnarray}
Comparing (\ref{calculoz3}) and (\ref{calculoz4}), we get the second property.
\end{proof}
We deduce from the last proposition the following important result.
\begin{prop}\label{o3z3}
The map $\partial$ acts as a derivation on the $\lambda$-bracket and the $j$-products, i.e.
$\partial [a_\lambda b] = [\partial a_\lambda b] + [a_\lambda \partial b]$ and $\partial(a_{(j)}b)=(\partial a)_{(j)}b + a_{(j)}\partial b$.
\end{prop}
\begin{proof}
Clearly, from the two previous propositions, we have
$$ \partial\left(a_\lambda b\right) = \partial a_\lambda b + a_\lambda \partial b \quad \text{et}\quad\partial(a_{(j)}b)=(\partial a)_{(j)}b + a_{(j)}\partial b \ , $$ 
which concludes the proof.
\end{proof}
\begin{remarque}
It is worth noticing that, according to the previous propositions, neither the $\lambda$-bracket, nor the $j$-products, are $\corps{C}[\partial]$-bilinear.  On the contrary, the map $a(z)\otimes b(w) \to [a(z),b(w)]$ is $\corps{C}[\partial]$-bilinear, and this comes from the fact that $\mathfrak{g}[[z^{\pm 1},w^{\pm 1}]]$ is a $\corps{C}[\partial]$-module.
\end{remarque}
\begin{remarque}
All the Fourier modes of any formal distribution are assumed to have the same parity.  Thus, for $a(z)=\sum_n a_{(n)} z^{-n-1}$, we pose $p(a)=p(a_{(n)})$ for all  $n\in \corps{Z}$.  If $b(z)=\sum_m b_{(m)} z^{-m-1}$, we pose $p(a,b)=p(a_{(n)},b_{(m)})$ for all $n,m \in \corps{Z}$.
\end{remarque}  
We now address the translation of the antisymmetry property and the Jacobi identity satisfied by $\mathfrak{g}$ in the language of the $\lambda$-bracket and $j$-products of $\mathfrak{g}$-valued formal distributions.
\begin{prop} \label{lprop}
Let $\mathfrak{g}$ be an associative Lie superalgebra.
The $\lambda$-bracket between $\mathfrak{g}$-valued formal distributions satisfies the following properties.
\begin{enumerate}
\renewcommand{\labelenumi}{\emph{(\theenumi)}}
\item \label{lprop2} $[b_\lambda a]=-p(a,b) [a_{-\lambda -\partial} b]$ \textbf{if} $(a,b)$ is a local pair;
\item \label{lprop3} $[a_\lambda [b_\mu c]]=[[a_\lambda b]_{\lambda+\mu} c] + p(a,b) [b_\mu [a_\lambda c]]$.
\end{enumerate}
\end{prop} 
\begin{proof}
The antisymmetry of the Lie superbracket implies ${[a(z),b(w)]=-p(a,b)[b(w),a(z)]}$. Applying first the Fourier transform in two indeterminate to both sides and then making use of property \ref{Four}(\ref{p003}), as $(a,b)$ is local, we get the first property.  

For the second property, we translate the Jacobi identity \ref{scp}(\ref{scp3}), for which $\mathfrak{g}$ has to be associative, in terms of formal distributions : $[a(z),[b(x),c(w)]] = [[a(z),b(x)],c(w)] + p(a,b) [b(x),[a(z),c(w)]]$.  
We then have
\begin{eqnarray*}
[a_\lambda [ b_\mu c ]] &=& \Fz \, [ a(z) , \F{x}{w}{\mu} [(b(x), c(w)] ] \\
&=& \Fz \F{x}{w}{\mu} \,\, [ a(z) , [(b(x), c(w)] ] \\
&=& \Fz \F{x}{w}{\mu} \,\, \Big( [[a(z),b(x)],c(w)] + p(a,b) [b(x),[a(z),c(w)]] \Big) \ .
\end{eqnarray*}
Since $[[a(z),b(x)],c(w)] \in \mathfrak{g}[[z^{\pm 1},x^{\pm 1},w^{\pm 1}]]$, we can make use of the property \ref{Four}(\ref{p004}), according to which $\Fz \F{x}{w}{\mu}=\F{x}{w}{\lambda + \mu} \F{z}{x}{\lambda}$.  Therefore,
\begin{eqnarray*}
[a_\lambda [ b_\mu c ]] &=& \F{x}{w}{\lambda + \mu} \, \F{z}{x}{\lambda}\,\, [[a(z),b(x)],c(w)] + p(a,b) \Fz \F{x}{w}{\mu}  [b(x),[a(z),c(w)]]\\
&=& \F{x}{w}{\lambda + \mu}\,\,[\F{z}{x}{\lambda} [a(z), b(x)] , c(w)] + p(a,b) \F{x}{w}{\mu} [b(x),\Fz [a(z),c(w)]] \\
&=& [[a_\lambda b]_{\lambda+\mu} c] + p(a,b) [b_\mu [a_\lambda c]] \ ,
\end{eqnarray*}
which proves the second property.
\end{proof}

In terms of $j$-products, we obtain the following proposition.
\begin{prop}\label{jprop}
Let $\mathfrak{g}$ be an associative Lie superalgebra.  The $j$-products between $\mathfrak{g}$-valued formal distributions satisfy the following properties ($j\in \corps{Z}^+$).
\begin{enumerate}
\renewcommand{\labelenumi}{\emph{(\theenumi)}}
\item $b_{(j)}a=-p(a,b)\sum_{l=0}^\infty (-1)^{j+l} \frac{\partial^l(a_{(j+l)}b)}{l!}$ \textbf{if} $(a,b)$ is a local pair;
\item $a_{(p)}(b_{(m)}c)=\sum_{k=0}^p \binom{p}{k}(a_{(k)}b)_{(p+m-k)}c + p(a,b) b_{(m)}(a_{(p)}c)$.
\end{enumerate}
\end{prop}
\begin{proof}
These properties are equivalent to those in the previous proposition.  Their translation in terms of $j$-products is as follows. We start from 
$$[b_\lambda a]=-p(a,b) [a_{-\lambda -\partial} b]$$
The LHS reads
\begin{eqnarray}
[b_\lambda a] = \sum_{j=0}^\infty \frac{\lambda^j}{j!}\left( b_{(j)}a \right) \ , \label{calculoz5a}
\end{eqnarray}
and the RHS, 
\begin{eqnarray}
-p(a,b)[a_{-\lambda-\partial} b] &=& -p(a,b)\sum_{k=0}^\infty \frac{(-\lambda-\partial)^k}{k!} a_{(k)}b \nonumber\\
&=& -p(a,b)\sum_{k=0}^\infty \frac{(-1)^k}{k!} \sum_{l=0}^k \binom{k}{l} \lambda^{k-l} \partial^l \left(a_{(k)}b\right) \nonumber\\
&=& -p(a,b) \sum_{k=0}^\infty \frac{(-1)^k}{l!(k-l)!} \sum_{l=0}^k  \lambda^{k-l} \partial^l \left(a_{(k)}b\right) \nonumber \\
&=& -p(a,b)\sum_{j=0}^\infty \frac{\lambda^{j}}{j!} \sum_{l=0}^\infty (-1)^{j+l} \frac{\partial^l \left(a_{(j+l)}b\right)}{l!}  \quad \text{where $j=k-l$} \label{calculoz5b} \ .
\end{eqnarray}
Comparing (\ref{calculoz5a}) and (\ref{calculoz5b}), we get the result.  The second property is analogous.
\end{proof}
\begin{remarque}
Note that the second property of this proposition is just \eqref{commexpfour}, where the Lie superbracket is defined in the usual way for an associative superalgebra.
\end{remarque}

\begin{remarque}
At section \ref{fourtrrr}, we studied properties of the Fourier transform in one indeterminate.  We are now able to derive another one, which stresses the symmetry of its action on a $\lambda$-bracket. 
\end{remarque}
\begin{prop}Let $\mathfrak{g}$ be an associative Lie superalgebra and $a(z)$, $b(z)$ $\mathfrak{g}$-valued formal distributions. The Fourier transform in one indeterminate acts on the $\lambda$-bracket as follows.
$$ \Fzd{z}{\lambda+\mu} [a(z)_\lambda b(z)]= [\Fzz a(z), \Fzd{z}{\mu} b(z)] $$
\end{prop}
\begin{proof}
This results from the following computation.
\begin{eqnarray*}
\Fzd{z}{\lambda + \mu} [a(z)_\lambda b(z)] &=& \Res{z} e^{(\lambda + \mu)z} [a(z)_\lambda b(z)] \\
&=& \Res{z} e^{(\lambda + \mu)z} \text{Res}_x e^{\lambda(x-z)}[a(x),b(z)] \\
&=& \Res{z} e^{(\lambda + \mu)z} \text{Res}_x e^{\lambda(x-z)} \left( a(x)b(z) - p(a,b) b(z) a(x) \right) \\
&=& \Res{z}\text{Res}_x \left( e^{\lambda x} a(x) e^{\mu z} b(z) - p(a,b) e^{\mu z} b(z) e^{\lambda x} a(x) \right) \\
&=& \left( \Res{x} e^{\lambda x} a(x) \right) \left( \Res{z} e^{\mu z} b(z) \right) - p(a,b) \left(  \Res{z} e^{\mu z} b(z) \right) \left(\Res{x} e^{\lambda x} a(x) \right) \\
&=& [\Fzz a(z), \Fzd{z}{\mu} b(z)] \ ,
\end{eqnarray*}
where we replaced $x$ by $z$ at the last line, keeping the parentheses. This is legitimate since the result does not depend on them.
\end{proof}

%
%
\section{Formal distributions Lie superalgebra}
The notion of local family of $\mathfrak{g}$-valued formal distributions allows us to characterize $\mathfrak{g}$ itself.
\begin{defin}
Let $\mathfrak{g}$ be a Lie superalgebra. If there exists a local family $F$ of $\mathfrak{g}$-valued formal distributions,  with their Fourier coefficients generating the whole $\mathfrak{g}$, then $F$ is said to endow $\mathfrak{g}$ with the structure of \emph{Formal distributions Lie superalgebra}\index{formal distributions Lie superalgebra}.  To insist on the role of $F$, we denote $\mathfrak{g}$ as ($\mathfrak{g}$,$F$).
\end{defin}

Starting from $F$, we denote by $\bar{F}$ the minimal subspace of $\mathfrak{g}[[z,z^{-1}]]$ containing $F$ and closed under their $j$-products.  $\bar{F}$ is called the \textit{minimal local family}\index{locality!minimal local family} of $F$.  It is indeed remarkable for $\bar{F}$ to be a local family as well.  
\begin{prop}\label{propdong}
Let $\mathfrak{g}$ be an associative and/or Lie superalgebra. Let $F\subset \mathfrak{g}[[z,z^{-1}]]$ be a local family.  Then $\bar{F}$ is a local family as well.
\end{prop}
\begin{proof}
Let $a(z),b(z),c(z)\in F$.  We have to show that $(a_{(n)}b,c)$ is then a local pair. But this results from Dong's lemma, which we will prove at a later section (after having defined the notion of generalized $j$-products).  
\end{proof}
\begin{exe}\label{constriv}
To see what is going on, let's take a simple example.
Starting from a Lie superalgebra $\mathfrak{g}$, let's define $F$ as follows :
$$F=\{ g(z)=\sum_{n\in\corps{Z}} g z^{-n-1} | g \in \mathfrak{g} \} \ .$$ 
This set is a local family.  Indeed, for $g,h\in\mathfrak{g}$, we have $[g(z),h(w)]=[g,h](w) \delta(z,w)$.  Posing $m=k-n$ at the first line below, we get 
\begin{eqnarray*}
[g(z),h(w)] &=& [g,h] \sum_{n,m\in\corps{Z}} z^{-n-1} w^{-m-1} \\
&=& [g,h] \sum_{n,k\in\corps{Z}} z^{-n-1} w^{-k+n-1} \\
&=& \sum_{k\in\corps{Z}} [g,h] w^{-k-1} \sum_{n\in\corps{Z}} z^{-1-n}w^n \\
&=& [g,h](w) \delta(z,w) \ .
\end{eqnarray*}
From this computation, we see that the $0$-product alone is non vanishing and already in $F$.  Therefore, $F$ is closed under the $j$-products : $\bar{F}=F$. Then $(\mathfrak{g},F)$ is a $\mathfrak{g}$-valued formal distributions Lie superalgebra.
\end{exe}

\begin{remarque}
The previous example shows that we can trivially associate to any Lie superalgebra $\mathfrak{g}$ a structure of formal distributions Lie superalgebra ($\mathfrak{g},F$).  However, in most cases, $F$ is restricted to contain a finite number of formal distributions.  As their Fourier coefficients have to generate the whole Lie superalgebra $\mathfrak{g}$, it is clear that the previous construction is useless in the case of an infinite dimensional algebra.  In other words, if $\mathfrak{g}$ is finite dimensional, then the structure of formal distributions Lie superalgebra does not bring any more good.  On the contrary, it leads to non trivial conditions in the case of an infinite dimensional Lie superalgebra.
\end{remarque}

\section{Lie conformal superalgebra}

Now, we take into consideration the invariance of the (minimal) local family under the derivation $\partial$.
\begin{defin}
Let ($\mathfrak{g}$,$F$) be a formal distribution Lie superalgebra and $\bar{F}$ the minimal local family containing $F$.  The subset of $\mathfrak{g}[[z^{\pm 1}]]$ containing $\bar{F}$ and being closed under the derivation $\partial_z$ (up to an arbitrary order) is called a \emph{conformal family}\index{conformal family}. This family is said to be generated by $\bar{F}$ as a $\corps{C}[\partial_z]$-module and is then denoted by $\corps{C}[\partial_z]\bar{F}$.
\end{defin}

The notions of $j$-product and $\lambda$-bracket were previously defined from $\mathfrak{g}$-valued formal distributions, with $\mathfrak{g}$ being a given superalgebra. Those products were shown to satisfy several properties, coming either from their definition, or from the nature of $\mathfrak{g}$ itself.  In this section, we take those properties as axioms of a new algebraic structure, defined intrinsically, without any reference neither to $\mathfrak{g}$, nor to $\mathfrak{g}$-valued formal distributions. For this reason, we write $\partial$ instead of $\partial_z$ in the following definition.



\begin{defin}\label{rappeldefin}
A $\corps{C}[\partial]$-module $\cal R$ is called a \emph{Lie conformal superalgebra}\index{Lie conformal superalgebra} if it is endowed with a $\corps{C}$-bilinear map, called the $\lambda$-bracket,
$$[\,_\lambda \,]\,:\,{\cal R}\otimes{\cal R} \to \corps{C}[\lambda]\otimes {\cal R} \ ,$$  
and satisfying the following properties, where $a,b,c\in{\cal R}$,
\begin{eqnarray}
\lbrack  \partial a_\lambda b\rbrack   &=& -\lambda \lbrack  a_\lambda b\rbrack   \label{cmod01}\\
\lbrack  a_\lambda \partial b\rbrack   &=& (\partial + \lambda)\lbrack  a_\lambda b\rbrack   \label{cmod02}\\
\lbrack  b_\lambda a\rbrack   &=& -p(a,b) \lbrack  a_{-\lambda - \partial}b\rbrack   \label{skcomml}\\
\lbrack  a_\lambda\lbrack  b_\mu c\rbrack  \rbrack   &=& \lbrack  \lbrack  a_\lambda b\rbrack  _{\lambda + \mu}c\rbrack   + p(a,b) \lbrack  b_\mu \lbrack  a_\lambda c\rbrack  \rbrack   \ .\label{jacol}
\end{eqnarray}
With 
$$[a_\lambda b]=\sum_{j\in\corps{Z}^+} \frac{\lambda^j}{j!}(a_{(j)}b) \ ,$$
those properties translate in terms of $j$-products as follows.
\begin{eqnarray}
(\partial a)_{(j)}b &=& -j a_{(j-1)}b   \label{cmod1}\\
a_{(j)}\partial b &=& \partial (a_{(j)}b)+ja_{(j-1)}b   \label{cmod2}\\
b_{(j)}a &=& -p(a,b)\sum_{l=0}^\infty (-1)^{j+l} \frac{\partial^l(a_{(j+l)}b)}{l!}   \label{skcomm}\\
a_{(p)}(b_{(m)}c) &=& \sum_{k=0}^p \binom{p}{k}(a_{(k)}b)_{(p+m-k)}c + p(a,b) b_{(m)}(a_{(p)}c)  \label{jaco} \ .
\end{eqnarray}
\end{defin}
\begin{remarque}
\cite{dandrea} The equalities (\ref{cmod01})-(\ref{skcomml}) are not independent.  Indeed, (\ref{cmod01}) and (\ref{skcomml}) together imply (\ref{cmod02}) : 
$$[a_\lambda \partial b] = - p(a,b)[\partial b_{-\lambda - \partial}\, a] = p(a,b)(-\lambda-\partial)[b_{-\lambda - \partial}\,a]=(\partial + \lambda)[a_\lambda b] \ .$$  
The same is true in terms of $j$-products.
\end{remarque}



\section{Maximal formal distributions Lie algebra}
We previously shown that any formal distribution Lie algebra ($\mathfrak{g},F$) was given the structure of Lie Conformal algebra ${\cal R}$, with ${\cal R}=\corps{C}[\partial_z]\bar{F}$, $\partial = \partial_z$ and $[a_\lambda b]=\Fz [a(z),b(w)]$. It turns out that the process can be reversed : to any Lie conformal algebra we can associate a formal distributions Lie algebra.  

According to the definition of a formal distribution Lie algebra, we first have to define a Lie algebra, denoted by $\text{Lie}({\cal R})$, and then associate to it a conformal family $\bar{{\cal R}}$ of $\text{Lie}({\cal R})$-valued formal distributions, whose coefficients span $\text{Lie}({\cal R})$, so that ($\text{Lie}({\cal R})$,$\bar{{\cal R}}$) is then the expected formal distribution Lie algebra.  The latter will then be shown to be \textit{maximal} in a sense to specify.

To define $\text{Lie}({\cal R})$, we have to consider a linear space (over $\corps{C}$) endowed with a suitable Lie bracket.  One way to find it is to recall proposition \ref{pairwise}, according to which
\begin{eqnarray}
[a_{(m)}, b_{(n)}]=\sum_{j\in\corps{Z}^+} \binom{m}{j}\left( a_{(j)}b \right)_{(m+n-j)} \ , \label{commre}
\end{eqnarray}
where $a(z)=\sum_{m}a_{(m)}z^{-m-1}$ and $b(z)=\sum_{n}b_{(n)}z^{-n-1}$ were mutually local $\mathfrak{g}$-valued formal distributions and $\mathfrak{g}$ was a \textit{given} Lie algebra.  Conversely, starting from a Lie conformal algebra $\cal R$, which is by definition endowed with $j$-products satisfying (\ref{cmod1})-(\ref{jaco}), equation \eqref{commre} will then be shown to provide a proper Lie bracket on a Lie algebra \textit{to be specified}.
To specify the suitable Lie algebra acted on by this bracket, we proceed in two steps.  

We first consider the space $\tilde{{\cal R}}={\cal R}\otimes\corps{C}[t,t^{-1}]$, with $\tilde{\partial}=\partial \otimes I + I \otimes \partial_t$, where $I$ appearing on the left (resp. right) of $\otimes$ is the identity operator acting on $\cal R$ (resp. $\corps{C}[t,t^{-1}]$). This space is called the \textit{affinization} of $\cal R$.  Its generating elements can be written $a\otimes t^m$, where $a\in{\cal R}$ and $m\in\corps{Z}$. For clarity, we will use the notation $\tilde{{\cal R}}={\cal R}[t,t^{-1}]$, $at^m$ for its elements and $\tilde{\partial}=\partial + \partial_t$.  We set $a_{(m)}=at^m$ and $b_{(n)}=bt^n$ in (\ref{commre}), so that we obtain a well defined commutation relation on $\tilde{{\cal R}}$ :
\begin{eqnarray}
[at^m, bt^n]=\sum_{j\in\corps{Z}^+} \binom{m}{j}\left( a_{(j)}b \right) t^{m+n-j} \ , \label{commre1}
\end{eqnarray}
which gives $\tilde{{\cal R}}$ a structure of algebra, denoted by ($\tilde{{\cal R}}$,$[\, ,\, ]$).  As expected,  $\tilde{{\cal R}}$ is a linear space over $\corps{C}$.  Indeed, by the definition of affinization, for $a_{(n)},b_{(n)} \in \tilde{{\cal R}}$, we have $\lambda a_{(n)}=(\lambda a)_{(n)}$ and $a_{(n)}+b_{(n)}=(a+b)_{(n)}$. Therefore, $\lambda a_{(n)} \in \tilde{{\cal R}}$ and $a_{(n)}+b_{(n)} \in \tilde{{\cal R}}$. 

Now the second step.  We have to check that the commutator verifies the antisymmetry and Jacobi identities, considering that the terms $a_{(j)}b$ of the RHS of (\ref{commre1}) satisfy the relations (\ref{skcomm}) and (\ref{jaco}).  The latter ones are not sufficient.  Indeed, as we will see it, another constraint has to be imposed on elements of $\tilde{{\cal R}}$, namely $\tilde{\partial}(at^m)=0$.  The algebraic formulation of the latter condition is as follows : the space $\tilde{{\cal R}}$ has to be quotiented by the subspace $I$ spanned by elements of the form $\{(\partial a)t^n + n at^{n-1} \,|\,n\in\corps{Z} \}$.  Using $\tilde{\partial}$, we can write $I=\tilde{\partial}\tilde{{\cal R}}$.
This process has thus two goals : first transfering on $\tilde{{\cal R}}/\tilde{\partial}\tilde{{\cal R}}$ the structure of algebra of ($\tilde{{\cal R}}$,$[\, ,\, ]$), and then endowing ($\tilde{{\cal R}}/\tilde{\partial}\tilde{{\cal R}}$,$[\, ,\, ]$) with the structure of Lie algebra.  The first goal is not direct.  Indeed, $\tilde{\partial}\tilde{{\cal R}}$ has to be a two-sided ideal of the algebra ($\tilde{{\cal R}}$,$[\, ,\, ]$), which is the case.
\begin{prop}
$\tilde{\partial}\tilde{{\cal R}}$ is a two-sided ideal of the algebra $(\tilde{{\cal R}}$,$[\, ,\, ])$, where $[\, ,\, ]$ is defined by (\ref{commre1}).
\end{prop}
\begin{proof}
$\tilde{\partial}\tilde{{\cal R}}$ is spanned by elements of the form $(\partial a)t^n + n at^{n-1}$.  Replacing $at^m$ by one of them in (\ref{commre1}), we obtain :
{\setlength\arraycolsep{0pt}
\begin{eqnarray*}
[(\partial a)t^m &&+ mat^{m-1},bt^n] = \sum_{j\in\corps{Z}^+} \binom{m}{j}(\partial a_{(j)}b)t^{m+n-j} + m \sum_{k\in\corps{Z}^+} \binom{m-1}{k}(a_{(k)}b)t^{m-1+n-k} \\
&&=\; - \sum_{j\geq 1} j \binom{m}{j}(a_{(j-1)}b)t^{m+n-j} + m \sum_{k\geq 0} \binom{m-1}{k}(a_{(k)}b)t^{m+n-k-1} \quad \text{by (\ref{cmod1})} \\
&&=\; - \sum_{j\geq 0} (j+1)\binom{m}{j+1}  (a_{(j)}b)t^{m+n-j-1} + m \sum_{k\geq 0} \binom{m-1}{k}(a_{(k)}b)t^{m+n-k-1}  \\
&&=\; - m \sum_{j\geq 0} \binom{m-1}{j} (a_{(j)}b)t^{m+n-j-1} + m \sum_{k\geq 0} \binom{m-1}{k}(a_{(k)}b)t^{m+n-k-1}  \\
&&=\; 0 \in \tilde{\partial}\tilde{{\cal R}}
\end{eqnarray*}}
Hence $\tilde{\partial}{\cal R}$ is a left ideal.  The following computation shows that it is a right ideal as well.
{\setlength\arraycolsep{0pt}
\begin{eqnarray*}
[at^m,&&(\partial b)t^n + nbt^{n-1}] = \sum_{j\geq 0} \binom{m}{j}(a_{(j)}\partial b)t^{m+n-j} + n \sum_{k\geq 0} \binom{m}{k}(a_{(k)}b)t^{m+n-1-k} \\
&&=\; \sum_{j\geq 0} \binom{m}{j} \partial(a_{(j)}b) t^{m+n-j} + \sum_{j\geq 0} j \binom{m}{j} (a_{(j-1)}b)t^{m+n-j}\\
&& \quad\quad+ \; n\sum_{j\geq 0} \binom{m}{j} (a_{(j)}b)t^{m+n-1-j}  \quad \text{by (\ref{cmod2})} \\
&&=\; \underbrace{\sum_{j\geq 0} \binom{m}{j} \partial(a_{(j)}b) t^{m+n-j} + \sum_{j\geq 0}\binom{m}{j}(m+n-j)(a_{(j)}b)t^{m+n-j-1}}_{\in\tilde{\partial}\tilde{{\cal R}} }\\
&&\quad\quad - \; \sum_{j\geq 0}\binom{m}{j}(m+n-j)(a_{(j)}b)t^{m+n-j-1} + \sum_{j\geq 0} j \binom{m}{j} (a_{(j-1)}b)t^{m+n-j}\\
&&\quad\quad + \; n\sum_{j\geq 0} \binom{m}{j} (a_{(j)}b)t^{m+n-1-j}\\
&&=\; \underbrace{\ldots}_{\in\tilde{\partial}\tilde{{\cal R}}} - \sum_{j\geq 0} \binom{m}{j}(m-j)(a_{(j)}b)t^{m+n-j-1} + \sum_{j\geq 0} j\binom{m}{j} (a_{(j-1)}b)t^{m+n-j} \\
&&=\; \underbrace{\ldots}_{\in\tilde{\partial}\tilde{{\cal R}}} - \sum_{j\geq 0} \binom{m}{j}(m-j)(a_{(j)}b)t^{m+n-j-1} + \sum_{j\geq 0} \binom{m}{j}(m-j)(a_{(j)}b)t^{m+n-j-1} \\
&&=\; \ldots \in\tilde{\partial}\tilde{{\cal R}}
\end{eqnarray*}}
\end{proof}
The next proposition shows that $(\tilde{{\cal R}}/\tilde{\partial}\tilde{{\cal R}}$,$[\, ,\, ])$ is a Lie algebra. We will make use of the following (well known) lemma.
\begin{lemme}[Vandermonde formula] Let $n,m\in\corps{Z}$.  Then
$$ \sum_{k=0}^\infty \binom{n}{k}\binom{m}{j-k}=\binom{n+m}{j} $$
\end{lemme}
\begin{proof}
Since $(x+1)^{m+n}=(x+1)^{n}(x+1)^{m}$,
\begin{eqnarray*}
\sum_{j\geq 0} \binom{m+n}{j}x^j &=& \sum_{k\geq0} \binom{n}{k}x^k \sum_{p\geq 0}\binom{m}{p}x^p \\
&=& \sum_k \sum_p \binom{n}{k}\binom{m}{p} x^{k+p}\\
&=& \sum_j \sum_k \binom{n}{k}\binom{m}{j-k} x^j 
\end{eqnarray*}
Therefore, comparing coefficients of the monomial $x^j$, we get the formula. Recall that the expansion $(x+1)^n = \sum_{k\geq0} \binom{n}{k}x^k$ holds for an integer $n<0$.  The Vandermonde formula thus holds for $n,m\in\corps{Z}$.
\end{proof}

Define the homomorphism $\phi \, :\, \tilde{{\cal R}} \to \tilde{{\cal R}}/\tilde{\partial}\tilde{{\cal R}}$, the commutator between two elements of $\tilde{{\cal R}}/\tilde{\partial}\tilde{{\cal R}}$ is defined by
\begin{eqnarray}
[\phi(at^m), \phi(bt^n)]=\sum_{j\in\corps{Z}^+} \binom{m}{j} \phi\left(\left( a_{(j)}b \right) t^{m+n-j}\right) \ , \text{ where $a,b\in{\cal R}$.}\label{commre2}
\end{eqnarray}
We are now all set to show that ($\tilde{{\cal R}}/\tilde{\partial}\tilde{{\cal R}}$,$[\, ,\, ]$) is a Lie algebra.
\begin{prop}\label{7z98d}
$(\tilde{{\cal R}}/\tilde{\partial}\tilde{{\cal R}}$,$[\, ,\, ])$ is a Lie algebra, where $[\, ,\, ]$ is defined by (\ref{commre2}).  Explicitly, we have
\begin{enumerate}
\renewcommand{\labelenumi}{\emph{(\theenumi)}}
\item $\tilde{{\cal R}}/\tilde{\partial}\tilde{{\cal R}}$ is linear space, closed under $[\, ,\, ]$,
\item $[\phi(bt^n),\phi(at^m)] = -p(a,b) [\phi(at^m),\phi(bt^n)]$,
\item $[\phi(at^m),[\phi(bt^n),\phi(ct^p)]=[[\phi(at^m),\phi(bt^n)],\phi(ct^p) + p(a,b)[\phi(bt^n),[\phi(at^m),\phi(ct^p)]$.
\end{enumerate}
\end{prop}
\begin{proof}
$\tilde{{\cal R}}/\tilde{\partial}\tilde{{\cal R}}$ is indeed a linear space over $\corps{C}$, as $\tilde{{\cal R}}$.  It is closed under $[\, ,\, ]$, since $\tilde{\partial}\tilde{{\cal R}}$ is a two-sided ideal of the algebra ($\tilde{{\cal R}}$,$[\, ,\, ]$), which proves (1).\\  
Let's prove (2).  By (\ref{skcomm}), we have
\begin{eqnarray*}
b_{(j)}a=-p(a,b) \sum_{l=0}^\infty (-1)^{j+l} \frac{\partial^l \left( a_{(j+l)}b \right)}{l!} \ .
\end{eqnarray*}
Using the latter in 
$$ [bt^n,at^m]= \sum_{j=0}^\infty \binom{n}{j} \left( b_{(j)}a \right) t^{m+n-j} \ , $$
we obtain
$$ [bt^n,at^m]=-p(a,b) \sum_{j=0}^\infty \sum_{l=0}^\infty \binom{n}{j}(-1)^{j+l} \frac{\partial^l \left( a_{(j+l)}b \right)}{l!} t^{n+m-j} \ . $$
But, since $(\partial a)t^n = -n a t^{n-1} + n a t^{n-1} + (\partial a)t^n = -n a t^{n-1} + \tilde{\partial}\tilde{{\cal R}}$, 
\begin{eqnarray*}
(\partial^l a)t^n &=& (-1)^l n(n-1)\cdots(n-l+1) a t^{n-l} + \tilde{\partial}\tilde{{\cal R}}\\
&=& (-1)^l \binom{n}{l} l! \; at^{n-1} + \tilde{\partial}\tilde{{\cal R}} \ . 
\end{eqnarray*} 
Hence,
$$ \partial^l \left( a_{(j+l)}b \right) t^{n+m-j} = (-1)^l \binom{n+m-j}{l} l! \left( a_{(j+l)}b \right) t^{n+m-j-l} + \tilde{\partial}\tilde{{\cal R}} \ . $$
Therefore, 
\begin{eqnarray*}
[bt^n,at^m] &=& -p(a,b) \sum_{j=0}^\infty \sum_{l=0}^\infty \binom{n}{j}(-1)^{j} \binom{n+m-j}{l}  \left( a_{(j+l)}b \right)  t^{n+m-j-l} + \tilde{\partial}\tilde{{\cal R}} \\
&=& -p(a,b) \sum_{\tilde{j}=0}^\infty \underbrace{\sum_{l=0}^\infty \binom{n}{\tilde{j}-l}(-1)^{\tilde{j}-l} \binom{n+m-\tilde{j}+l}{l}}_{(*)}  \left( a_{(\tilde{j})}b \right)  t^{n+m-\tilde{j}} + \tilde{\partial}\tilde{{\cal R}} \ .
\end{eqnarray*}
Replacing $(-1)^{\tilde{j}-l}$ by $(-1)^{\tilde{j}+l}$, we can write
\begin{eqnarray*}
(*) &=& (-1)^{\tilde{j}} \sum_{l=0}^\infty \binom{n}{\tilde{j}-l}(-1)^l \binom{n+m-\tilde{j}+l}{l} \\
&=& (-1)^{\tilde{j}} \sum_{l=0}^\infty \binom{n}{\tilde{j}-l} \binom{-n-m+\tilde{j}-1}{l} \ ,
\end{eqnarray*}
where we made use of the following result (easy to prove) : $\binom{r}{s}=(-1)^s \binom{s-r-1}{s}$ ($**$). \\ By Vandermonde formula, we can write
\begin{eqnarray*}
(*) &=& (-1)^{\tilde{j}} \sum_{l=0}^\infty \binom{n}{\tilde{j}-l} \binom{-n-m+\tilde{j}-1}{l} \\
&=& (-1)^{\tilde{j}} \binom{-m+\tilde{j}-1}{\tilde{j}} \\
&=& \binom{m}{\tilde{j}} \quad \text{by ($**$)}\ .
\end{eqnarray*}
Hence
\begin{eqnarray*}
[bt^n,at^m] &=& -p(a,b) \sum_{\tilde{j}=0}^\infty \binom{m}{\tilde{j}} \left( a_{(\tilde{j})}b \right)  t^{n+m-\tilde{j}}  + \tilde{\partial}\tilde{{\cal R}}\\
&=& -p(a,b)[at^m,bt^n] + \tilde{\partial}\tilde{{\cal R}} \ .
\end{eqnarray*}
Applying $\phi$ to both members, we obtain (2) since $\phi$ is an homomorphism of algebras and $\phi(\tilde{\partial}\tilde{{\cal R}})=0$.\\ Let's prove (3).  By definition, 
$$ [at^m,[bt^n,ct^p]] = \sum_{j=0}^\infty \sum_{k=0}^\infty \binom{n}{j}\binom{m}{k} \left( a_{(k)}\left( b_{(j)}c \right) \right) t^{m+n+p-j-k} \ .$$
Using (\ref{jaco}), we obtain
\begin{eqnarray} 
[at^m,[bt^n,ct^p]] &=& \underbrace{\sum_{j=0}^\infty \sum_{k=0}^\infty \sum_{l=0}^k \binom{n}{j}\binom{m}{k} \binom{k}{l} \left( \left( a_{(l)}b \right)_{(k+j-l)}c \right) t^{m+n+p-j-k}}_{(*)} \nonumber \\
&&\quad\quad + p(a,b) \sum_{j=0}^\infty \sum_{k=0}^\infty \binom{n}{j}\binom{m}{k} \left( b_{(j)}\left( a_{(k)}c \right) \right) t^{m+n+p-j-k} \label{compeq1}
\end{eqnarray}
It is easy to see that the second term of the RHS is $p(a,b)[bt^n,[at^m,ct^p]]$.  Moreover, setting $\tilde{k}=j+k-l$, ($*$) can be written
\begin{eqnarray*}
(*) &=& \sum_{l=0}^\infty \sum_{\tilde{k}=0}^\infty \sum_{j=0}^\infty \binom{n}{j}\binom{m}{\tilde{k}-j+l}\binom{\tilde{k}-j+l}{l} \left( \left( a_{(l)}b \right)_{(\tilde{k})}c \right) t^{m+n+p-\tilde{k}-l} \\
&=& \sum_{l=0}^\infty \sum_{\tilde{k}=0}^\infty \binom{m}{l} \sum_{j=0}^\infty \binom{n}{j}\binom{m-l}{\tilde{k}-j} \left( \left( a_{(l)}b \right)_{(\tilde{k})}c \right) t^{m+n+p-\tilde{k}-l} \\
&=& \sum_{l=0}^\infty \sum_{\tilde{k}=0}^\infty \binom{m}{l} \binom{n+m-l}{\tilde{k}}\left( \left( a_{(l)}b \right)_{(\tilde{k})}c \right) t^{m+n+p-\tilde{k}-l} \\
&=& [[at^m,bt^n],ct^p]
\end{eqnarray*}
Therefore, (\ref{compeq1}) is equivalent to $[at^m,[bt^n,ct^p]=[[at^m,bt^n],ct^p + p(a,b)[bt^n,[at^m,ct^p]$. We then prove (3) by applying $\phi$ to both members.
\end{proof}
According to the last proposition, ($\tilde{{\cal R}}/\tilde{\partial}\tilde{{\cal R}}$,$[\, ,\, ]$) is the expected Lie algebra.  So $\text{Lie}({\cal R})=\tilde{{\cal R}}/\tilde{\partial}\tilde{{\cal R}}$. 
The family $\bar{{\cal R}}$ of $\text{Lie}({\cal R})$-valued formal distributions, whose coefficients span $\text{Lie}({\cal R})$ is defined by 
\begin{eqnarray} 
\bar{{\cal R}}= \{ \sum_{n\in\corps{Z}} \phi(at^n) z^{-n-1} \,|\, a\in{\cal R} \} \label{rbarrrr} 
\end{eqnarray}

The algebraic structure of (Lie($\cal R$),$[\, , \,]$) gives $\bar{{\cal R}}$ the following properties. 
\begin{prop}
Let $\cal R$ be a Lie conformal algebra and the space $\bar{{\cal R}}$ defined by \eqref{rbarrrr}.
\begin{enumerate}
\renewcommand{\labelenumi}{\emph{(\theenumi)}}
\item $\bar{{\cal R}}$ is a local family; 
\item $\bar{{\cal R}}$ is a conformal family.
\end{enumerate}
\end{prop}
\begin{proof}
Let's prove (1). Let 
$$\bar{a}(z)=\sum_{n\in\corps{Z}} \phi(at^n)z^{-n-1} \in \text{Lie}({\cal R})[[z^{\pm 1}]] \quad \text{et}\quad \bar{b}(w)=\sum_{m\in\corps{Z}} \phi(bt^m)w^{-m-1} \in \text{Lie}({\cal R})[[w^{\pm 1}]] \ .$$
(\ref{commre2}) implies
$$ [\bar{a}(z),\bar{b}(w)]=\sum_{j\in \corps{Z}^+} \phi\left( a_{(j)}b \right)(w) \frac{\partial_w^j \delta(z,w)}{j!} \ ,  $$
as already shown. $\bar{{\cal R}}$ is thus a local family.\\
Let's prove (2). $\bar{{\cal R}}$ is the \textit{minimal} local family since $\cal R$ is.  Moreover, since $\cal R$ is a $\corps{C}[\partial]$-module, $\partial a \in {\cal R}$ (for all $a\in {\cal R}$).  So $(\overline{\partial a})(z) \in \bar{{\cal R}}$.  Now, since $\phi((\partial a)t^n)=-n \phi(at^{n-1})$, we have
\begin{eqnarray*}
(\overline{\partial a})(z) &=& \sum_{n\in\corps{Z}} \phi\left((\partial a)t^n\right) z^{-n-1} \\
&=& \sum_{n\in\corps{Z}} -n \phi(at^{n-1}) z^{-n-1} \\
&=& \sum_{n\in\corps{Z}} -(\tilde{n}-1) \phi(at^{\tilde{n}}) z^{\tilde{n}-2} \quad \text{where $\tilde{n}=n-1$} \\
&=& \partial_z \bar{a}(z) \ ,
\end{eqnarray*}  
so that $\bar{{\cal R}}$ is indeed $\partial_z$-invariant.
\end{proof}

The first property shows that (Lie($\cal R$),$\bar{R}$) is the expected formal distribution Lie algebra.
By the second property, $\bar{{\cal R}}$ is a $\corps{C}[\bar{\partial}]$-module.  Since the latter can be endowed with a $\lambda$-bracket (resp. $j$-products) satisfying its (resp. their) defining relations, $\bar{{\cal R}}$ is endowed with a structure of Lie conformal algebra.  

Note that the spaces $\cal R$ and $\bar{{\cal R}}$ isomorphic by construction.  They are indeed mapped to one another by a bijective application $\psi$ defined by
$$ a \in {\cal R} \stackrel{\psi}{\rightarrow} \bar{a}(z)=\sum_{n\in\corps{Z}} \phi(at^n) z^{-n-1} \in \bar{{\cal R}} \ .$$ 
To conclude, we mention that there exist formal distributions Lie algebras, other than (Lie($\cal R$),$\bar{R}$), leading to corresponding \textit{isomorphic} Lie conformal algebras as well.  They are related to (Lie($\cal R$),$\bar{R}$) in a specific way.  Indeed, they appear as the quotient of the latter by an \textit{irregular} ideal, as shown in proposition \ref{questionposee}. First, we define this notion of \textit{irregular} ideal below.
\begin{defin} Let $({\mathfrak g},F)$ a formal distribution Lie algebra and $\bar{F}$ the minimal local family of $F$.  An ideal $I\subset{\mathfrak g}$ is said \emph{irregular}\index{irregular ideal} if and only if there are no non-vanishing elements $b(z)=\sum_n b_n z^n \in\bar{F}$ such that $b_n\in I \ \forall n\in\corps{Z}$. 
\end{defin}
In other words, an ideal $I\subset{\mathfrak g}$ is irregular if and only if the only element $b(z)\in\bar{F}$ whose coefficients $b_n$ lies in $I$ is $b(z)=0$. 
\begin{prop}\label{questionposee}
Let's consider the formal distribution Lie algebras $(\text{Lie}({\cal R}),\bar{{\cal R}})$ and $(\text{Lie}({\cal R})/I,\bar{{\cal R}}_I)$ where $I \subset\text{Lie}({\cal R})$ is an irregular ideal.  The homomorphism $\phi \,:\, \text{Lie}({\cal R}) \to \text{Lie}({\cal R})/I$ induces an isomorphism $\bar{\phi}\,:\,\bar{{\cal R}} \to \bar{{\cal R}}_I$.
\end{prop}
\begin{proof}
Surjectivity is clear, so let's show that if $I$ is an irregular ideal, $\bar{\phi}$ is injective and vice versa.  Consider $\text{Ker}\bar{\phi}=\{ b(z)\in \bar{{\cal R}} \, \mid \, \bar{\phi}(b(z))=0 \}$.  Let $b(z)\in\text{Ker}\bar{\phi}$. By definition, $\bar{\phi}(b(z))=0$ if and only if $b(z)\in I$, \textit{i.e.} $b_n \in I \ \forall n\in\corps{Z}$, where $b(z)=\sum_n b_n z^n$.  But since $I$ is irregular, $b_n \in I \ \forall n\in\corps{Z}$ if and only if  $b_n =0 \ \forall n\in\corps{Z}$, so that $b(z)=0$ and hence the injectivity.
\end{proof}
Since all formal distribution Lie algebras leading to corresponding isomorphic Lie conformal algebras arise as quotients of $(\text{Lie}({\cal R}),\bar{{\cal R}})$ by an irregular ideal, or as $(\text{Lie}({\cal R}),\bar{{\cal R}})$ itself, the latter is said \textit{maximal}\index{formal distributions Lie superalgebra!maximal}.  

Define two formal distributions Lie algebras as \textit{equivalent} if their corresponding Lie conformal algebras are isomorphic.  Each equivalence class of formal distributions Lie algebras is thus shown to be \textit{bijectively} mapped to the corresponding isomorphism class of Lie conformal algebras.  This result is referred to as the \textbf{canonical bijection}\index{canonical bijection}.

We previously constructed the maximal formal distribution Lie algebra $\text{Lie}({\cal R})$ associated to the Lie conformal algebra $\cal R$.  We will soon illustrate the connection between both notions on important examples (see next section). In those examples, $\cal R$ will be seen to have the structure of a $\corps{C}[\partial]$-module \textit{finitely} and \textit{freely} generated by the elements of a set $F=\{a^j,C ; j\in J \}$, where $J$ is an index set and $C$ such that $\partial C = 0$.  Before going any further, we define what is meant by a \textit{finitely} and \textit{freely} generated $R$-module\index{$R$-module} by a set $E$ in general.
\begin{defin}
Let $M$ be a $R$-module and $n<\infty$.\\A set $E=\{e_1,e_2,\ldots,e_n \}$ is a finite\index{$R$-module!finitely generated} and \textit{free}\index{$R$-module!freely generated} generating set of $M$ if
\begin{enumerate}
\renewcommand{\labelenumi}{(\theenumi)}
\item $E$ generates $M$, i.e. any element of $M$ is an $R$-linear combination of elements of $E$, 
\item $E$ is a free generating set, i.e. for $r_1,r_2,\ldots,r_n\in R$, 
$$r_1 e_1 + r_2 e_2 + \cdots + r_n e_n = 0 \Longrightarrow r_1=r_2=\cdots=r_n=0 \ . $$ 
\end{enumerate}
An $R$-module $M$ endowed with a finite and free generating set is said \textit{finitely and freely generated}.
\end{defin} 
We show in the next proposition that the knowledge of the generating set of $\cal R$ is equivalent to the knowledge of a basis of $\text{Lie}({\cal R})$.  
\begin{prop} The following statements are equivalent.
\begin{enumerate}
\renewcommand{\labelenumi}{\emph{(\theenumi)}}
\item ${\cal R}$ is finitely and freely generated as a $\corps{C}[\partial]$-module by the set $F=\{a^j,C ; j\in J \}$, where $J$ is an index set and $C$ such that $\partial C = 0$, i.e.
$${\cal R}=\bigoplus_{i\in I} \corps{C}[\partial]a^i + \corps{C}C \ \ , \ \partial C=0 \ .$$
\item $\{a^j_{(n)}, C_{(-1)} \mid n\in\corps{Z} \ , \ j\in J \}$ is a basis of Lie($\cal R$). 
\end{enumerate}
\end{prop}
\begin{proof}
Let's prove (1)$\Leftarrow$(2).
Since, on the one hand, $\{a^j_{(n)}, C_{(-1)} \mid n\in\corps{Z} \ , \ j\in J \}$ is a basis of Lie($\cal R$), and on the other hand ${\cal R}\cong\bar{\cal R}$, as previously shown, the elements of ${\cal R}$ are bijectively mapped to elements of the form $a^j(z)=\sum_{n\in\corps{Z}} a^j_{(n)}z^{-n-1}$, with $a^j_{(n)} \in \text{Lie}({\cal R})$, and $C=\sum_{n\in\corps{Z}} C_{(n)} z^{-n-1}$, with $C_{(n)}=C\delta_{-1,n}$. We will show that those elements finitely and freely generate ${\cal R}$ as a $\corps{C}[\partial]$-module, \textit{i.e.} by the action of differential operators of finite order or, formally, by multiplication -- in the sense of a $\corps{C}[\partial]$-module -- by polynomials in $\partial$ with constant complex coefficients, the action of $\partial$ on formal distributions being defined by $(\partial a)(z)=\partial_z a(z)$.  Let the polynomials be $p_j(\partial)=\sum_{k_j=0}^{K_j} \alpha_{k_j}^{j} \partial^{k_j}$, with $\alpha_{k_j}^{j}\in\corps{C}$ and $j\in J$.  We have to show the following implication : 
$$p_1(\partial) a^1(z) + \cdots + p_m(\partial) a^m (z) = 0 \Rightarrow p_j(\partial)=0 \ \forall j\in J \ .$$  
The equality $p_j(\partial)=0$ means that $\alpha_{k_j}^j=0 \ \forall k_j\in\corps{Z}\mid 0\leq k_j\leq K_j$.  Suppose there exist values $\alpha_{k_j}^{j}\neq 0$ such that $p_1(\partial) a^1(z) + \cdots + p_m(\partial) a^m (z) = 0$.  In that case, we would have
\begin{eqnarray}
\text{Res}_z z^l \left( p_1(\partial) a^1(z) + \cdots + p_m(\partial) a^m (z)\right) = 0 \quad \forall l\in\corps{Z} . \label{q32e6}
\end{eqnarray} 
But, $(\partial^k a)_{(l)}=\text{Res}_z z^l \partial^k a(z)=(-1)^k \binom{l}{k}k! a_{(l-k)}$, $\forall l\in\corps{Z}$. The expression (\ref{q32e6}) then reads
\begin{eqnarray} 
\sum_{j=1}^m \sum_{k_j=0}^{K_j} \alpha_{k_j}^{j} (-1)^{k_j} \binom{l}{k_j}k_j! a_{(l-k_j)} = 0 \quad \forall l\in\corps{Z} \label{eold}
\end{eqnarray}
Since $\alpha_{k_j}^{j}\neq 0$, this shows that $a^j_{(n)}$ are linearly dependant, hence a contradiction.  Therefore, ${\cal R}$ is indeed freely generated by the elements $C=C_{(-1)}$ and $a^j (z)=\sum_{n\in\corps{Z}} a^j_{(n)} z^{-n-1}$ ($j\in J$) as a $\corps{C}[\partial]$-module. We can then write ${\cal R} = \bigoplus_{j\in J} \corps{C}[\partial] a^j + \corps{C}C$. Moreover, we clearly have $\partial C=0$.

Let's prove (1)$\Rightarrow$(2). By assumption, ${\cal R} = \bigoplus_{j\in J} \corps{C}[\partial] a^j + \corps{C}C$ and $\partial C=0$,
which means $\lbrace a^j,\, C \ ; j\in J \rbrace$ is a generating set of $\cal R$ as a $\corps{C}[\partial]$-module, and $C$ is the only torsion element.  We have to show that $\lbrace a^j_{(n)},\, C_{(-1)} \ ; n\in\corps{Z} \ , j\in J \rbrace$ is then a basis of $\text{Lie}({\cal R})$.  We use the notation $a^j_{(n)}\doteq \phi(a^j t^n)$, where $\phi\,:\,\tilde{\cal R} \to \text{Lie}({\cal R})$.  For clarity, we replace $\phi(a^j t^n)$ by $a^j t^n$.  The elements $a^j_{(n)}$ ($\forall n\in\corps{Z}$, $\forall j\in J$) and $C_{(-1)}$ generate $\text{Lie}({\cal R})$ by the action of differential operators of finite order -- or ``polynomial'' in $\corps{C}[\partial]$ --
since, in $\text{Lie}({\cal R})$, we have $(\partial a^j)_{(m)}=(\partial a^j)t^m=-m a^j t^{m-1}=-m a^j_{(m-1)}$, hence
$(\partial^k a^j)_{(m)} = (\partial^k a^j)t^m = (-1)^k \binom{m}{k} k! a^j_{(m-k)}$. 
Now, since $(\partial C)t^n + n C t^{n-1} = 0$ in $\text{Lie}({\cal R})$ and $\partial C=0$ by assumption, we deduce that $C_{(n-1)}=0 \ \forall n \in\corps{Z}_0$, so that $C_{(-1)}$ alone is non vanishing.  
Therefore, $\lbrace a^j_{(n)},\, C_{(-1)} \ ; n\in\corps{Z} \ , j\in J \rbrace$ is a basis of $\text{Lie}({\cal R})$.  It is moreover a \textit{free} basis, as proved by contradiction. Suppose there exist $\alpha^j_{k_j}\neq 0$ such that $\sum_{j,k_j}\alpha^j_{k_j} a_{(k_j)}^j=0$.  Such an expression can be written as (\ref{eold}) and therefore as (\ref{q32e6}). 
Since $\alpha^j_{k_j}\neq 0$, $\cal R$ would not be freely generated as a $\corps{C}[\partial]$-module by the set $\lbrace a^j,\, C \rbrace$, which would lead to a contradiction.
\end{proof}

\section{Examples}
This section illustrates on relevant examples in Theoretical Physics the construction of a formal distributions Lie superalgebra and the associated Lie conformal algebra.

\subsection{Virasoro algebra}
The Virasoro algebra\index{Virasoro algebra}, denoted by $\mathfrak{Vir}$, is an infinite-dimensional complex Lie algebra. A basis is given by the set $\{ L_m, \, C \, ; m\in\corps{Z} \}$, where $C$ is the central element (i.e. commuting with the other generators).  Those elements satisfy the following commutation relation.
\begin{eqnarray}
[L_m,L_n]=(m-n)L_{m+n} + \delta_{m,-n} \frac{m^3-m}{12} C \quad \text{and} \quad [L_m,C]=0 \quad \forall m,n\in\corps{Z} \label{vira} \ .
\end{eqnarray}
The central element put aside, we obtain the complex Lie algebra of derivations of the space $\corps{C}[z,z^{-1}]$, called the \textit{Witt algebra}, generated by the elements
$$ d_n = - z^{n+1} \frac{d}{dz} \ , $$
where $n\in\corps{Z}$.

Starting from $\mathfrak{Vir}$, we construct $\mathfrak{Vir}$-valued formal distributions by setting conventionally
$$L(z)=\sum_{n\in\corps{Z}} L_{(n)} z^{-n-1} \quad \text{with $L_{(n)}\doteq L_{n-1}$.}$$
Hence, $L(z)=\sum_{n\in\corps{Z}} L_{n} z^{-n-2}$.  This convention will acquire a certain meaning when we will treat the notion of \textit{weight} of an \textit{eigendistribution}. 

Now, we translate the commutation relations of the Virasoro Algebra in terms of the formal distributions previously defined. 
\begin{prop}
The commutation relations (\ref{vira}) are equivalent to 
\begin{eqnarray}
[L(z),L(w)]=\partial_w L(w) \delta(z,w) + 2 L(w) \partial_w \delta(z,w) + \frac{C}{12} \partial_w^3 \delta(z,w) \ . \label{viradis}
\end{eqnarray}
As $[L(z),C]=0$, $\{ L(z),C \}$ forms a local family of $\mathfrak{Vir}$-valued formal distributions.
\end{prop}
\begin{proof}
Using (\ref{vira}), we compute
\begin{eqnarray}
[L(z),L(w)] &=& \sum_{n,m\in\corps{Z}} [L_{(n)},L_{(m)}] z^{-n-1}w^{-m-1} \nonumber\\
&=& \sum_{n,m\in\corps{Z}} [L_{n-1},L_{m-1}] z^{-n-1}w^{-m-1} \nonumber\\
&=& \sum_{n,m\in\corps{Z}} \left( (n-m)L_{n+m-2} + \delta_{n+m,2} \frac{n(n-1)(n-2)}{12}C \right) z^{-n-1}w^{-m-1}\label{ekk}
\end{eqnarray}
The first term of the RHS reads
{\setlength\arraycolsep{0pt}
\begin{eqnarray*}
\sum_{n,m\in\corps{Z}} (n-m)&&L_{n+m-2}z^{-n-1}w^{-m-1} \\
&&=\; \sum_{n,m\in\corps{Z}} (-m) L_{n+m-2} z^{-1-n} w^{-m-1} + \sum_{n,m\in\corps{Z}} n L_{n+m-2} z^{-1-n} w^{-m-1} \\
&&=\; \partial_w \left( \sum_{n,m\in\corps{Z}} L_{n+m-2} z^{-1-n} w^{-m} \right) + \sum_{l\in\corps{Z}} L_{l-1} w^{-l-1} \sum_{n\in\corps{Z}} n z^{-1-n}w^{n-1}  \\
&&=\; \partial_w \left( \sum_{k\in\corps{Z}} L_{k-1} w^{-k-1} \sum_{n\in\corps{Z}} z^{-1-n} w^{n} \right) + L(w) \partial_w \delta(z,w) \\
&&=\; \partial_w \left( L(w) \delta(z,w) \right) + L(w) \partial_w \delta(z,w) \;=\; \partial_w L(w) \delta(z,w) + 2L(w)\partial_w \delta(z,w)
\end{eqnarray*}}
where $m=l-n+1$ at the second line and $m=k-n+1$ at the third.\\
Considering (\ref{ekk}) again, the second term of the RHS reads
{\setlength\arraycolsep{0pt}
\begin{eqnarray*}
\sum_{n,m\in\corps{Z}} \delta_{n+m,2} &&\frac{n(n-1)(n-2)}{12}C z^{-n-1}w^{-m-1} \\
&&=\; \sum_{n\in\corps{Z}} \frac{n(n-1)(n-2)}{12}C z^{-n-1} w^{n-3} \\
&&=\; \frac{C}{12}\partial_w^3 \sum_{n\in\corps{Z}}z^{-n-1} w^n \;=\; \frac{C}{12}\partial_w^3 \delta(z,w) \ .
\end{eqnarray*}}
So (\ref{viradis}) is verified. Since the RHS is a finite sum of $\delta$'s and its derivatives (weighted by formal distributions in one indeterminate), $L(z)$ is a local formal distribution with respect to itself.  As $[L(z),C]=0$, the set $\{ L(z),C \}$ is indeed a local family, hence the proposition.
\end{proof}
\begin{prop}\label{s6c59}
In term of the $\lambda$-bracket, (\ref{viradis}) translates as follows :
\begin{eqnarray}
[L_\lambda L]=(\partial + 2\lambda) L  + \frac{\lambda^3}{12}C \ , \label{lprodvira}
\end{eqnarray}
In terms of $j$-products, we get
$$L_{(0)}L=\partial L,\quad  L_{(1)}L=2 L, \quad L_{(3)}L=\frac{C}{2}, \quad L_{(j)}L=0 \ \forall j\neq\{0,1,3\},\quad L_{(j)}C=0 \ \forall j\in\corps{Z} \ .$$
\end{prop}
\begin{proof}
Starting from (\ref{viradis}), we compute
\begin{eqnarray*}
[L_\lambda L] &\doteq & \Fz [L(z),L(w)] \\
&=& \text{Res}_z e^{\lambda(z-w)} \left( \partial_w L(w) \delta(z,w) + 2 L(w) \partial_w \delta(z,w) + \frac{C}{12} \partial_w^3 \delta(z,w) \right) 
\end{eqnarray*}
The first term reads
\begin{eqnarray*}
\text{Res}_z e^{\lambda(z-w)}  \partial_w L(w) \delta(z,w) &=& \partial_w L(w) \text{Res}_z \sum_{j=0}^\infty \frac{\lambda^j}{j!} (z-w)^j \delta(z,w) \\
&=& \partial_w L(w) \sum_{j=0}^\infty \frac{\lambda^j}{j!} \text{Res}_z \underbrace{(z-w)^j \delta(z,w)}_{=0 \; \forall j>0} \;=\; \partial_w L(w)
\end{eqnarray*}
The second reads
\begin{eqnarray*}
\text{Res}_z e^{\lambda(z-w)}  2 L(w) \partial_w \delta(z,w) &=& 2 L(w) \sum_{j=0}^\infty \frac{\lambda^j}{j!} \text{Res}_z \underbrace{(z-w)^j \partial_w \delta(z,w)}_{=0 \; \forall j>1} \\
&=& 2L(w) \left( \underbrace{\text{Res}_z \partial_w \delta(z,w)}_{=0} + \lambda \text{Res}_z \underbrace{(z-w) \partial_w \delta(z,w)}_{=\delta(z,w)} \right) \\
&=& 2\lambda\, L(w)
\end{eqnarray*}
At last, the third reads
\begin{eqnarray*}
\text{Res}_z e^{\lambda(z-w)}\frac{C}{12} \partial_w^3 \delta(z,w) &=& \frac{C}{12} \sum_{j=0}^\infty \frac{\lambda^j}{j!} \text{Res}_z (z-w)^j \partial_w^3 \delta(z,w) \\
&=& \frac{C}{12} \sum_{j=0}^\infty \frac{\lambda^j}{j!} 3! \text{Res}_z \frac{(z-w)^j\partial_w^3 \delta(z,w)}{3!} \\
&=& \frac{C}{12} \sum_{j=0}^\infty \frac{\lambda^j}{j!} 3! \text{Res}_z \frac{\partial_w^{3-j}\delta(z,w)}{(3-j)!} \\
&=& \frac{C}{12} \sum_{j=0}^\infty \frac{\lambda^j}{j!} 3! \delta_{3,j} \\
&=& \frac{C}{12} \lambda^3 \ ,
\end{eqnarray*}
which proves (\ref{lprodvira}).  The translation in terms of $j$-products is immediate by (\ref{viradis}).  This can be done in another way : since $[L_\lambda L]\doteq \sum_{j=0}^\infty \frac{\lambda^j}{j!} \left( L_{(j)}L\right)$, we have
$[L_\lambda L] = \frac{\lambda^0}{0!} \partial_w L(w) + \frac{\lambda^1}{1!} 2L(w) + \frac{\lambda^3}{3!} \frac{3!C}{12}$.
Proceeding by identification, we get the expected result.  Now, $[L(z),C]=0$ implies $L_{(j)}C=0$ for $j\in\corps{Z}$.
\end{proof}
The minimal local family $\bar{F}$ is just $\corps{C}[\partial]F$.  It contains $F$ and, according to \eqref{cmod1} and \eqref{cmod2}, is closed under $j$-product. The conformal family $\cal R$ is in this case obtained as a $\corps{C}[\partial]$-module generated by $F$ itself, i.e. ${\cal R}=\corps{C}[\partial]F$.  This reads ${\cal R}=\corps{C}[\partial]L + \corps{C}C$, with $\partial C=0$ and $\partial\equiv\partial_z$.  Endowed with the $\lambda$-bracket of the previous proposition, we get the Lie conformal algebra\index{Virasoro Lie conformal algebra} by
$$\text{Vir} = \corps{C}[\partial]L + \corps{C}C, \quad \partial C=0 \ ,$$ as verified in the next proposition.
\begin{prop}\label{virrr}
Vir$=\corps{C}[\partial]L + \corps{C}C$, with $\partial C=0$ and endowed with the $\lambda$-bracket
$$[L_\lambda L]=(\partial + 2\lambda) L  + \frac{\lambda^3}{12}C \quad\text{et}\quad [L_\lambda C]=0$$ 
is a Lie conformal algebra.
\end{prop}
\begin{proof}
Vir is indeed a Lie conformal algebra.  Four properties have to be satisfied (see the definition).  For the first two, we  just note that the $\lambda$-bracket was defined from formal distributions in $\mathfrak{Vir}[[z^{\pm 1}]]$, which is a $\corps{C}[\partial]$-module and had precisely permitted us to introduce those two properties.  Now the last two properties.  The antisymmetry is satisfied : as $\partial C=0$,
\begin{eqnarray*}
[L_\lambda L]=(\partial + 2\lambda) L + \frac{\lambda^3}{12} C &=& -\left( \partial + 2(-\lambda - \partial) \right) L - \frac{(-\lambda - \partial)^3}{12}C \\
&=& -[L_{-\lambda - \partial} L] \ .
\end{eqnarray*}
At last, Jacobi is true as well ($p(L,L)=1$).  We compute separately
\begin{eqnarray}
\lbrack L_\lambda \lbrack L_\mu L\rbrack \rbrack  &=& \lbrack L_\lambda \partial L\rbrack  + 2\mu \lbrack L_\lambda L\rbrack  \nonumber\\
&=&(\partial + \lambda + 2\mu)(\partial + 2\lambda) L \nonumber \\
&=& (\partial^2 + (3\lambda+2\mu)\partial + 2\lambda^2 + 4\mu \lambda) L \label{j0801} \\
\lbrack \lbrack L_\lambda L\rbrack _{\lambda + \mu} L\rbrack  &=& \lbrack \partial L_{\lambda+\mu}L\rbrack  + 2\lambda \lbrack L_{\lambda + \mu}L\rbrack  \nonumber\\
&=& (\lambda - \mu) (\partial + 2(\lambda+\mu))L \nonumber \\
&=& ((\lambda-\mu)\partial + 2\lambda^2 -2 \mu^2) L \label{j0802} \\
\lbrack L_\mu \lbrack L_\lambda L\rbrack \rbrack  &=& (\partial^2 + (3\mu+2\lambda)\partial + 2\mu^2 + 4\mu \lambda) L \label{j0803}
\end{eqnarray}
Subtracting (\ref{j0801}) by (\ref{j0802}), we get (\ref{j0803}), hence $$\lbrack L_\lambda \lbrack L_\mu L\rbrack \rbrack=\lbrack \lbrack L_\lambda L\rbrack _{\lambda + \mu} L\rbrack + p(L,L)\lbrack L_\mu \lbrack L_\lambda L\rbrack \rbrack \ .$$
\end{proof}

The Virasoro Algebra can be recovered from the Lie conformal algebra.
\begin{prop} The Lie conformal algebra $\text{Vir}$ encodes the algebraic structure of $(\mathfrak{Vir},\text{Vir})$.  (\ref{vira}) is equivalent to
$$ [L_{(m)},L_{(n)}] = \sum_{j\in\corps{Z}^+} \binom{m}{j} \left( L_{(j)}L \right)_{(m+n-j)} \ .$$
\end{prop}
\begin{proof}
By \ref{s6c59}, we have
\begin{eqnarray}
[L_{(m)},L_{(n)}] &=& \sum_{j\in\corps{Z}^+} \binom{m}{j} \left( L_{(j)}L \right)_{(m+n-j)} \nonumber\\
&=& \binom{m}{0} \left( L_{(0)}L \right)_{(m+n)} + \binom{m}{1} \left( L_{(1)}L \right)_{(m+n-1)} + \binom{m}{3} \left( L_{(3)}L \right)_{(m+n-3)} \nonumber\\
&=& \left( \partial L \right)_{(m+n)} + 2m L_{(m+n-1)} + \frac{m(m-1)(m-2)}{6} \left( \frac{1}{2}C_{(m+n-3)} \right) \ . \label{qsdf64}
\end{eqnarray}
But, $(\partial L)_{(n)}=-n L_{(n-1)}$ and $C_{(m+n-3)}=\delta_{m+n-3,-1} C = \delta_{m+n,2} C$. Thus (\ref{qsdf64}) becomes
\begin{eqnarray}
[L_{(m)},L_{(n)}] &=& \left( \partial L \right)_{(m+n)} + 2m L_{(m+n-1)} + \frac{m(m-1)(m-2)}{6} \left( \frac{1}{2}C_{(m+n-3)} \right) \nonumber \\
&=& -(n+m) L_{(n+m-1)} + 2m L_{(m+n-1)} + \frac{m(m-1)(m-2)}{12} C\delta_{m+n,2} \nonumber\\
&=& (m-n)L_{(n+m-1)} + \delta_{m+n,2} \frac{m(m-1)(m-2)}{12} C \ . \label{qf94}
\end{eqnarray}
As $L_{(n)}=L_{n-1}$, we just have to make the shifts $m\to m+1$ and $n\to n+1$ in (\ref{qf94}) to recover
$$[L_m,L_n] = (m-n) L_{n+m} + \delta_{m,-n} \frac{m^3-m}{12} C \ ,$$
which concludes the proof.
\end{proof}

\subsection{Neveu-Schwarz superalgebra}
The Neveu-Schwarz Lie superalgebra\index{Neveu-Schwarz Lie superalgebra}, denoted by $\mathfrak{ns}$, is closely related to the Virasoro algebra, as it appears as a supersymmetric extension of the latter.  We begin by defining its Lie conformal superalgebra and the associated formal distributions Lie superalgebra.  We then derive the commutation relations satisfied by their Fourier modes.

The Neveu-Schwarz Lie Conformal algebra\index{Neveu-Schwarz Lie conformal algebra}, denoted by NS, is the $\corps{C}[\partial]$-module 
$$\text{NS}=\corps{C}[\partial]L \bigoplus \corps{C}[\partial]G \bigoplus \corps{C}C \ , \ \ \partial C=0 \ ,$$ 
where $C$ is a central element.  The elements $L$ and $C$ are even, but $G$ is odd.  This $\corps{C}[\partial]$-module is endowed with the $\lambda$-brackets
\begin{eqnarray}
\lbrack L_\lambda L\rbrack  &=& (\partial + 2\lambda) L + \frac{\lambda^3}{12}C \nonumber\\
\lbrack L_\lambda G\rbrack  &=& \left( \partial + \frac{3}{2} \lambda \right) G \nonumber\\
\lbrack G_\lambda L\rbrack  &=& \left( \frac{1}{2}\partial + \frac{3}{2} \lambda \right) G \label{nslambda}\\
\lbrack G_\lambda G\rbrack  &=& L + \frac{\lambda^2}{6}C \nonumber\\
\lbrack C_\lambda \,\cdot\,\rbrack  &=& 0 \ .\nonumber
\end{eqnarray}

\begin{prop}
$\text{NS}=\corps{C}[\partial]L \bigoplus \corps{C}[\partial]G \bigoplus \corps{C}C$, $\partial C=0$, endowed with the $\lambda$-brackets (\ref{nslambda}), is a Lie conformal superalgebra.
\end{prop}
\begin{proof}
As the $\lambda$-brackets are defined on generators of a $\corps{C}[\partial]$-module, the first two conditions are automatically satisfied.  Let's verify the anticommutativity.  As $\partial C=0$, $p(L,L)=p(G,L)=1$ and $p(G,G)=-1$, we can write
\begin{eqnarray*}
\lbrack L_\lambda L\rbrack  &=& (\partial + 2\lambda) L + \frac{\lambda^3}{12}C = -(\partial + 2(-\lambda - \partial))L - \frac{(-\lambda - \partial)^3}{12}C \\ 
&=& -p(L,L)\lbrack L_{-\lambda - \partial}L\rbrack  \ , \\
\lbrack G_\lambda G\rbrack  &=& L+ \frac{\lambda^2}{6}C = L + \frac{(-\lambda - \partial)^2}{6} C \\
&=& -p(G,G)\lbrack G_{-\lambda - \partial} G\rbrack  \ , \\
\lbrack G_\lambda L\rbrack  &=& \left( \frac{1}{2}\partial + \frac{3}{2} \lambda \right) G = -\left(\partial +  \frac{3}{2}(-\lambda-\partial) \right) G \\
&=& -p(G,L)\lbrack L_{-\lambda-\partial} G\rbrack 
\end{eqnarray*}

Now the Jacobi identity.  As $\lbrack C_\lambda \cdot\rbrack = 0$, there are eight identities to be checked, mixing the elements $L$ and $G$. It can be shown that it suffices to verify only four of them. Without going into much details, this results from the fact that the structures of Lie conformal algebra and formal distributions Lie algebra are in bijection.

As $p(G,G)=-1$,
\begin{eqnarray}
\lbrack G_{\lambda} \lbrack G_{\mu}  G\rbrack \rbrack = \lbrack \lbrack G_{\lambda} G \rbrack_{\lambda + \mu}  G\rbrack  - \lbrack G_{\mu} \lbrack G_{\lambda}  G\rbrack \rbrack \ . \label{permut1}
\end{eqnarray}  
Since $[G_\lambda G]=L + \frac{\lambda^2}{6}C$ and $[C_{\lambda}\,\cdot\,]=[\,\cdot\,_{\lambda}C]=0$, (\ref{permut1}) is equivalent to 
\begin{eqnarray*}
[G_\lambda L] = [L_{\lambda+\mu} G] - [G_\mu L] \label{permut2} \ .
\end{eqnarray*}  
This equality is clearly satisfied, since by (\ref{nslambda}), it is equivalent to 
$$ \frac{\partial + 3\lambda}{2} G = (\partial + \frac{3}{2}(\lambda + \mu))G - \frac{1}{2}(\partial + 3\mu)G \ . $$

As $p(L,L)=1$,
\begin{eqnarray*}
\lbrack L_\lambda\lbrack L_\mu G\rbrack \rbrack  - \lbrack L_\mu\lbrack L_\lambda G\rbrack \rbrack  &=& (  \partial + \lambda + \frac{3}{2} \mu ) \lbrack L_\lambda G\rbrack  -   (  \partial + \mu + \frac{3}{2} \lambda ) \lbrack L_\mu G\rbrack  \quad \text{by (\ref{cmod02})} \\
&=& (\lambda - \mu) \partial G + \frac{3}{2}(\lambda^2 - \mu^2) \\
&=& (\lambda - \mu)( \partial + \frac{3}{2}(\lambda + \mu)) G \\
&=& (\lambda - \mu) \lbrack L_{\lambda+\mu}G\rbrack  \\
&=& -(\lambda+\mu)\lbrack L_{\lambda+\mu}G\rbrack  + 2\lambda \lbrack L_{\lambda+\mu}G\rbrack  \\
&=& \lbrack (\partial + 2\lambda) L_{\lambda + \mu} G\rbrack  \quad \text{by (\ref{cmod01})} \\
&=& \lbrack \lbrack L_\lambda L\rbrack _{\lambda + \mu} G\rbrack \ .
\end{eqnarray*}

As $p(L,G)=1$,
\begin{eqnarray*}
\lbrack \lbrack L_\lambda G\rbrack _{\lambda + \mu}G\rbrack  + \lbrack G_\mu \lbrack L_\lambda G\rbrack \rbrack  &=& \lbrack ( \partial + \frac{3}{2}\lambda )G_{\lambda + \mu}G\rbrack  + \lbrack G_\mu (  \partial + \frac{3}{2}\lambda )G\rbrack  \\
&=& (  \frac{\lambda}{2} - \mu ) \lbrack G_{\lambda + \mu}G\rbrack  + (  \partial + \frac{3}{2}\lambda + \mu ) \lbrack G_\mu G\rbrack  \\
&=& (  \partial + 2\lambda )L + \frac{\lambda^3}{12}C \\
&=& \lbrack L_\lambda L\rbrack \\
&=& \lbrack L_\lambda\lbrack G_\mu G\rbrack \rbrack  \quad \text{since $\lbrack \cdot_{\lambda}C\rbrack =0$} \ .
\end{eqnarray*}
Finally, the last one, containing only $L$, has already been proved in the Virasoro case.
\end{proof}

\begin{prop}
In terms of $\mathfrak{ns}$-valued formal distributions, the non trivial commutation relations of (\ref{nslambda}) translate as follows :
\begin{eqnarray*}
\lbrack L(z),L(w)\rbrack  &=& \partial_w L(w) \delta(z,w) + 2L(w) \partial_w \delta(z,w) + \frac{C}{12}\partial_w^3 \delta(z,w) \\
\lbrack L(z),G(w)\rbrack  &=& \partial_w G(w) \partial(z,w) + \frac{3}{2} G(w) \partial_w \delta(z,w) \\
\lbrack G(z),L(w)\rbrack  &=& \frac{1}{2}\partial_w G(w) + \frac{3}{2} G(w) \delta(z,w) \\
\lbrack G(z),G(w)\rbrack  &=& L(w)\partial(z,w) + \frac{C}{6} \partial_w^2 \delta(z,w)
\end{eqnarray*}
\end{prop}
\begin{proof}
Starting from the $\lambda$-brackets, the $j$-products are directly derived. We then use
\begin{eqnarray*}
\lbrack a(z),b(w)\rbrack =\sum_{j\in\corps{Z}^+} (a(w)_{(j)} b(w)) \frac{\partial_w^j \delta(z,w)}{j!} \ ,
\end{eqnarray*}
with $a,b \in \{ L,G,C \}$.
\end{proof}

The Fourier modes satisfy the following commutation relations.
\begin{prop}
Setting
\begin{eqnarray}
L(z) &=& \sum_{n\in\corps{Z}} L_n z^{-n-2} \label{mf0001}\\
G(z) &=& \sum_{n\in\frac{1}{2}+\corps{Z}} G_n z^{-n-3/2}\label{mf0002} \ ,
\end{eqnarray}
the non trivial commutation relations of (\ref{nslambda}) translate as follows :
\begin{eqnarray}
\lbrack L_m,L_n\rbrack  &=& (m-n) L_{m+n} + \delta_{n,-m} \frac{m^3-m}{12}C \\
\lbrack L_m,G_n\rbrack  &=& \left( \frac{m}{2}-n \right) G_{m+n} \\
\lbrack G_m,G_n\rbrack  &=& L_{m+n} + \frac{C}{6} \left( m^2 - \frac{1}{4} \right) \delta_{n,-m}
\end{eqnarray}
\end{prop}
\begin{proof}
By definition, we have $L_{(n)}=L_{n-1}$ and $G_{(n)}=G_{n-\frac{1}{2}}$.
The first relation has already been proved in the Virasoro case.  For the second one, by proposition \ref{s6c59}, we have
\begin{eqnarray*}
\lbrack L_{(m)},G_{(n)}\rbrack  &=& \sum_{j\in\corps{Z}^+} \binom{m}{j} \left( L_{(j)}G \right)_{(m+n-j)} \,=\, (\partial G)_{(m+n)} + m \frac{3}{2} G_{(m+n-1)} \\
&=& -(n+m) G_{(m+n-1)} + m \frac{3}{2} G_{(m+n-1)} \,=\, \left( \frac{m}{2} - n \right) G_{(m+n-1)}
\end{eqnarray*}
Hence, 
$ \lbrack L_{m-1}, G_{n-\frac{1}{2}} \rbrack = \left( \frac{m}{2} - n \right) G_{m+n-\frac{3}{2}}$.
By the shifts $m\to m+1$ and $n\to n+\frac{1}{2}$, we get
$$ \lbrack L_m,G_n\rbrack  = \left( \frac{m}{2}-n \right) G_{m+n} \ ,$$ 
which proves the second relation.
For the third one,
\begin{eqnarray*}
[G_{(m)},G_{(n)}] &=& \sum_{j\in\corps{Z}^+} \binom{m}{j} \left( G_{(j)}G \right)_{(m+n-j)} \\
&=& \binom{m}{0} L_{(m+n)} + \binom{m}{2} \frac{1}{3} C_{(m+n-2)} \\
&=& L_{(m+n)} + \frac{m(m-1)}{6} \delta_{m+n-2,-1} C 
\end{eqnarray*}
Hence,
$ [G_{m-\frac{1}{2}}, G_{n-\frac{1}{2}}] = L_{m+n-1} + \frac{m(m-1)}{6} \delta_{m-1,-n} C$.
By the shifts $m\to m+\frac{1}{2}$ and $n\to n+\frac{1}{2}$, we finally get
$ \lbrack G_m,G_n\rbrack  = L_{m+n} + \frac{C}{6} \left( m^2 - \frac{1}{4} \right) \delta_{n,-m} \ . $ 
\end{proof}

\subsection{Algebras of free bosons and free fermions}

We start by briefly defining the notion of \textit{loop algebra}\index{loop algebra}.
The loop algebra is the Lie algebra of the loop group $LG$\index{loop group}, the group of continuous maps from the circle to a Lie group $G$, its product being given the pointwise composition. We start by introducing the notion of invariant supersymmetric bilinear form\index{supersymmetric bilinear form}.
\begin{defin}[Supersymmetric bilinear form] 
Let $V=V_{\bar{0}} \oplus V_{\bar{1}}$ be a superspace.  A \textit{supersymmetric} bilinear form $(\cdot | \cdot)$ is map $(\cdot | \cdot)\,:\,V\otimes V\to \corps{C}$ such that $(a |  b)=(-1)^{p(a)}(b |  a)$, where $a,b$ are homogeneous elements of $V$.
\end{defin}
A supersymmetric bilinear form admits the following straightforward properties.
\begin{prop}
A supersymmetric bilinear form vanishes on $V_{\bar{0}} \oplus V_{\bar{1}}$ and $V_{\bar{1}} \oplus V_{\bar{0}}$, symmetric on $V_{\bar{0}} \oplus V_{\bar{0}}$ and antisymmetric on $V_{\bar{1}} \oplus V_{\bar{1}}$.
\end{prop}
\begin{proof}
Let $a,b\in V$. Since $(a |  b)=(-1)^{p(a)}(b | a)= (-1)^{p(a)+p(b)}(a |  b)$, we have $(a |  b)=0$ when $a$ and $b$ have opposite parities.  The remainder of the proposition is easily proved.
\end{proof}
\begin{defin}[Invariant bilinear form]
Let $V=V_{\bar{0}} \oplus V_{\bar{1}}$ be a superspace and carrying a representation of a Lie algebra $\mathfrak{g}$.   A bilinear form $(\cdot | \cdot)$ is said \emph{invariant}\index{invariant bilinear form} if, for any homogeneous element $X \in\mathfrak{g}$ and $v,w \in V$,
$(Xv|w) + p(v,X)(v|Xw) = 0$.
Moreover, if $V$ is the representation space of the adjoint representation, then the bilinear form is said \textit{adjoint} invariant\index{invariant bilinear form!adjoint invariant}.
\end{defin}
\begin{defin}
Let $\mathfrak{g}$ be a Lie superalgebra endowed with an adjoint invariant bilinear form $(\cdot | \cdot)$.  The associated loop algebra, denoted by $\tilde{\mathfrak{g}}$, is the algebra $\tilde{\mathfrak{g}}\doteq \mathfrak{g}\oplus\corps{C}[t,t^{-1}]$, endowed with the superbracket satisfying
$[ a\otimes f(t) , b\otimes g(t) ] = [a,b]\otimes f(t) g(t)$, where $a,b\in\mathfrak{g}$ and $f(t),g(t)\in \corps{C}[t,t^{-1}]$.
\end{defin}
Clearly, $\mathfrak{g}$ induces on $\tilde{\mathfrak{g}}$ a Lie superalgebra structure.
The commutation relations can be defined on the elements $f(t)=t^n$ and $g(t)=t^m$, with $n,m \in \corps{Z}$.  Then, replacing $a\otimes t^n$ by $at^n$ for clarity, those commutation relations read $[at^n,bt^m]=[a,b]t^{n+m}$.

We now define the notion of \textit{affinization}\index{affinization of a Lie algebra} of a Lie superalgebra.  But first, the central extension\index{central extension} of the loop algebra, denoted by $\hat{\mathfrak{g}}$, is the algebra $\hat{\mathfrak{g}}\doteq \tilde{\mathfrak{g}} \oplus \corps{C}K$, with the superbracket 
\begin{eqnarray} [at^n,bt^m]=[a,b]t^{n+m} + m\delta_{m,-n} (a|b) K \ , \label{affcrochet} \end{eqnarray}
where $a,b\in\mathfrak{g}$, $n,m\in \corps{Z}$ and $K$ is a central element\index{central element}.
Starting from a Lie superalgebra $\mathfrak{g}$, the central extension of the associated loop algebra $\tilde{\mathfrak{g}}$, by a one-dimensional center $\corps{C}K$, will be called the \textit{affinization} of $\mathfrak{g}$ and denoted by $\hat{\mathfrak{g}}$. The latter is also called the \textit{current algebra}\index{current algebra}. The superbracket (\ref{affcrochet}) gives $\hat{\mathfrak{g}}$ the structure of Lie superalgebra, as it can easily be shown.
If $\mathfrak{g}$ is a finite-dimensional simple Lie superalgebra, then the affinization of $\mathfrak{g}$ leads to a Kac-Moody algebra\index{Kac-Moody algebra}.

Starting from the current algebra $\hat{\mathfrak{g}}=\mathfrak{g}[t,t^{-1}] \oplus \corps{C}K$, where $\mathfrak{g}[t,t^{-1}]=\mathfrak{g}\otimes \corps{C}[t,t^{-1}]=\tilde{\mathfrak{g}}$, we now derive its underlying Lie conformal structure.  The $\hat{\mathfrak{g}}$-valued formal distributions are themselves called \textit{currents} and defined, for $a\in\mathfrak{g}$ by $a(z)=\sum_{n\in\corps{Z}} at^n z^{-n-1}$.  In terms of currents, the relation (\ref{affcrochet}) translates as follows. 
\begin{prop}
Setting $a(z)=\sum_{n\in\corps{Z}} at^n z^{-n-1}$, the commutation relation (\ref{affcrochet}) translates as
$[a(z),b(w)]=[a,b](w)\delta(z,w) + K (a|b) \partial_w \delta(z,w)$,  where $a,b\in\mathfrak{g}$.
\end{prop}
\begin{proof}
We just compute :
\begin{eqnarray*}
[a(z),b(w)] &=& \sum_{m,n\in\corps{Z}} [at^m,at^n] z^{-1-m}w^{-1-n} \\
&=& \sum_{m,n\in\corps{Z}} [a,b] t^{m+n} z^{-1-m} w^{-1-n} + \sum_{m,n\in\corps{Z}} K(a|b) m \delta_{m,-n} z^{-1-m}w^{-1-n} \\
&=& \sum_{k\in\corps{Z}} [a,b] t^{k} z^{-1-k} \sum_{m\in\corps{Z}} w^{m} z^{-1-m} + \sum_{m\in\corps{Z}} K(a|b) m z^{-1-m}w^{m-1} \\ 
&=& [a,b](w) \delta(z,w) + K(a|b)\partial_w \delta(z,w)
\end{eqnarray*}
\end{proof}

\begin{prop}\label{propalgcour}
The structure of Lie conformal superalgebra, associated to the current algebra, is given by 
$$\text{Cur}\mathfrak{g}=\corps{C}[\partial]\mathfrak{g}\oplus \corps{C}K \ ,$$
endowed with the $\lambda$-bracket defined by
\begin{eqnarray*}
\lbrack  a_\lambda b\rbrack  &=& \lbrack a,b\rbrack  + (a|b)K \lambda \, \quad \text{where $a,b\in\mathfrak{g}$, } \\
\lbrack  K_\lambda a \rbrack  &=& 0 \quad \forall a\in\mathfrak{g} \ . 
\end{eqnarray*}
\end{prop} 
\begin{proof}
The translation in terms of the $\lambda$-bracket is straightforward.
\end{proof}

Let $\mathfrak{g}$ be an abelian Lie superalgebra. In that case, Cur$\mathfrak{g}$ is known as the \textbf{algebra of free bosons}\index{algebra of free bosons}, associated to the free bosons algebra $\hat{\mathfrak{g}}$, the latter being endowed with the relations $[at^m,bt^n]=m(a|b) \delta_{m,-n}K$.   

We now come to the notion of \textit{Clifford affinization}\index{affinization of a Lie algebra!Clifford} and its underlying conformal structure.
We previously defined the notion of supersymmetric bilinear form.  Now we define its \textit{antisupersymmetric} analogue.
\begin{defin}[Antisupersymmetric bilinear form]\index{antisupersymmetric bilinear form}
Let $V=V_{\bar{0}}\oplus V_{\bar{1}}$ be a superspace.  A bilinear form $<\cdot,\cdot>$ is \textit{antisupersymmetric} if it satisfies the relation $<a,b>=-(-1)^{p(a)}<b,a>$ for homogeneous elements $a,b\in V$.  
\end{defin}
It is easy to show that an antisymmetric bilinear form vanishes on $V_{\bar{0}}\oplus V_{\bar{1}}$ and $V_{\bar{1}}\oplus V_{\bar{0}}$, is antisymmetric on the even part $V_{\bar{0}}\oplus V_{\bar{0}}$ and symmetric on the odd part $V_{\bar{1}}\oplus V_{\bar{1}}$.
\begin{defin}
Let $\mathfrak{g}$ be a Lie superalgebra. The central extension $\hat{\mathfrak{g}}$ is the algebra $\hat{\mathfrak{g}}=\mathfrak{g}[t,t^{-1}]+\corps{C}K$, endowed with the following superbracket : 
\begin{eqnarray} [\varphi t^m, \psi t^n ] = K \delta_{m,-1-n} <\varphi,\psi> \label{affclifford} \ , \end{eqnarray} 
where $\varphi,\psi \in \mathfrak{g}$.  This is called the \textit{Clifford affinization} of $\mathfrak{g}$.
\end{defin}
It is not difficult to show that $\hat{\mathfrak{g}}$, endowed with the superbracket (\ref{affclifford}), is a Lie superalgebra.

Now, in terms of formal distributions, we have the following proposition.
\begin{prop}
Setting $\varphi(z)=\sum_{m\in\corps{Z}} \varphi t^m z^{-m-1}$, the relation (\ref{affclifford}) reads
$$ [\varphi(z),\psi(w)] = K <\varphi,\psi> \delta(z,w) \ .$$
\end{prop}
\begin{proof}
The proof is straightforward.
\end{proof}
The Lie conformal superalgebra associated to the Clifford affinization is called the \textbf{algebra of free fermions}\index{algebra of free fermions} and defined in the next proposition.
\begin{prop}\label{ferlibr}
The structure of Lie conformal superalgebra associated to the Clifford affinization is given by $F(\mathfrak{g})=\corps{C}[\partial] \mathfrak{g} + \corps{C}K$, with $\partial K=0$ and the $\lambda$-bracket defined by $[a_\lambda b]=<a,b> K$, where $a,b\in\mathfrak{g}$.
\end{prop}
\begin{proof}
The translation in terms of the $\lambda$-bracket is straightforward.
\end{proof}

\chapter{Toward vertex algebras}

We address the last necessary tools needed to define the notion of vertex algebra and reveal the interplay between the algebraic structures studied so far.
Common aspects of Conformal Field Theory\index{Conformal Field Theory} (CFT) will be seen to emerge naturally in the framework of local formal distributions, such as \textit{operator product expansions} and \textit{normal products}.

\section{Formal Cauchy's formulas}
The following will be used to define the notion of \textit{operator product expansion}.
\begin{defin}
Let $a(z)$ be a formal distribution.
\begin{eqnarray*}
a(z)_+ &\doteq& \sum_{n\leq -1} a_{(n)} z^{-1-n} \\
a(z)_- &\doteq& \sum_{n\geq 0} a_{(n)} z^{-1-n} \ .
\end{eqnarray*}
\end{defin}

\begin{prop}
The partition defined by $a(z)_{\pm}$ admits the following properties :
\begin{enumerate}
\renewcommand{\labelenumi}{\emph{(\theenumi)}}
\item Let the map $\pm$ be defined by $\pm\left( a(z) \right)\doteq a(z)_\pm$.  We have $[\pm , \partial_z]=0$, i.e. $\left( \partial_z a(z) \right)_\pm = \partial_z \left( a(z)_\pm \right)$.
\item \textit{Formal Cauchy's formulas}\index{formal Cauchy's formulas} : 
\begin{eqnarray*}
\frac{\partial_w^{n} a(w)_+}{n!} &=& \text{Res}_z a(z)i_{z,w} \frac{1}{(z-w)^{n+1}} \\
\frac{\partial_w^{n} a(w)_-}{n!} &=& - \text{Res}_z a(z)i_{w,z} \frac{1}{(z-w)^{n+1}}
\end{eqnarray*}
\end{enumerate}
\end{prop}
\begin{proof}
As $a(z)=\sum_{n} a_{(n)} z^{-1-n}$,
\begin{eqnarray*}
\partial_z a(z) &=& \sum_{n\in \corps{Z}} a_{(n)} (-n-1) z^{-n-2} \\
&=& \sum_{n\leq -2} a_{(n)} (-n-1) z^{-n-2} + \sum_{n\geq -1} a_{(n)} (-n-1) z^{-n-2} \\
&=& \sum_{n\leq -2} a_{(n)} (-n-1) z^{-n-2} + \sum_{n\geq 0} a_{(n)} (-n-1) z^{-n-2}
\end{eqnarray*}
and let's set
\begin{eqnarray}
\left( \partial_z a(z) \right)_+ &\doteq& \sum_{n\leq -2} a_{(n)} (-n-1) z^{-n-2} \label{partieup} \\
\left( \partial_z a(z) \right)_- &\doteq& \sum_{n\geq 0} a_{(n)} (-n-1) z^{-n-2} \label{partiedown} \ .
\end{eqnarray}
Now,
\begin{eqnarray*}
\partial_z \left( a(z)_+ \right) &=& \sum_{n\leq-1} a_{(n)} (-n-1) z^{-n-2} \,=\, \sum_{n\leq-2} a_{(n)} (-n-1) z^{-n-2} \\
\partial_z \left( a(z)_- \right) &=& \sum_{n\geq 0} a_{(n)} (-n-1) z^{-n-2} \,=\, \sum_{n\geq -1} a_{(n)} (-n-1) z^{-n-2} \ \\
\end{eqnarray*}
which, by (\ref{partieup}) and (\ref{partiedown}), shows the first property. 

For the second one, we compute
\begin{eqnarray*}
\text{Res}_z a(z) i_{z,w}\frac{1}{(z-w)} &=& \text{Res}_z \sum_n a_{(n)} z^{-1-n} z^{-1}\sum_{m\geq 0} \left( \frac{w}{z} \right)^m \\
&=& \text{Res}_z \sum_n a_{(n)} z^{-2-n-m} \sum_{m\geq 0} w^m \\
&=&  \sum_{m\geq 0} a_{(-m-1)} w^m \,=\,  a(w)_+
\end{eqnarray*}
Applying $\partial_w^n$ to both members, we get the first Cauchy's formula.  The proof of the second one is analogous.\\

\end{proof}
\begin{remarque}
In the second property, the analogy with the Cauchy's formula in Complex Analysis is clear.  Recall that for a holomorphic function in a domain $\overline{D}$ ($D$ and its boundary), we have
$$ \frac{\partial_z^n f(a)}{n!} = \frac{1}{2\pi i} \int_{\partial D} \frac{f(z)}{(z-a)^{n+1}}dz \ . $$
\end{remarque}

\section{Operator Product Expansion} \label{dev4687}
Let $\mathfrak{g}$ be a Lie superalgebra.   As we have seen, for a local pair of $\mathfrak{g}$-valued formal distributions,  according to the decomposition theorem and the definition of $j$-products, we can write
\begin{eqnarray}
[a(z),b(w)]=\sum_{j\in\corps{Z}^+} (a(w)_{(j)} b(w)) \frac{\partial_w^j \delta(z,w)}{j!} \ 
\end{eqnarray}
This shows that the $j$-products control the \textit{singularity} (for $z\to w$) of their commutator.  It turns out that the notion of \textit{Operator Product Expansion}\index{operator product expansion} is nothing more than a corollary of this theorem.

The notion of operator product expansion is introduced in this context via a particular formal distribution.
\begin{defin}
Let $a(z)$, $b(z)$ be two $\mathfrak{g}$-valued formal distributions, where $\mathfrak{g}$ is an associative superalgebra.  We define the following $\mathfrak{g}$-valued formal distribution in \emph{two} indeterminates, denoted by $:a(z)b(w):$, by
\begin{eqnarray}
:a(z)b(w):  \doteq a(z)_+ b(w) + p(a,b)b(w)a(z)_- \ . \label{defpop1}
\end{eqnarray}
\end{defin}
This definition leads to the notion of \textit{normal product} between formal distributions.

\begin{remarque}
Using formal Cauchy's formulas, we can write
\begin{eqnarray*}
:a(z)b(w): \, \doteq \text{Res}_x \left( a(x)b(w) i_{x,z} \frac{1}{x-z} - p(a,b)b(w)a(x)i_{z,x}\frac{1}{x-z}\right) \ . \label{defpop2}
\end{eqnarray*}
\end{remarque}

Let $\mathfrak{g}$ be an associative algebra.
The following theorem defines the notion of Operator Product Expansion.  Its proof stresses the crucial role played by the assumption of locality (by pair) of formal distributions.
\begin{theo}
Let $a(z)$, $b(w)$ be a pair of mutually local $\mathfrak{g}$-valued formal distributions.
Then their product reads
\begin{eqnarray*}
a(z)b(w) &=& \sum_{j\in \corps{Z}^+} a(w)_{(j)}b(w) i_{z,w} \frac{1}{(z-w)^{j+1}} \, + \,  :a(z)b(w): \\
p(a,b)b(w)a(z) &=& \sum_{j\in \corps{Z}^+} a(w)_{(j)}b(w) i_{w,z} \frac{1}{(z-w)^{j+1}} \, + \,  :a(z)b(w):
\end{eqnarray*}
Those expressions define the notion of \textit{Operator product expansion}.
\end{theo}
\begin{proof}
Using the decomposition theorem and the expression of $\partial_w^j \delta(z,w)/j!$ in terms of expansions $i_{z,w}$ and $i_{w,z}$, we obtain
\begin{eqnarray}
[a(z),b(w)] &=& \sum_{j\in\corps{Z}^+} a(w)_{(j)}b(w) \frac{\partial_w^j \delta(z,w)}{j!} \\
&=& \sum_{j\in\corps{Z}^+} a(w)_{(j)}b(w) \left( i_{z,w} \frac{1}{(z-w)^{j+1}} - i_{w,z} \frac{1}{(z-w)^{j+1}} \right) \label{epo1}
\end{eqnarray}
with, for $j\in\corps{Z}^+$, $a(w)_{(j)}b(w)=\Res{z} (z-w)^j [a(z),b(w)] \, \in \mathfrak{g}[[w^{\pm 1}]]$.

We then separate negative and positive powers in $z$ in both members of (\ref{epo1}) to obtain
\begin{eqnarray}
\lbrack a_-(z),b(w) \rbrack &=& \sum_{j\in\corps{Z}^+} a(w)_{(j)}b(w) \, i_{z,w} \frac{1}{(z-w)^{j+1}} \label{epo2} \\
\lbrack a_+(z),b(w)\rbrack &=& - \sum_{j\in\corps{Z}^+} a(w)_{(j)}b(w) \, i_{w,z} \frac{1}{(z-w)^{j+1}}  \label{epo3}\ .
\end{eqnarray}
Indeed, recall $i_{z,w}$ and $i_{w,z}$ gives an expansion with positive powers in $w/z$ (the domain being $|z|>|w|$) and $z/w$ (the domain being $|w|>|z|$) respectively.

Now, we also have
\begin{eqnarray*}
\lbrack a_-(z),b(w) \rbrack &=& a(z)_- b(w) - p(a,b) b(w)a(z)_- \\
\lbrack a_+(z),b(w)\rbrack &=&  a(z)_+ b(w) - p(a,b) b(w)a(z)_+ \ , 
\end{eqnarray*}
from which we isolate the terms $a(z)_- b(w)$ and $p(a,b) b(w)a(z)_+$, and insert them into
\begin{eqnarray*}
a(z) b(w)  &=& a(z)_- b(w) + a(z)_+ b(w) \\
p(a,b) b(w) a(z)  &=& p(a,b) b(w) a(z)_- + p(a,b) b(w)a(z)_+ \ .
\end{eqnarray*}
Finally,
\begin{eqnarray*}
a(z) b(w)  &=& \lbrack a_-(z),b(w) \rbrack + a(z)_+ b(w) + p(a,b) b(w)a(z)_- \,=\, \lbrack a_-(z),b(w) \rbrack + :a(z)b(w): \\
p(a,b) b(w) a(z)  &=& - \lbrack a_+(z),b(w) \rbrack + a(z)_+ b(w) + p(a,b) b(w)a(z)_- \,=\, - \lbrack a_+(z),b(w) \rbrack + :a(z)b(w): \ ,
\end{eqnarray*}
in which we insert (\ref{epo2}) and (\ref{epo3}).
\end{proof}

\begin{remarque}[Notation in Physics]
In general, the Operator Product expansion admits a slightly different notation in Physics\index{operator product expansion!notation in Physics}.  The regular part is often omitted, as well as the expansion $i_{z,w}$.  This simply reduces to :
$$ a(z)b(w) \sim \sum_{j\in \corps{Z}^+} \frac{a(w)_{(j)}b(w)}{(z-w)^{j+1}} \ .$$
\end{remarque}

\begin{exe}[Virasoro] Starting from (\ref{viradis}), we previously deduced the non vanishing $j$-products. Hence, we obtain\index{operator product expansion!Virasoro} :
$$ L(z)L(w) \sim \frac{c/2}{(z-w)^4} + \frac{2L(w)}{(z-w)^2} + \frac{\partial_w L(w)}{z-w} \ . $$
\end{exe}
\begin{exe}[Neveu-Schwartz]\index{operator product expansion!Neveu-Schwartz}
The product $L(z)L(w)$ is the same as above.\\
Starting from (\ref{nslambda}) : $[L_\lambda G]= \partial G + \frac{3}{2} G \lambda$, we get : 
$$ L_{(0)}G = \partial G \quad \text{et} \quad L_{(1)}G=\frac{3}{2} G \ , $$ the other $j$-products being zero. Hence
$$ L(z) G(w) \sim \frac{3 G(w)/2}{(z-w)^2} + \frac{\partial_w G(w)}{z-w} \ .$$
Similarly, we start from : $[G_\lambda L]=\frac{1}{2} \partial G + \frac{3}{2} G \lambda$, and we get :
$$ G(z) L(w) \sim \frac{3 G(w)/2}{(z-w)^2} + \frac{\partial_w G(w)/2}{z-w} \ . $$
Finally, since $[G_\lambda G]=L + \frac{c}{6} \lambda^2 = L + \frac{c}{3} \frac{\lambda^2}{2!}$, we get :
$$ G(z)G(w) \sim \frac{c/3}{(z-w)^3} + \frac{L(w)}{z-w} \ .$$
\end{exe}

\section{Normal ordered product}
The operator product expansion is thus composed of two parts.  One part is \textit{singular} in $z=w$ and the other part is \textit{regular} in $z=w$ in a sense to be clarified. The latter is the term $:a(z)b(w):$.  Its limit for $z\to w$ leads to
$$ :a(w)b(w): \, \doteq a(w)_+ b(w) + p(a,b)b(w)a(w)_- \ ,$$
which is not \textit{a priori} a well defined formal distribution, due to the forbidden products in both terms of the RHS.
The apparent problem is solved by introducing the notion of \textbf{field representation}\index{field representation}.
\begin{defin}
Let $(\mathfrak{g},F)$ be a formal distributions Lie superalgebra, where $F=\{a(z)\}$ is a local family of $\mathfrak{g}$-valued formal distributions, whose Fourier coefficients generate $\mathfrak{g}$. Let $(\rho,V)$ be a representation of $\mathfrak{g}$ in a linear space $V$ and $A(z)$ the $End(V)$-valued formal distribution associated to $a(z)$ in the following sense :
$$ a_{(n)} \mapsto A_{(n)} \doteq \rho\left( a_{(n)} \right) \in End(V) \ .$$
The Lie superalgebra $\mathfrak{g}$ is said to be \emph{represented} by means of \emph{fields} $A(z)$ acting on $V$ if the following condition is satisfied :
\begin{eqnarray}
\forall v \in V \,\, \exists N \in \corps{Z} \, : \, A_{(n)} v=0 \,\,\forall n\geq N \ .
\end{eqnarray}
Therefore, $A(z)v\in V[[z]][z^{-1}]$ $\forall v\in V$.
\end{defin} 
Despite the slight abuse of notation, we often write $a(z)$ instead of $A(z)$. To avoid any trouble, we will always assume to work within this framework when dealing with normal products\index{normal product}.

The following proposition shows that the normal product of two fields is again a field.
\begin{prop}\label{ponchamps}
Let $a(z),\, b(z) \in \text{End}V[[z,z^{-1}]]$ be two fields.  Then $:a(z)b(z): \, \in \text{End}V[[z,z^{-1}]]$ is again a field. The space $\text{End}V[[z,z^{-1}]]$ is then given a structure of algebra with respect to the normal product.
\end{prop}
\begin{proof}
We have to check that for $v\in V$, $:a(z)b(z):v=(a(z)_+ b(z) + p(a,b)b(z)a(z)_-)v$ is a Laurent series in $V[z^{-1}][[z]]$.  
But the first term is product between a series expansion with positive powers and a Laurent series ($b(z)$ is a field).  Hence the result is also a Laurent series.
For the second term, since $a(z)$ is a field, $a(z)_- v \in V$ is a finite sum.  Applying $b(z)$ on each term, we get a finite sum of Laurent series, which leads to a Laurent series again.
\end{proof}
We thus defined a notion of \textit{normal ordered product} between fields $a(w)$, $b(w)\in\text{End}V[[z,z^{-1}]]$ as the expression of the formal distribution $:a(z)b(w):$, taken at the limit for $z\to w$, which leads to a well-defined field.

%

The normal product between $a(z)$ and $b(z)$ is explicitly given in the following proposition.
\begin{prop}\label{prodordcomp}
Let $a(z)$ and $b(z)$ be two fields. Their normal ordered product reads explicitly
$$ :a(z)b(z): = \sum_{j\in \corps{Z}} :ab:_{(j)} z^{-j-1} \ ,$$
with
$$ :ab:_{(j)} = \sum_{n\leq -1} a_{(n)} b_{(j-n+1)} + p(a,b) \sum_{n\geq 0} b_{(j-n+1)}a_{(n)} \ .$$
\end{prop}
\begin{proof}
Starting from the definition, we have
\beqns
 :a(z)b(z) : &=& a(z)_+ b(z) + p(a,b) b(z) a(z)_- \\
  &=& \sum_{n=-1}^{-\infty} \sum_{m\in \corps{Z}} a_{(n)} b_{(m)} z^{-n-m-2} + p(a,b) \sum_{m\in \corps{Z}} \sum_{n=0}^\infty b_{(m)} a_{(n)} z^{-n-m-2}
\eeqns
Setting $j=n+m-1$, we obtain
\beqns
 :a(z)b(z) : &=& \sum_{n=-1}^{-\infty} \sum_{m\in \corps{Z}} a_{(n)} b_{(m)} z^{-n-m-2} + p(a,b) \sum_{m\in \corps{Z}} \sum_{n=0}^\infty b_{(m)} a_{(n)} z^{-n-m-2} \\
 &=& \sum_{j\in \corps{Z}} \left( \sum_{n=-1}^{-\infty} a_{(n)} b_{(j-n+1)} + p(a,b) \sum_{n=0}^\infty b_{(j-n+1)} a_{(n)} \right) z^{-j-1} \ ,
\eeqns
hence the proposition.
\end{proof}

The normal product will now be shown to be \textit{quasi-commutative}\index{normal product!quasi-commutativity}.  This property will turn out to enter the definition of a vertex algebra (see last chapter) and is also useful in applications.

For a local pair $(a,b)$ of formal distributions, the normal ordered product is commutative by an integral term.
\begin{prop}
For a local pair $(a,b)$, we have :
\beqn
:ab: - p(a,b) :ba: = \int_{-\partial}^0 [a_\lambda b] d\lambda
\eeqn
\end{prop}
\begin{proof}
Define $a(z,w)$ by 
$$ a(z,w) \doteq a(z) b(w) i_{z,w} \frac{1}{z-w} - p(a,b) b(w) a(z) i_{w,z} \frac{1}{z-w} \ . $$
By locality, there exists $N>>0$ such that $(z-w)^N [a(z),b(w)]=0$.  But
$$ [a(z),b(w)] = (z-w) a(z,w) \ ,$$
hence $a(z,w)$ is local since there exists $N>>0$ for which $(z-w)^{N+1} a(z,w)=0$.

In addition, we clearly have $:a(w)b(w):=\Res{z} a(z,w)$.  On the other hand,
$$ :b(w) a(w): = p(a,b) \Res{z} a(w,z) \ .$$
Indeed, 
\beqns
 p(a,b) \Res{z}a(w,z) &=& p(a,b) \Res{z} \left( a(w)b(z) i_{w,z} \frac{1}{z-w} - p(a,b) b(z) a(w) i_{z,w} \frac{1}{w-z} \right) \\
 &=& \Res{z} \left( -p(a,b) a(w) b(z) i_{w,z} \frac{1}{z-w} + b(z) a(w) i_{z,w} \frac{1}{z-w} \right) \\
 &=& :b(w)a(w):
\eeqns

Now let's compute
\beqns
:ab: - p(a,b) :ba: &=& \Res{z} \left( a(z,w) - a(w,z) \right) \\
&=& \left.\left( \Res{z} e^{\lambda(z-w)} \left( a(z,w) - a(w,z) \right) \right)\right|_{\lambda=0} \\
&=& \left.\left( \Fz a(z,w) - \Fz a(w,z) \right)\right|_{\lambda=0} \\
&=& \left.\left( \Fz a(z,w) - \F{z}{w}{-\lambda-\partial_w} a(z,w)\right)\right|_{\lambda=0} 
\eeqns
since $a(z,w)$ is local.  Therefore,
\beqns
:ab: - p(a,b) :ba: &=& \Res{z} \left( e^{\lambda(z-w)} - e^{(-\lambda - \partial_w)(z-w)}  \right)|_{\lambda=0} a(z,w) \\
&=& \Res{z} \left( 1 - e^{-\partial_w (z-w)} \right) a(z,w) \\
&=& \Res{z} \left( 1 - \sum_{n\geq 0} \frac{(-\partial_w)^n}{n!} (z-w)^n \right) a(z,w) \\
&=& - \Res{z} \sum_{n\geq 1} \frac{(-\partial_w)^n}{n!}(z-w)^n a(z,w) \\
&=& - \Res{z} \sum_{n\geq 1} \frac{(-\partial_w)^n}{n!} (z-w)^{n-1} [a(z),b(w)] \\
&=& - \sum_{n\geq 1} \frac{(-\partial_w)^n}{n!} a_{(n-1)}b \\
&=& - \sum_{n\geq 0} \frac{(-\partial_w)^{n+1}}{(n+1)!} a_{(n)}b \\
&=& \left.\left( \sum_{n\geq 0} \frac{\lambda^{n+1}}{(n+1)!} a_{(n)}b \right)\right|^{\lambda=0}_{\lambda=-\partial}   \\
&=& \int_{-\partial}^0 [a_\lambda b] d\lambda
\eeqns
\end{proof}

\section{Generalized $j$-products}
Recall the operator product expansion
$$ a(z)b(w) = \sum_{j\in \corps{Z}^+} a(w)_{(j)}b(w) i_{z,w} \frac{1}{(z-w)^{j+1}} \, + \,  :a(z)b(w): \ . $$
The singular part appears as a sum weighted by the $j$-products (defined for $j\geq 0$).  We will see that the regular part can be expressed as a certain expansion as well, weighted by \textit{new} $j$-products (for $j<0$).  Those two kinds of products, non related at first sight, will be shown to appear as \textit{special cases} of a \textit{generalized $j$-product}\index{$j$-product!generalized} (for $j\in \corps{Z}$).

The expansion of the regular part to be found is analogous to the famous Taylor expansion. The \textit{formal} Taylor expansion\index{formal Taylor expansion} is based on the following proposition.

\begin{prop}
Let $a(z)=\sum_n a_n z^n$ be any formal distribution. We introduce the following formal distribution in two indeterminates $z$ and $w$ as follows :
\begin{eqnarray*}
i_{z,w} a(z-w) &\doteq& \sum_n a_n i_{z,w}(z-w)^n \ .
\end{eqnarray*}
Then
\begin{eqnarray}
i_{z,w}a(z+w)=\sum_{j=0}^\infty \frac{\partial^j a(z)}{j!} w^j \label{tay} \ .
\end{eqnarray}
Both members are understood to be equal as formal distributions in $z$ and $w$ in the domain $|z|>|w|$.  
\end{prop}
For clarity, we will drop the symbol $i_{z,w}$.
\begin{proof}
Let $a(z)=\sum_n a_n z^n$.  Then, $\partial^j a(z)/j!=\sum_n \binom{n}{j}a_n z^{n-j}$.  But by definition, $$a(z+w)\doteq\sum_n a_n (z+w)^n$$ in the domain $|z|>|w|$.  We then obtain the expected result if
\begin{eqnarray}
(z+w)^n=\sum_{j=0}^\infty \binom{n}{j} z^{n-j}w^j \ . \label{s8z9}
\end{eqnarray}
For $n>0$, the expression is clearly verified, regardless of domain.  For $n<0$, we just get the binomial expansion in the region $|z|>|w|$, by replacing $w$ by $-w$ in definition \ref{exppp}.  
\end{proof}
This leads us to the notion of formal Taylor expansion of a formal distribution $a(z)$ \textit{around} $w$. First replacing $z$ by $w$ and then making the change of variable $z\to z-w$ in (\ref{tay}), the domain to be consider becomes $|z-w|<|w|$ and the equality serves as the needed definition.
\begin{defin}
Let $a(z)$ be any formal distribution.  We define its \textbf{formal Taylor expansion} by
$$ a(z) = \sum_{j=0}^\infty \frac{\partial^j a(w)}{j!} (z-w)^j \ ,$$
valid in $|z-w|<|w|$, with $\partial^j a(w)=\partial_z^j a(z)|_{z=w}$.
\end{defin}
Now applying this to $:a(z)b(w):$, we get
\begin{eqnarray}
:a(z)b(w): = \sum_{j\geq 0} \frac{:(\partial^j a(w))b(w):}{j!} (z-w)^j
\end{eqnarray}
in $|z-w|<|w|$. Setting, for $j\in\corps{Z}^+$, 
\begin{eqnarray}
a(w)_{(-1-j)}b(w) \doteq \frac{:(\partial^j a(w))b(w):}{j!} \ ,\label{jproddown}
\end{eqnarray}
which extends the notion of $j$-\textit{product} to negative values of $j$, and urges us to write
\begin{eqnarray}
a(z)b(w) = \sum_{j\in\corps{Z}} \frac{a(w)_{(j)}b(w)}{(z-w)^{j+1}} \ , \label{tayup}
\end{eqnarray}
which is unfortunatly \textbf{erroneous}. Indeed, the sum on the RHS separates in two parts, both valid, but in a \textit{different} region : $|z|>|w|$ for the singular part and $|z-w|<|w|$ for the regular part.  To solve the issue, it can be shown that the expansions have to be restricted to an (arbitrary) order $N\in\corps{Z}^+$ (see \cite{kac1}).  

Let's focus on the $j$-products.  Those are defined as follows ($j\in \corps{Z}^+$) :
\begin{eqnarray}
a(w)_{(j)}b(w) &=& \Res{z} (z-w)^j [a(z),b(w)] \label{nprod1} \\
a(w)_{(-1-j)}b(w) &=& \frac{:(\partial_w^j a(w))b(w):}{j!} \label{nprod2} \ .
\end{eqnarray}
At first sight, these are not related at all. It turns out that, on the contrary, they appear as special cases of a \textit{generalized $j$-product}\index{$j$-product!generalized}.
\begin{prop}
The $j$-products (\ref{nprod1}) and (\ref{nprod2}) are special cases of the following generalized $j$-product defined by
\begin{eqnarray}
a(w)_{(j)}b(w) = \Res{z} \left( i_{z,w} (z-w)^j a(z)b(w) - i_{w,z} (z-w)^j p(a,b) b(w) a(z) \right) \ . \label{nprodgene}
\end{eqnarray}
\end{prop}
\begin{proof}
Since $j\geq 0$, $i_{z,w} (z-w)^j = i_{w,z} (z-w)^j = (z-w)^j$.  Therefore,
\beqns
a(w)_{(j)}b(w) &=& \Res{z} \left( i_{z,w} (z-w)^j a(z)b(w) - i_{w,z} (z-w)^j p(a,b) b(w) a(z) \right)\\
&=&  \Res{z} (z-w)^j [a(z),b(w)] \ .
\eeqns
We thus recover (\ref{nprod1}).

For the other case,
\begin{eqnarray*}
a(w)_{(-1-j)}b(w) &=& \Res{z} \left( i_{z,w} (z-w)^{-1-j}  a(z)b(w) - i_{w,z} (z-w)^{-1-j} p(a,b) b(w) a(z) \right) \\
&=& \left( \Res{z} a(z) i_{z,w} \frac{1}{(z-w)^{j+1}} \right) b(w) + p(a,b) b(w) \left( -\Res{z} a(z) i_{z,w} \frac{1}{(z-w)^{j+1}} \right) \\
&=&  \left( \frac{\partial_w^j a(w)_+}{j!} \right) b(w) + p(a,b) b(w) \left( \frac{\partial_w^j a(w)_-}{j!} \right) \\
&=&  \frac{: \left( \partial_w^j a(w) \right) b(w) :}{j!} \ , 
\end{eqnarray*}
where we used formal Cauchy's formulas.  We thus recover (\ref{nprod2}).
\end{proof}

The following properties are satisfied by generalized $j$-product. They are proved here in the general case.
\begin{prop}
We have
\begin{enumerate}
\renewcommand{\labelenumi}{\emph{(\theenumi)}}
\item $(\partial a)_{(j)}b  = -j a_{(j-1)} b$
\item $\partial \left( a_{(j)} b \right) = \partial a_{(j)}b + a_{(j)} \partial b$
\end{enumerate}
\end{prop}
\begin{proof}
For (1) :
\beqns
(\partial a)_{(j)}b &=& \Res{z} \left( i_{z,w} (z-w)^j \partial_z a(z) b(w) \right) - \Res{z} \left( i_{w,z} (z-w)^j p(a,b) b(w) \partial_z a(z) \right) \\
&=& -\Res{z} \partial_z \left( i_{z,w} (z-w)^j \right) a(z) b(w) + \Res{z} \partial_z \left( i_{w,z} (z-w)^j  \right) p(a,b) b(w) a(z)
\eeqns
But $[\partial_\cdot,i_{\cdot,\cdot}]=0$, so
\beqns
\partial_z \left( i_{z,w} (z-w)^j  \right) &=& j i_{z,w} (z-w)^{j-1} \\
\partial_z \left( i_{w,z} (z-w)^j  \right) &=& j i_{w,z} (z-w)^{j-1}
\eeqns
Hence,
\beqns
(\partial a)_{(j)}b &=& -j \left( \Res{z} \left( i_{z,w} (z-w)^{j-1} a(z) b(w) - i_{w,z} (z-w)^{j-1} p(a,b) b(w) a(z) \right) \right) \\
&=& -j a_{(j-1)}b
\eeqns
For (2), using (1) :
\beqns
\partial \left( a_{(j)} b \right) &=& \partial_w \left( \Res{z} \left( i_{z,w} (z-w)^{j} a(z) b(w) - i_{w,z} (z-w)^{j} p(a,b) b(w) a(z) \right) \right) \\
&=& \Res{z}\left( -j i_{z,w} (z-w)^{j-1} a(z)b(w) + i_{z,w} (z-w)^j a(z) \partial_w b(w) \right. \\
&& \quad\quad + \left. j i_{w,z} (z-w)^{j-1} p(a,b) b(w) a(z) - i_{w,z} (z-w)^j p(a,b) \partial_w b(w) a(z)  \right) \\
&=& -j \Res{z} \left( i_{z,w} (z-w)^{j-1} a(z) b(w) - i_{w,z} (z-w)^{j-1} p(a,b) b(w) a(z) \right) \\
&& \quad\quad + \Res{z} \left( i_{z,w} (z-w)^{j} a(z) b(w) - i_{w,z} (z-w)^{j} p(a,b) b(w) a(z) \right) \\
&=& -j a_{(j-1)} b + a_{(j)}\partial b \\
&=& \partial a_{(j)} b + a_{(j)}\partial b
\eeqns
\end{proof}

\section{Dong's lemma}\label{dong}
Dong's lemma\index{Dong's lemma} was used in a previous proposition.  We now prove it in the case of generalized $j$-products.
\begin{prop}
If $a(w),b(w)$ and $c(w)$ are mutually local formal distributions, then $(a,b_{(n)}c)$ is a local pair as well for $n\in\corps{Z}$.
\end{prop}
\begin{proof}
Since $a(w),b(w)$ and $c(w)$ are mutually local formal distributions, there exists $r\in \corps{Z}^+$ such that
\beqns
(z-x)^r a(z) b(x) &=& p(a,b) (z-x)^r b(x) a(z) \\
(x-w)^r c(w) b(x) &=& p(c,b) (x-w)^r b(x) c(w) \\
(z-w)^r a(z) c(w) &=& p(a,c) (z-w)^r c(w) a(z)
\eeqns

What we want to show is that there exists $m\in \corps{Z}^+$ such that
\beqn
\underbrace{(z-w)^m a(z) \left( b(w)_{(n)}c(w) \right)}_{(*)} = \underbrace{p(a,[b,c]) (z-w)^m \left( b(w)_{(n)} c(w) \right) a(z)}_{(**)} \ , \label{equaaa}
\eeqn
where $b(w)_{(n)} c(w)=\Res{x}\left( b(x) c(w) i_{x,w} (x-w)^n - p(b,c) c(w) b(x) i_{w,x} (x-w)^n \right)$.

We can assume $r$ to be large enough so that $r+n>0$.  Setting $m=4r$ and using the identity
$$z-w = (z-x) + (x-w)$$ in $(*)$,  we can write
\beqn
(z-w)^{4r} = (z-w)^r \sum_{i=0}^{3r} \binom{3r}{i} (z-x)^i (x-w)^{3r-i} \ . \label{eqnaaabbb}
\eeqn
We then have
\beqns
(*) &=& (z-w)^m a(z)\Res{x}\left( b(x) c(w) i_{x,w} (x-w)^n - p(b,c) c(w) b(x) i_{w,x} (x-w)^n \right) \\
&=& \sum_{i=0}^{3r} \binom{3r}{i} A_i 
\eeqns
where 
$$A_i \doteq \Res{x} \left( (z-w)^r   (z-x)^i (x-w)^{3r-i} a(z) \left( b(x) c(w) i_{x,w} (x-w)^n - p(b,c) c(w) b(x) i_{w,x} (x-w)^n \right) \right) \ . $$
Since $3r-i\geq 0$, $(x-w)^{3r-i} i_{x,w} (x-w)^n = i_{x,w} (x-w)^{3r-i+n}$ and the same for $i_{w,x}$.  

Hence
$$A_i = \Res{x} \left( (z-w)^r   (z-x)^i a(z) \left( b(x) c(w) i_{x,w} (x-w)^{3r-i+n} - p(b,c) c(w) b(x) i_{w,x} (x-w)^{3r-i+n} \right) \right) \ .$$

If $i\leq r$, then $3r-i \geq 2r$ and so $3r-i+n \geq 2r+n > r$. Therefore,
$$A_i = \Res{x} \left( (z-w)^r   (z-x)^i a(z) (x-w)^{3r-i+n} [ b(x), c(w)] \right)$$
and $A_i = 0$ by the locality of the pair $(b,c)$.

If $i>r$, then
\beqns
A_i &=& \Res{x} \left( (z-w)^r   (z-x)^i (x-w)^{3r-i} a(z) \left( b(x) c(w) i_{x,w} (x-w)^n - p(b,c) c(w) b(x) i_{w,x} (x-w)^n \right) \right) \\
&=& \Res{x} \left( (z-w)^r (x-w)^{3r-i}  \underbrace{(z-x)^i a(z)b(x)}_{=p(a,b)(z-x)^i b(x) a(z)} c(w) i_{x,w} (x-w)^n \right. \\
&&  \quad \quad \left. - p(b,c) (z-x)^i \underbrace{(z-w)^r a(z) c(w)}_{=p(a,c) (z-w)^r c(w) a(z)} b(x) i_{w,x} (x-w)^n \right) \\
&=& \Res{x} \left(  (x-w)^{3r-i} p(a,b)(z-x)^i b(x) \underbrace{(z-w)^r a(z)c(w)}_{=p(a,c)(z-w)^r c(w) a(z)} i_{x,w} (x-w)^n \right. \\
&&  \quad \quad \left. - p(b,c)  p(a,c) (z-w)^r c(w) \underbrace{(z-x)^i a(z) b(x)}_{=p(a,b)(z-x)^i b(x) a(z)} i_{w,x} (x-w)^n \right) \\
&=& p(a,[b,c]) \Res{x} \left( (x-w)^{3r-i} (z-x)^i b(x)(z-w)^r c(w) a(z) i_{x,w} (x-w)^n \right. \\
&&  \quad \quad \left. - p(b,c) (z-w)^r c(w)(z-x)^i b(x) a(z) i_{w,x} (x-w)^n \right) \\
&=& p(a,[b,c]) \Res{x} \left( (z-w)^r (z-x)^i (x-w)^{3r-i} \left( b(x) c(w) i_{x,w} (x-w)^n - p(b,c) c(w) b(x) i_{w,x} (x-w)^n \right) a(z)  \right)
\eeqns

and we can write 
$$ (*) = \sum_{i=r+1}^{3r} \binom{3r}{i} A_i \ .$$

Let's check that
$$ \sum_{i=r+1}^{3r}\binom{3r}{i} A_i = (**) \ .$$  
But, starting from $(**)$ and using (\ref{eqnaaabbb}), we can write 
$$(**)=\sum_{i=0}^{3r} \binom{3r}{i} B_i \ , $$ 
where $B_i = p(a,[b,c]) \Res{x} \left( (z-w)^r (z-x)^i (x-w)^{3r-i} \left( b(x) c(w) i_{x,w}(x-w)^n - p(b,c) c(w) b(x) i_{w,x} (x-w)^n \right) a(z) \right)$.

For $i\leq r$, we get $B_i = 0$ as previously.  
For $i>r$, we have
$B_i = A_i$.
\end{proof}

\section{Wick's formulas}
As the quasi-commutativity of the normal product, Wick's formulas\index{Wick's formulas} will turn out to enter the definition of a vertex algebra in the last chapter.  They are also useful in applications (in computing $\lambda$-brackets), as shown at the end of this section.
    \subsection{General formula}
    \begin{prop} 
    The general Wick's formula\index{Wick's formulas!general} is the following equality, valid for $\mathfrak{g}$-valued formal distributions, where $\mathfrak{g}$ is an \textbf{associative Lie superalgebra}.
    \beqn
    [a_{\lambda} \left( b_{(n)}c \right)] = \sum_{k\in \corps{Z}^+} \frac{\lambda^k}{k!} [a_\lambda b]_{(n+k)} c + p(a,b) b_{(n)} [a_\lambda c] \label{wick1}
    \eeqn
    \end{prop}
    \begin{proof}
    \beqns
    [a(w)_\lambda \left( b(w)_{(n)}c(w) \right)] &=& \Res{z} e^{\lambda(z-w)} [a(z),b(w)_{(n)}c(w)] \\
    &=& \underbrace{\Res{z} e^{\lambda(z-w)} \Res{x}[a(z),b(x)c(w)] i_{x,w} (x-w)^n}_{(*)} \\
    && \quad \quad \underbrace{- p(b,c) \Res{z} e^{\lambda(z-w)} \Res{x} [a(z),c(w)b(x)]i_{w,x} (x-w)^n}_{(**)}
    \eeqns
    
    Since $\mathfrak{g}$ is associative, we have 
    $$[a,bc]=[a,b]c + p(a,b) b[a,c]$$
    which is easily proved.
    
    Using it in $(*)$ and $(**)$, we can write
    \beqns
    (*) &=& \Res{z} \Res{x} e^{\lambda(z-w)} [a(z),b(x)] c(w) i_{x,w} (x-w)^n \\
        && \quad \quad + p(a,b) \Res{z} \Res{x} e^{\lambda(z-w)} b(x)[a(z),c(w)] i_{x,w} (x-w)^n
    \eeqns
    Using the identity $e^{\lambda(z-w)}=e^{\lambda(z-x)}e^{\lambda(x-w)}$, we obtain
    \beqns
    (*) &=& \Res{x} e^{\lambda(x-w)} i_{x,w} (x-w)^n \underbrace{\Res{z} e^{\lambda(z-x)} [a(z),b(x)]}_{=[a(x)_\lambda b(x)]} c(w) \\
        &&  + p(a,b) \Res{x} b(x) i_{x,w} (x-w)^n \underbrace{\Res{z} e^{\lambda(z-w)} [a(z),c(w)]}_{=[a(w)_\lambda c(w)]} \\
        &=& \Res{x} \sum_{j\geq 0} \frac{\lambda^j}{j!} i_{x,w} (x-w)^{n+j} [a(x)_\lambda b(x)] c + p(a,b) \Res{x} b(x) i_{x,w} (x-w)^n [a(w)_\lambda c(w)] \ ,
    \eeqns
    where we used at the last line that $[(z-w)^j, i_{x,w}]=0$ for $j\geq 0$.
    \beqns
    (**) &=& -p(b,c) \Res{z} e^{\lambda(z-w)} \Res{x} [a(z), c(w)b(x)] i_{w,x} (x-w)^n \\
         &=& -p(b,c) \Res{x} \underbrace{\Res{z} e^{\lambda(z-w)} [a(z),c(w)]}_{=[a(w)_\lambda b(w)]} b(x) i_{w,x} (x-w)^n \\
          &&  -p(b,c) p(a,c) \Res{x} e^{\lambda(x-w)} c(w) \underbrace{\Res{z} e^{\lambda(z-x)}[a(z),b(x)]}_{=[a(x)_\lambda b(x)]} i_{w,x}(x-w)^n \\
         &=& -p(b,c) \Res{x} [a(w)_\lambda c(w)] b(x) i_{w,x}(x-w)^n - p(ab,c) \Res{z} \sum_{j\geq 0} \frac{\lambda^j}{j!} c(w) [a(x)_\lambda b(x)] i_{w,x} (x-w)^{n+j}
    \eeqns
    Hence,
    \beqns
    (*) + (**) &=& \Res{x} \sum_{j\geq 0} \frac{\lambda^j}{j!} i_{x,w} (x-w)^{n+j} [a(x)_\lambda b(x)] c + p(a,b) \Res{x} b(x) i_{x,w} (x-w)^n [a(w)_\lambda c(w)] \\
    &&   -p(b,c) \Res{x} [a(w)_\lambda c(w)] b(x) i_{w,x}(x-w)^n - p(ab,c) \Res{z} \sum_{j\geq 0} \frac{\lambda^j}{j!} c(w) [a(x)_\lambda b(x)] i_{w,x} (x-w)^{n+j} \\
    &=& \sum_{j\geq 0} \frac{\lambda^j}{j!} \Res{x} \left( i_{x,w} (x-w)^{n+j} [a(x)_\lambda b(x)] c(w) - p(ab,c) i_{w,x} (x-w)^{n+j} c(w) [a(x)_\lambda b(x)] \right) \\
    &&  + p(a,b) \Res{x} \left( i_{x,w} (x-w)^n b(x) [a(w)_\lambda c(w)] - \underbrace{p(b,c)p(a,b)}_{=p(b,ac)} i_{w,x}(x-w)^n [a(w)_\lambda c(w)] b(x) \right) \\
    &=& \sum_{j\geq 0} \frac{\lambda^j}{j!} [a_\lambda b]_{(n+j)}c + p(a,b) b_{(n)}[a_\lambda c]
    \eeqns
    \end{proof}
    
    \begin{remarque}
    Notice that Jacobi identity for the $\lambda$-bracket is hidden in this result.  We just have to apply $\sum_{n\in\corps{Z}^+} \frac{\mu^n}{n!}$ to both members of (\ref{wick1}), as computed below.
    \begin{proof}
    The LHS becomes : 
    \beqns
    \sum_{n\geq 0} \frac{\mu^n}{n!}  \lbrack a_\lambda \left( b_{(n)}c \right) \rbrack  &=&  \lbrack a_\lambda \left( \sum_{n\geq0}\frac{\mu^n}{n!} b_{(n)}c \right) \rbrack  \\
    &=&  \lbrack a_\lambda  \lbrack b_\mu c \rbrack  \rbrack 
    \eeqns
    
    The RHS :
    \beqns
    \underbrace{\sum_{n\geq 0} \frac{\mu^n}{n!}  \sum_{k\geq 0} \frac{\lambda^k}{k!} [a_\lambda b]_{(n+k)}c}_{(*)} + p(a,b) \underbrace{\sum_{n\geq 0} \frac{\mu^n}{n!}  b_{(n)}[a_\lambda c]}_{(**)} \ ,
    \eeqns
    with, setting $k=m-n$ (hence $m\geq n$),
    \beqns
      (*) &=& \sum_{m\geq n} \sum_{n\geq 0} \frac{\mu^n \lambda^{m-n}}{n! (m-n)!} [a_\lambda b]_{(m)}c  \\
      &=&  \sum_{m\geq 0} \sum_{n\geq 0} \binom{m}{n} \frac{\mu^n \lambda^{m-n}}{m!} [a_\lambda b]_{(m)}c \\
      &=& \sum_{m\geq 0} \frac{(\mu + \lambda)^m}{m!} [a_\lambda b]_{(m)}c \\
      &=& [[a_\lambda b]_{\mu+\lambda}c]
    \eeqns
    \beqns
      (**) = p(a,b) [b_\mu[a_\lambda c]]
    \eeqns
    Therefore,
    $$\lbrack a_\lambda  \lbrack b_\mu c \rbrack  \rbrack = [[a_\lambda b]_{\mu+\lambda}c] + p(a,b) [b_\mu[a_\lambda c]] $$
    \end{proof}
    \end{remarque}
    
    \subsection{Non abelian Wick's formula}
    We get the so-called \textit{non abelian Wick's formula}\index{Wick's formulas!non abelian} by setting $n=-1$ in (\ref{wick1}) and expressing the sum by an integral.
    \begin{prop}
    The non abelian Wick's formula is the following equality :
    \beqn
    [a_\lambda : bc :] = : [a_\lambda b] c : + p(a,b) : b[a_\lambda c] : + \int_0^\lambda [[a_\lambda b]_\mu c] d\mu \label{wick2}
    \eeqn
    \end{prop}
    \begin{proof}
    Setting $n=-1$ in the general Wick's formula,
    \beqns
    [a_\lambda :bc:] &=& \sum_{k=0}^\infty \frac{\lambda^k}{k!} [a_\lambda b]_{(-1+k)} c + p(a,b) :b[a_\lambda c]: \\
    &=& :[a_\lambda b] c : + p(a,b) : b[a_\lambda c] : + \sum_{k\geq 1} \frac{\lambda^k}{k!} [a_\lambda b]_{(k-1)} c
    \eeqns
    In this expression, the sum can be written as a formal integral :
    \beqns
    \sum_{k\geq 1} \frac{\lambda^k}{k!} [a_\lambda b]_{(k-1)} c &=& \sum_{k\geq 0} \frac{\lambda^{k+1}}{(k+1)!} [a_\lambda b]_{(k)} c \\
    &=& \int_0^\lambda [[a_\lambda b]_\mu c] d\mu \ ,
    \eeqns
    hence the result.
    \end{proof}

\subsection{Application : uncharged free superfermions}    
    We show now how Wick's formulas can be used in computing $\lambda$-brackets appearing here in the case of ``uncharged free superfermions''\index{uncharged free superfermions}.
    The setting is as follows.  Let $A$ be a superspace, \textit{i.e.} $A=A_{\bar{0}} \oplus A_{\bar{1}}$, whose \textit{superdimension} is $sdim(A)=dim(A_{\bar{0}}) - dim(A_{\bar{1}})$. We endow $A$ with a non degenerate and antisupersymmetric bilinear form $\langle \cdot,\cdot \rangle$.   We start from the Lie conformal superalgebra 
$$F(A)=\corps{C}[\partial]A + \corps{C}K$$ 
defined by 
\begin{eqnarray*}
\lbrack a_\lambda b\rbrack &=& \langle a,b \rangle K \ , \quad  a,b\in A \\
\lbrack K_\lambda a\rbrack &=& 0 \quad \forall a\in A \ ,
\end{eqnarray*}
The associated formal distributions Lie superalgebra $(\hat{A},\bar{A})$ can be constructed as $(\text{Lie}(F(A)),\overline{F(A)})$, or we can proceed as follows. We define $\hat{A}$ as the space of Laurent polynomials on $A$, with an extra central element $K$
$$\hat{A}=A[t,t^{-1}] + \corps{C}K$$
endowed with the Lie superbracket
$$ [\varphi t^m, \psi t^n ] = K \delta_{m,-1-n} \langle \varphi,\psi \rangle \quad \text{where $\varphi,\psi \in A$.} \ $$
Then, $\bar{A}$ is defined as the $\partial_z$-invariant space of $\hat{A}$-valued formal distributions. Moreover, posing $\partial = -\frac{d}{dt}$, so that $\partial_z a(z)=(\partial a)(z)$, $(\hat{A},\bar{A})$ is said \textit{regular}\index{regular} [for having a derivation $\partial$ on $\hat{A}$ satisfying the previous relation] and we can give $\bar{A}$ the structure of $\corps{C}[\partial]$-module. This allows us to define $\lambda$-brackets on $\bar{A}$ as above (on $F(A)$).  We will deal with normal ordered products, so we work in the framework of representations by fields, acting on a linear space $V$.  The central element $K$ acts as $kI_V$, where $k\in\corps{C}$ and $I_V$ is the identity endomorphism on $V$.  To simplify, we set $k=1$, so that $[a_\lambda b]=\langle a,b \rangle $ where $a,b\in A$.

Now, let $\{ \phi_i \}$ be a basis of $A$ and $\{ \phi^i \}$ the dual basis : $\langle \phi_i , \phi^j \rangle = \delta_{ij}$.  Then $\{\phi, \phi_i \}$ is a local family of formal distributions, since 
$$ [\phi(z),\phi_i (w)]=\langle \phi, \phi_i \rangle \delta(z,w) \ . $$ 
We define
$$ L(z)= \frac{1}{2} \sum_i : \partial_z \phi^i(z) \phi_i (z) : \ . $$
By Dong's lemma, $\phi(z)$ are $L(z)$ are mutually local.  Therefore, their $\lambda$-bracket has to be polynomials in $\lambda$. What we want to do is computing them explicitly.  We state and prove the result in the following proposition.

	\begin{prop} The $\lambda$-brackets are given by
	\begin{enumerate}
	\renewcommand{\labelenumi}{\emph{(\theenumi)}}
	\item $[\phi_\lambda L] = \frac{1}{2}\left(\lambda - \partial  \right) \phi$
	\item $[L_\lambda \phi] = \left( \partial + \frac{\lambda}{2} \right) \phi$
	\item $[L_\lambda L] = - \frac{sdim(A)}{24} \lambda^3 + 2\lambda L + \partial L$, where $sdim(A)=dim(A_{\bar{0}})-dim(A_{\bar{1}})$
	\end{enumerate}
	\end{prop}
	\begin{proof}
	Those three equalities are based on the non abelian Wick's formula.  Let's stress the fact that we assume working within the framework of a representation of the Lie superalgebra $\hat{A}$ by means of $End(V)$-valued fields, so that their Fourier coefficients belongs to an associative Lie superalgebra, as required.
	
	Let's prove (1).  By the non abelian Wick's formula, we can write 
	\beqns
	  [\phi_\lambda L] = \frac{1}{2} \sum_i : [\phi_\lambda \partial \phi^i] \phi_i : + \frac{1}{2} \sum_i p(\phi,\phi_i) : \partial \phi^i [\phi_\lambda \phi_i] : + \frac{1}{2} \sum_i \int_0^\lambda [[\phi_\lambda \partial \phi^i]_\mu \phi_i] d\mu
	\eeqns
	The integral term vanishes, since $[\phi_\lambda \partial \phi^i]=(\partial + \lambda) [\phi_\lambda \phi^i]=\lambda \braket{\phi}{\phi^i}$.
	
	The first term reads
	\beqns
	  \frac{1}{2} \sum_i : [\phi_\lambda \partial \phi^i] \phi_i : &=& \frac{1}{2} \sum_i : \left( (\partial + \lambda) [\phi_\lambda \phi^i] \right) \phi_i : \\
	  &=& \frac{1}{2} \sum_i \lambda  \langle \phi , \phi^i \rangle  : \phi_i : \\
	  &=& \frac{\lambda}{2} \sum_i  \langle \phi , \phi^i \rangle  \phi_i\\
	  &=& \frac{\lambda}{2} \phi
	\eeqns
	For the second, we notice that
	\beqns
	 p(\phi,\phi_i)[\phi_\lambda \phi_i] &=& -[{\phi_i}_{-\lambda-\partial} \phi] \\
	 &=& - [{\phi_i}_{\mu} \phi]|_{\mu=-\lambda-\partial} \\
	 &=& -  \langle \phi_i , \phi \rangle 
	\eeqns
	Hence
	\beqns
	\frac{1}{2} \sum_i p(\phi,\phi_i) : \partial \phi^i [\phi_\lambda \phi_i] : &=& \frac{1}{2} \sum_i (-1)  \langle \phi_i , \phi \rangle \partial \phi^i \\
	&=& -\frac{\partial}{2} \phi
	\eeqns
	Summing both terms, we get the result.
	
	Let's prove (2).  We use anticommutativity.  $L$ being even,
	\beqns
	 [L_\lambda \phi] &=& -[\phi_{-\lambda-\partial}L] \\
	 &=& - \left( \frac{-\lambda - \partial}{2} \phi - \frac{\partial}{2} \phi \right) \\
	 &=& \left( \frac{\lambda}{2} + \partial \right) \phi
	\eeqns
	
	Let's prove (3). $L$ being still even,
	\beqns
	\lbrack L_\lambda L\rbrack  &=& \frac{1}{2} \sum_i \lbrack L_\lambda :\partial \phi^i \phi_i:\rbrack  \\
	&=& \frac{1}{2} \sum_i :\lbrack L_\lambda \partial \phi^i\rbrack  \phi_i : + \frac{1}{2} \sum_i :\partial \phi^i \lbrack L_\lambda \phi_i\rbrack : + \frac{1}{2} \sum_i \int_0^\lambda \lbrack \lbrack L_\lambda \partial \phi^i\rbrack _\mu \phi_i\rbrack  d\mu \\
	&=& \frac{1}{2} \sum_i :(\partial + \lambda) \lbrack L_\lambda \phi^i\rbrack  \phi_i : + \frac{1}{2} \sum_i : \partial \phi^i \lbrack L_\lambda \phi_i\rbrack  : + \frac{1}{2} \sum_i \int_{0}^\lambda \lbrack (\partial + \lambda)\lbrack L_\lambda \phi^i\rbrack _\mu \phi_i\rbrack  d\mu
	\eeqns
	
	Using (1) and (2),
	\beqns
	\lbrack L_\lambda L\rbrack  &=& A + B + C
	\eeqns
  with 
  \beqns
    A &=& \frac{1}{2} \sum_i :(\partial + \lambda)(\partial + \frac{\lambda}{2}) \phi^i \phi_i: \\
    B &=& \frac{1}{2} \sum_i : \partial \phi^i (\partial + \frac{\lambda}{2}) \phi_i : \\
    C &=& \frac{1}{2} \sum_i \int_0^\lambda \partial [L_\lambda \phi^i]_\mu \phi_i d\mu + \frac{1}{2} \sum_i \int_0^\lambda \lambda [[L_\lambda \phi^i]_\mu \phi_i ] d\mu
  \eeqns
	
	Developing each term, we obtain
	\beqns
	  A &=& \frac{1}{2} \sum_i : \partial^2 \phi^i \phi_i : + \frac{1}{2} \sum_i \frac{3}{2} \lambda :\partial \phi^i \phi_i : + \frac{1}{2} \sum_i \frac{\lambda^2}{2} : \phi^i \phi_i : \\
	  B &=& \frac{1}{2} \sum_i : \partial \phi^i \partial \phi_i : + \frac{\lambda}{4} \sum_i : \partial \phi^i \phi_i : \\
	  C &=& \frac{1}{2} \sum_i \int_{0}^\lambda (\lambda - \mu) [(\partial + \frac{\lambda}{2})\phi^i_\mu \phi_i] d\mu
	\eeqns
	
	Moreover,
	\beqns
	  : \partial^2 \phi^i \phi_i : &=& \partial^2 \phi^i_{(-1)} \phi \\
	  &=& \partial \left( \partial \phi^i_{(-1)}\phi_i \right) - \partial \phi^i_{(-1)} \partial \phi_i
	\eeqns
	
	Hence, gathering everything,
	\beqns
	\lbrack L_\lambda L\rbrack  &=&  \partial \left( \underbrace{\frac{1}{2}\sum_i :\partial \phi^i \phi_i:}_{= L}\right) - \frac{1}{2} \sum_i :\partial \phi^i \partial \phi_i:
	+ 2\lambda \left( \underbrace{\frac{1}{2} \sum_i :\partial \phi^i \phi_i :}_{= L} \right) \\
	&& + \quad \frac{\lambda^2}{4} \underbrace{\sum_i :\phi^i \phi_i:}_{(*)} \quad + \quad \frac{1}{2} \sum_i :\partial \phi^i \partial \phi_i: \quad + \quad C
	\eeqns
	
	By symmetry, $(*)=0$.  Indeed, since $(*)$ is a \textit{Casimir} type expression, it does not depend on the basis chosen, so that we can write
	\beqns
	  \sum_i :\phi^i \phi_i : &=& - \sum_i  (-1)^{p(\phi_i)} :\phi_i \phi^i :
	\eeqns
	since the basis $\{\phi^i\}$ and $\{ (-1)^{1+p(\phi_i)} \phi_i \}$ are dual.  [This comes from the fact that the bilinear form is antisupersymmetric and $\langle \phi_i , \phi^j \rangle = \delta_{ij}$, from which we also deduce that $p(\phi_i)=p(\phi^i)$.]

But the normal product being quasi-commutative, i.e. $:\phi_i \phi^i : = (-1)^{p(\phi_i)p(\phi^i)} :\phi^i \phi_i :$ and since $p(\phi_i)=p(\phi^i)$, we have
  \beqns
    :\phi_i \phi^i : = (-1)^{p(\phi_i)} :\phi^i \phi_i :
  \eeqns
	hence
	\beqns
	  \sum_i :\phi^i \phi_i : &=& - \sum_i  (-1)^{p(\phi_i)} :\phi_i \phi^i : \\
	  &=& - \sum_i :\phi^i \phi_i :
	\eeqns
	i.e.
	\beqns
	 \sum_i :\phi^i \phi_i : &=& 0
	\eeqns

Let's compute the contribution of the term $C$ :
  \beqns
    C &=& \frac{1}{2} \sum_i \int_{0}^\lambda (\lambda - \mu) [(\partial + \frac{\lambda}{2})\phi^i_\mu \phi_i] d\mu \\
      &=&  \frac{1}{2} \sum_i \int_{0}^\lambda (\lambda - \mu) \left( [\partial \phi^i_\mu \phi_i] + \frac{\lambda}{2} [\phi^i_\mu \phi_i] \right) d\mu \\
      &=& \frac{1}{2} \sum_i \int_{0}^\lambda (\lambda - \mu)(\frac{\lambda}{2}- \mu) [\phi^i_\mu \phi_i] d\mu \\
      &=& \frac{1}{2} \int_{0}^\lambda (\lambda - \mu)(\frac{\lambda}{2}- \mu) d\mu \sum_i \langle \phi^i , \phi_i \rangle 
  \eeqns
  By antisupersymmetry,
  \beqns
    \sum_i \langle \phi^i , \phi_i \rangle &=& \sum_i (-1)(-1)^{p(\phi^i)} \langle \phi_i , \phi^i \rangle \\
                                           &=& \sum_i (-1)(-1)^{p(\phi^i)} \delta_{ii} \\
                                           &=& -\sum_i (-1)^{p(\phi^i)} \\
                                           &=& - \left( dim A_{\bar{0}} - dim A_{\bar{1}} \right) \\
                                           &=& - sdim A
  \eeqns
	
	The integral gives
	\beqns
	 \int_{0}^\lambda (\lambda - \mu)(\frac{\lambda}{2}- \mu)d\mu &=& \frac{\lambda^3}{12}
	\eeqns
	
	Finally we get the expected result :
	\beqns
	\lbrack L_\lambda L\rbrack  &=& \partial L + 2\lambda L - \frac{sdim A}{24} \lambda^3
	\eeqns
	
	\end{proof}

\section{Eigendistributions and the notion of weight}
As OPE's and NOP's, those notions appear in Conformal Field Theory.  We see now how they translate in this formalism.

\begin{defin}
In the language of $\lambda$-bracket, an eigendistribution\index{eigendistribution} of weight\index{eigendistribution!weight} $\Delta_a$ with respect to $L$ is a formal distribution $a(z)$ satisfying
\begin{eqnarray} 
[L_\lambda a ] = \left( \partial + \Delta_a \lambda \right) a + \cal{O}(\lambda) \ . \label{distrprop}
\end{eqnarray} 
If $L$ is a Virasoro formal distribution, $\Delta_a$ is called the \textit{conformal weight}\index{eigendistribution!conformal weight}.
\end{defin}

\begin{remarque}
Clearly, $a(z)$ is an eigendistribution in the sense that it is an eigenvector of the endomorphism $L_{(1)}$, whose action is defined by $L_{(1)}(b)\doteq L_{(1)}b$.  Indeed, (\ref{distrprop}) implies $L_{(1)}a = \Delta_a a$.
\end{remarque}

\begin{exe}
The distribution $\phi(z)$, introduced previously, is an eigendistribution of weight $\frac{1}{2}$ with respect to $L(z)$.  In addition, $L(z)$ is a Virasoro eigendistribution\index{eigendistribution!Virasoro} of conformal weight $2$ with respect to itself.
\end{exe}

When the respective weights of two formal eigendistributions are known, the we can determine the weight of their $n$-product.
\begin{prop}\label{poidsnproduits}
If $a$ and $b$ have weights $\Delta_a$ and $\Delta_b$ with an even formal distribution $L$, then $a_{(n)}b$ has weight $\Delta_a + \Delta_b - n - 1$ with respect to $L$.  In particular, $\Delta_{:ab:} = \Delta_a + \Delta_b$ and $\Delta_{\partial a} = \Delta_a + 1$.
\end{prop}
\begin{proof}
The proof is, once again, based on the non abelian Wick's formula.  We will then work in the framework of distributions valued in an associative Lie superalgebra.

As $L$ is even, 
\beqns
  [L_\lambda \left(a_{(n)}b \right)] &=& \sum_{k\in \corps{Z}^+} \frac{\lambda^k}{k!} [L_\lambda a]_{(n+k)} b + a_{(n)}[L_\lambda b] \\
  &=& \sum_{k\in \corps{Z}^+} \frac{\lambda^k}{k!} \left( \partial a + \Delta_a \lambda a + \mathcal{O}(\lambda^2) \right)_{(n+k)} b + a_{(n)}\left( \partial b + \Delta_b \lambda b + \mathcal{O}(\lambda^2)\right) \\
  &=& \partial a_{(n)} b + \lambda \partial a_{(n+1)} b + \Delta_a \lambda a_{(n)}b + a_{(n)} \partial b + \Delta_b \lambda a_{(n)}b + \mathcal{O}(\lambda^2) \\
  &=& \left( \partial + (\Delta_a + \Delta_b - n - 1)\lambda \right) \left( a_{(n)}b \right) + \mathcal{O}(\lambda^2)
\eeqns 
This shows $a_{(n)}b$ has weight $\Delta_a + \Delta_b - n - 1$ with respect to $L$.

Setting $n=-1$, we get $\Delta_{:ab:} = \Delta_a + \Delta_b$.

The following computation
\beqns
 [L_\lambda \partial a] &=& (\lambda + \partial) [L_\lambda a] \\
 &=& (\lambda + \partial) (\partial a + \Delta_a \lambda a + \mathcal{O}(\lambda^2)) \\
 &=& \partial^2 a + \Delta_a \lambda \partial a + \lambda \partial a + \mathcal{O}(\lambda^2) \\
 &=& \left( \partial + (\Delta_a + 1)\lambda \right) \partial a + \mathcal{O}(\lambda^2)
\eeqns
shows that $\Delta_{\partial a} = \Delta_a + 1$.
\end{proof}

The expansion of an eigendistribution $a(z)$ of weight $\Delta_a$ is often adapted as follows :
$$ a(z) = \sum_{n\in -\Delta_a + \corps{Z}} a_n z^{-n-\Delta_a} \ .$$
By comparison with the Fourier modes of $a(z)=\sum_{m\in \corps{Z}} a_{(m)} z^{-m-1}$, we must have $m=n+\Delta_a - 1$.  By the way, this confirms that $m\in \corps{Z}$ when $n\in  -\Delta_a + \corps{Z}$.  In terms of the Fourier modes, we get
\beqns
 a(z)=\sum_m a_{m-\Delta_a + 1} z^{-m-1}
\eeqns
hence
$$ a_{(m)} = a_{m-\Delta_a + 1} \ . $$
In the other way,
$$ a_{n} = a_{(n+\Delta_a - 1)} \ .$$
\begin{remarque}
Notice that since $\Delta_{\partial a} = \Delta_a + 1$, we can write :
\beqns
 \partial_z a(z) = \sum_n -(n+\Delta_a) a_n z^{-n-\Delta_a -1} = \sum_n -(n+\Delta_a) a_n z^{-n-\Delta_{\partial a}}
\eeqns
As $(\partial a)(z) = \partial_z a(z)$, we get 
$(\partial a)_n = -(n+\Delta_a) a_n$.
\end{remarque}

One of the interesting features of this change of notation is that it reveals the gradation of the superbracket.
\begin{prop}
In the new notation, we can write
\beqns
 [a_m, b_n] = \sum_{j\in \corps{Z}^+} \binom{m+\Delta_a -1}{j} \left( a_{(j)}b \right)_{m+n} \label{nouvellerelcomm}
\eeqns
\end{prop}
\begin{proof}
This results from the following computation :
\beqns
 [a_m, b_n] &=& [a_{(m+\Delta_a - 1)},b_{(n+\Delta_b - 1)}] \\
            &=& \sum_{j\geq 0} \binom{m+\Delta_a - 1}{j} \left( a_{(j)}b \right)_{(m+\Delta_a-1 + n+\Delta_b - 1 - j)} \\
            &=& \sum_{j\geq 0} \binom{m+\Delta_a - 1}{j} \left( a_{(j)}b \right)_{m+\Delta_a-2 + n+\Delta_b - j - \Delta_{a_{(j)}b} + 1} \\
            &=& \sum_{j\geq 0} \binom{m+\Delta_a - 1}{j} \left( a_{(j)}b \right)_{m+n}
\eeqns
where we used $\Delta_{a_{(j)}b}=\Delta_a + \Delta_b -j -1$ (see proposition \ref{poidsnproduits}).
\end{proof}

\begin{defin}
A formal distribution $a(z)$ is called a primary eigendistribution\index{eigendistribution!primary} of weight $\Delta_a$ if
$$ [L_\lambda a] = (\partial + \Delta_a \lambda) a \ ,$$
where $L$ is Virasoro eigendistribution (hence of weight $2$ with respect to itself).
\end{defin}

\begin{prop}\label{Lacomm}
Let $a$ be a primary eigendistribution of weight $\Delta_a$ with respect to $L$.  Then we have, $\forall n,m \in \corps{Z}$ :
$$ [L_m, a_n] = (m(\Delta_a - 1)-n) a_{m+n} \ .$$
\end{prop}
\begin{proof}
By assumption, 
$$[L_\lambda a] = (\partial + \Delta_a \lambda)a \ ,$$ 
hence
\beqns
  L_{(0)} a &=& \partial a \\
  L_{(1)} a &=& \Delta_a a \\
  L_{(j)} a &=& 0 \quad \forall j\neq \{0,1\}
\eeqns
Knowing that $\Delta_L = 2$, 
\beqns
 [L_m , a_n ] &=& \sum_{j\geq 0} \binom{m+\Delta_L - 1}{j} \left( L_{(j)}a \right)_{m+n} \\
              &=& \left( \partial a \right)_{m+n} + (m+1) \left( \Delta_a a \right)_{m+n} \\
              &=& -(m+n+\Delta_a) a_{m+n} + (m+1) \Delta_a a_{m+n} \\
              &=& (m(\Delta_a - 1) - n) a_{m+n} \ ,
\eeqns
hence the statement.
\end{proof}

\begin{exe}
For free superfermions, we have found $L$ to be a Virasoro eigendistribution and $\phi$ satisfied
$$ [L_\lambda \phi] = (\partial + \frac{\lambda}{2})\phi \ .$$
$\phi$ is thus a primary eigendistribution of weight $\frac{1}{2}$ and its expansion will be written
$$ \phi(z) = \sum_{n\in -\frac{1}{2}+\corps{Z}} \phi_n z^{-n-\frac{1}{2}} \ . $$

\begin{prop}
In the framework of free superfermions and assuming that the weight of $1$ is $0$, then we have :
\begin{enumerate}
\renewcommand{\labelenumi}{\emph{(\theenumi)}}
 \item $[\phi_m, \psi_n] = \langle \phi, \psi \rangle \delta_{m,-n}$
 \item $[L_m, \phi_n] = -\left( \frac{m}{2} + n \right) \phi_{m+n}$
\end{enumerate}
\end{prop}
\begin{proof}
The first relation :
\beqns
  [\phi_m, \psi_n] &=& \sum_j \binom{m+ \Delta_\phi - 1}{j} \left( \phi_{(j)}\psi \right)_{m+n}
\eeqns
where $\phi_{(0)}\psi = \langle \phi, \psi \rangle$ and $\phi_{(j)}\psi = 0$ $\forall j\neq 0$.
Assuming the weight of $1$ to be $0$ means that we define its expansion as
$$ 1 = \sum_{n\in -0 + \corps{Z}} \delta_{n,0} z^{-n-0} \ ,$$
so that $1_{m+n}= \delta_{m+n,0}=\delta_{m,-n}$.

Thus we have
$$ [\phi_m , \psi_n ] = \langle \phi , \psi \rangle (1)_{m+n} =  \langle \phi, \psi \rangle \delta_{m,-n}$$
For the second relation, we use proposition \ref{Lacomm} to obtain
\beqns
 [L_m, \phi_n] &=& (m (\frac{1}{2}-1)-n)\phi_{m+n} \\
               &=& -\left( \frac{m}{2} + n \right) \phi_{m+n}
\eeqns
\end{proof}
\end{exe}

\begin{remarque}[Energy operator]
Setting $m=0$ in the expression of $[L_m, a_n]$, we obtain
\beqns
 [L_0,a_n] &=& \sum_j \binom{1}{j} \left( L_{(j)}a \right)_n \\
           &=&  \left( L_{(0)}a \right)_n + \left( L_{(1)}a \right)_n
\eeqns
If $a$ is an eigendistribution of weight $\Delta_a$ with respect to $L$, then 
\beqns
 L_{(1)}a &=& \Delta_a a \\
 L_{(0)}a &=& \partial a
\eeqns
and hence
\beqns 
  [L_0, a_n] &=& (\partial a)_n + \Delta_a a_n \\
             &=& -(n+\Delta_a) a_n + \Delta_a a_n \\
             &=& -n a_n
\eeqns
Therefore,
$L_0$ is often called the \textit{energy operator}\index{energy operator}.
\end{remarque}

%
%

\chapter{Vertex algebras}



First introduced by Richard Borcherds in a different context, the structure of Vertex Algebra turns out to provide a rigorous mathematical formulation of the conformal theory of chiral (i.e. holomorphic) fields in two dimensions. Even though it \textit{might} be more natural to investigate physical aspects first and then study its intrinsic algebraic properties, we deliberately adopt the opposite point of view here.  We will start from the main definition of Vertex Algebra and then derive its properties in light of the algebraic structures studied in the previous chapters.

\section{Associative, commutative and unital algebra}

We first review the algebraic structure appearing as the ``classical'' limit of a vertex algebra : the structure of \textit{associative, commutative and unital algebra}\index{associative, commutative and unital algebra}. We formulate its definition in a way that serves our purpose.

\begin{defin}
An associative, commutative and unital algebra is defined as a linear space $V$ on $\corps{C}$, endowed with a linear map $Y\,:\, V\to \text{End}(V)$ satisfying
$$[Y(v),Y(w)]=0 \quad \forall v,w\in V \ .$$  The unit element, denoted by $\mathbf{1}$, is such that
$$ Y(\mathbf{1})=\text{I}_V \quad \text{and} \quad Y(v)\mathbf{1}=v \quad \forall v\in V \ .$$
\end{defin}
\begin{prop} This formulation is equivalent to the standard one.
\end{prop}  
\begin{proof} We can indeed define the product law by $v\cdot w = Y(v) w$, $\forall v,w\in V$. 
Since 
$$ v \cdot \mathbf{1} = Y(v) \mathbf{1} = v = \text{I}_V v = Y(\mathbf{1}) v = \mathbf{1} \cdot v \ ,$$
$\mathbf{1}$ is the unit element.
The following
\begin{eqnarray*} 
v \cdot w &=& Y(v)Y(w) \mathbf{1} \\
&=& Y(w)Y(v) \mathbf{1} = w \cdot v 
\end{eqnarray*}
shows that the product is commutative, and 
\begin{eqnarray*} v\cdot (w\cdot x) &=& Y(v)Y(w)x \\
&=& Y(v)Y(x)w \\
&=& Y(x)Y(v) w \\
&=& Y(Y(v)w)x = (v\cdot w)\cdot x 
\end{eqnarray*}
shows associativity.
\end{proof}

\section{Two equivalent definitions}
The previous formulation consisted to defining an associative, commutative and unital algebra by axioms satisfied by its set ${\{ Y(v) \mid v\in V \}}$ of left translations\index{left translation}.  A vertex algebra is defined in an analogous way, although its left translations depend on a formal variable $z$, so that, roughly speaking, a vertex algebra can be seen as an algebra with a parameter-dependant product law given by 
$$ a_z b = Y(a,z)b=\sum_{n\in \corps{Z}} a_{(n)}b z^{-n-1} \ .$$
We now state a first precise definition.
\begin{defin}
A vertex algebra\index{vertex algebra} is the data consisting of four elements $(V,\ket{0},T,Y)$ satisfying several axioms. For the data, $V$ is the superspace\index{superspace} $V=V_{\bar{0}}\oplus V_{\bar{1}}$, called the \emph{state space}\index{state space}, $\ket{0} \in V_{\bar{0}}$ is called the \textit{vacuum vector}\index{vacuum vector},  $T$ is a linear and even endomorphism acting on $V$, called the \emph{infinitesimal translation covariant operator}\index{translation covariant operator} and finally $Y$ is a linear and parity preserving map, from $V$ to the space of $End(V)$-valued formal distributions, acting as fields : 
$$ Y(\cdot, z) \,:\, V \to \text{End}(V)[[z,z^{-1}]] \,:\, a \to Y(a,z) \ ,$$
with, $\forall v\in V$, $Y(a,z)v \in V[[z]][z^{-1}]$. The map $Y$ is said to realize the \emph{state-field correspondence}\index{state-field correspondence}. For historical reasons, the field $Y(a,z)$ is called a \emph{vertex operator}\index{vertex operator}.\\
Those four elements $(V,\ket{0},T,Y)$ endow $V$ with the structure of vertex algebra when the following axioms are satisfied :
\begin{enumerate}
\renewcommand{\labelenumi}{(\theenumi)}
\item Vacuum\index{vacuum axioms} :  
\beqn
Y(\ket{0},z) &=& \text{Id}_V \\
Y(a,z)\ket{0} &=& a \mod z V[[z]]\\
T \ket{0} &=& 0 \ ,
\eeqn
\item translation covariance\index{translation covariance axiom} :
\beqn
[T,Y(a,z)]=\partial_z Y(a,z) \ ,
\eeqn
\item Locality\index{locality axiom} : $\{Y(a,z) \mid a\in V \}$ is a local family of fields, i.e. for $a,b\in V$,
\beqn
  (z-w)^n [Y(a,z),Y(b,w)] = 0 \quad n>>0 \ ,
\eeqn
where $[Y(a,z),Y(b,z)]\doteq Y(a,z)Y(b,w) - p(a,b)Y(b,w)Y(a,z)$.
\end{enumerate}
\end{defin}

Those axioms can be written in terms of Fourier modes.
\begin{prop}
Let $Y(a,z)=\sum_{n} a_{(n)} z^{-n-1}$. The space $V$ is endowed with an $n$-product defined by $a_{(n)}b=a_{(n)}(b)$.  Then the axioms of a Vertex Algebra read
\begin{enumerate}
\renewcommand{\labelenumi}{\emph{(\theenumi)}}
\item Vacuum\index{vacuum axioms} : 
\beqn
 \ket{0}_{(n)}a &=& \delta_{n,-1} a \quad (n\in \corps{Z}) \\
 a_{(n)} \ket{0} &=& \delta_{n,-1} a \quad (n\in \corps{Z}^+ \cup \{-1\})
\eeqn
\item translation covariance\index{translation covariance axiom} :
\beqn
 [T,a_{(n)}] = -n a_{(n-1)} \quad (n\in \corps{Z})
\eeqn
\end{enumerate}
\end{prop}
\begin{proof}
Applied to $a\in V$, the vacuum axiom gives :
\beqns
Y(\ket{0},z) a = \sum_{n\in \corps{Z}} \left( \ket{0}_{(n)}a \right) z^{-n-1} = a = \sum_{n\in\corps{Z}} a \delta_{n,-1} z^{-1-n}
\eeqns
so that $\ket{0}_{(n)}a = \delta_{n,-1} a \quad (n\in \corps{Z})$.  

The second axiom :
\beqns
\sum_{n\in\corps{Z}} a_{(n)}\ket{0} z^{-n-1} = a \mod z V[[z]]
\eeqns
so that $a_{(n)} \ket{0} = \delta_{n,-1} a \quad (n\in \corps{Z}^+ \cup \{-1\})$.

translation covariance reads :
\beqns
 \sum_n [T,a_{(n)}] z^{-1-n} = \sum_n (-n a_{(n-1)}) z^{-1-n}
\eeqns
so $[T,a_{(n)}] = -n a_{(n-1)} \quad (n\in \corps{Z})$.
\end{proof}

\begin{remarque}
When $Y(a,z)$ are formal series in positive powers of $z$, the vertex algebra is said to be \textit{holomorphic}\index{holomorphic}.
\end{remarque}

\begin{remarque}\label{actionTa}
The condition $T\ket{0}=0$ and the translation covariant axiom\index{translation covariant operator} permit to express $T$ by
$$ T(a) = a_{-2} \ket{0} $$
\begin{proof}
Applying the axiom on $\ket{0}$, we get
\beqns
 \sum_{n\in\corps{Z}} T(a_{(n)}\ket{0}) z^{-1-n} = \sum_{n\in\corps{Z}} a_{(n)} (-1-n) z^{-n-2}
\eeqns
Now, since $Y(a,z)\ket{0} \in V[[z]]$, we can set $z=0$ to both members, so that
\beqns
T(a_{(-1)}\ket{0}) = T(a) = a_{(-2)}\ket{0} \ .
\eeqns
\end{proof}
\end{remarque}

\begin{remarque}[Equivalent definition]
As a consequence, a second definition of a vertex algebra\index{vertex algebra} can be given by the data of three elements $(V,\ket{0},Y)$ satisfying the axioms of the first definition, but with the condition $T\ket{0}=0$ replaced by $T(a)=a_{(-2)}\ket{0}$ $\forall a\in V$. The latter defines the action of $T$ on $V$ and clearly implies $T\ket{0}=0$.  Indeed, for $a=\ket{0}$, we obtain $T\ket{0}=\ket{0}_{(-2)}\ket{0}=0$.
\end{remarque}

As an illustration, we will consider two examples of vertex algebra.
\begin{exe}
We start from an associative, commutative and unital superalgebra, endowed with an even derivation $T$.  Let $1$ be the unit element of $V$.  Then setting 
\beqns
\ket{0} &\doteq& 1 \\
Y(a,z) &\doteq & e^{zT} a \ ,
\eeqns
$V$ is then given a structure of vertex algebra.  Moreover,
$a_{(-1-n)}b = \left( \frac{T^{n} a}{n!} \right)b$ and $a_{(n)}b=0$ for $n \in \corps{Z}^+$.
\end{exe}
\begin{proof}
Let's check the vacuum axiom. As $T$ is an even derivation and $1$ the unit element,  
\beqns 
T(1)=T(1\cdot 1) &=& T(1)\cdot 1 + (-1)^{p(1)p(T)} 1\cdot T(1) \\
                 &=& 2 T(1)
\eeqns
then $T(1)=0$ so $T\ket{0}=0$.  Hence, 
\beqns
Y(\ket{0},z)b=\left(e^{zT} 1\right) b = b 
\eeqns
and $Y(\ket{0},z)=I_V$.  We also have
\beqns
Y(a,z) \ket{0} = \left(e^{zT} a \right) 1 = a \mod z V[[z]]
\eeqns
Let's check the translation covariant axiom :
\beqns
[T,Y(a,z)]b &=& T(Y(a,z) b) - Y(a,z) T(b) \\
            &=& T(e^{zT}a) b - (e^{zT}a) T(b) \\
            &=& T(ab + T(a)b z + T^2(a) b \frac{z^2}{2}+\cdots)-aT(b) - T(a)T(b) \frac{z^2}{2} - \cdots \\
            &=& (T(a) + T^2(a) z + T^3(a) \frac{z^2}{2} + \cdots) b \\
            &=& \partial_z (a+ T(a) z + T^2(a) \frac{z^2}{2} + T^3(a) \frac{z^3}{6}+\cdots) b \\
            &=& (\partial_z Y(a,z))b
\eeqns
The locality axiom is clearly satisfied since $V$ is commutative, so that 
\beqns
[Y(a,z),Y(b,z)]=0
\eeqns
Finally, since 
$$Y(a,z)b=\sum_{n\in\corps{Z}}a_{(n)}b z^{-1-n}=\left( e^{zT} a \right) b \ , $$ we clearly have $a_{(-1-n)}b = \left( \frac{T^{n} a}{n!} \right)b$ and $a_{(n)}b=0$ for $n \in \corps{Z}^+$.
\end{proof}

\begin{exe}
Another example of vertex algebra is provided by taking $V$ as the space of $End(U)$-valued fields, where $U$ is a linear superspace.  
\begin{theo} Let $V \subset End(U)[[z,z^{-1}]]$ be a local family of $End(U)$-valued fields.  If $V$ contains the unit element $1$, if it is $\partial_z$-invariant and closed under all $n$-products ($n\in\corps{Z}$) between formal distributions and if $T$ is the infinitesimal translation generator, whose action on a $End(U)$-valued field is defined by $T(a(w))= \partial_w a(w)$, then $V$ is vertex algebra, whose vacuum vector is $1$.
\end{theo}
\begin{proof}
Let's check the vacuum and the translation covariant axioms in terms of Fourier modes.  By assumption, their action on $V$ corresponds to the $n$-products between formal distributions.

Vacuum axiom : for $n\geq 0$, we have
\beqns
 1_{(-1-n)}a(z) &=& :\frac{(\partial^n 1)}{n!}a(z) : = a(z) \delta_{n,0} \\
 1_{(n)}a(z) &=& \Res{x} (x-z)^n[1,a(z)] = 0
\eeqns
Hence $1_{(n)}a(z)=\delta_{-1,n} a(z) \quad \forall n\in\corps{Z}$. 

For $n\geq 0$ :
\beqns
 a(z)_{(n)}1 &=& 0 \\
 a(z)_{(-1-n)}1 &=& : \frac{(\partial_z^n a(z))}{n!} 1 : \\
                &=& \frac{(\partial_z^n a(z)_+)}{n!} 1 + p(a,1) 1 \frac{(\partial_z^n a(z)_-)}{n!} \\
                &=& \frac{\partial_z^n a(z)}{n!} \quad \text{since $p(1)=\bar{0}$}
\eeqns
hence $a(z)_{(n)}1 = \delta_{-1,n} a(z) \quad \forall n\in \corps{Z}^+ \cup \{-1\}$.  The third reads $T(1)=\partial_z (1)=0$.

translation covariant axiom :
\beqns
 [T,a_{(n)}]b &=& \partial_z \left( a(z)_{(n)}b(z) \right) - a(z)_{(n)}\partial_z b(z) \\
              &=& -na(z)_{(n-1)}b(z) \\
\eeqns

Let's check the locality axiom. We proceed in successive steps.  The first step consists to using the generalized $n$-product formula to obtain
\beqns
Y(a(w),x) b(w) &=& \sum_{n\in\corps{Z}} a(w)_{(n)}b(w) x^{-1-n} \\
               &=& \Res{z} \left( a(z)b(w) i_{z,w} \sum_{n\in\corps{Z}} x^{-1-n} (z-w)^n - p(a,b) b(w) a(z) i_{w,z} \sum_{n\in\corps{Z}} x^{-1-n} (z-w)^n   \right) \\
               &=& \Res{z} \left( a(z)b(w) i_{z,w} \delta(z-w,x) - p(a,b) b(w) a(z) i_{w,z} \delta(z-w,x) \right)
\eeqns
The second step consists to computing the expression of the commutator $[Y(a(z),x),Y(b(w),y)]$ when applied on an element $c(w)$ of $V$ :  
\beqns
  [Y(a(z),x),Y(b(w),y)]c(w) = \underbrace{Y(a(z),x)\left( Y(b(w),y) c(w)\right)}_{(*)} - p(a,b) \underbrace{Y(b(w),y) \left( Y(a(z),x)c(w)\right)}_{(**)}
\eeqns
Using the expression above, we get :
\beqns
  (*) &=& Y(a(z),x)d(w) \\
      &=& \Res{z_1} \left( a(z_1) d(w) i_{z_1,w} \delta(z_1 - w,x) - p(a,bc) d(w) a(z_1) i_{w,z_1} \delta(z_1 - w,x) \right) 
\eeqns
where
\beqns
    d(w) &=& Y(b(w),y) c(w) \\
         &=& \Res{z_2} \left( b(z_2) c(w) i_{z_2,w} \delta(z_2-w,y) - p(b,c) c(w) b(z_2) i_{w,z_2} \delta(z_2-w,y)) \right)
\eeqns
Inserting the latter in the above expression, we get
\beqns
   (*) &=& \Res{z_1}\Res{z_2} \left[ \left( a(z_1)b(z_2) c(w) i_{z_2,w} i_{z_1,w} - p(b,c)a(z_1)c(w)b(z_2)i_{w,z_2} i_{z_1,w} \right.\right. \\
      &&  \left. -p(a,b)p(a,c) b(z_2)c(w) a(z_1)i_{z_2,w} i_{w,z_1} + p(a,b)p(a,c)p(b,c) c(w) b(z_2) a(z_1) i_{w,z_2} i_{w,z_1}\right) \\         
      &&  \left. \quad\quad\quad\quad\quad\quad\quad\quad\quad\quad\quad\quad\quad\quad\quad\quad\quad\quad\quad\quad\quad\quad \times \quad \delta(z_2-w,x)\delta(z_1-w,y)\right]
\eeqns
Similarly,
\beqns
   (**) &=& \Res{z_1}\Res{z_2} \left[ \left( b(z_2)a(z_1) c(w) i_{z_1,w} i_{z_2,w} - p(a,c)b(z_2)c(w)a(z_1)i_{w,z_1} i_{z_2,w} \right.\right. \\
      &&  \left. -p(b,a)p(b,c) a(z_1)c(w) b(z_2)i_{z_1,w} i_{w,z_2} + p(b,a)p(b,c)p(a,c) c(w) a(z_1) b(z_2) i_{w,z_1} i_{w,z_2}\right) \\         
      &&  \left. \quad\quad\quad\quad\quad\quad\quad\quad\quad\quad\quad\quad\quad\quad\quad\quad\quad\quad\quad\quad\quad\quad \times \quad \delta(z_1-w,x)\delta(z_2-w,y)\right]
\eeqns
Inserting those expression in the starting one,
\beqns
 [Y(a(z),x),Y(b(w),y)]c(w) &=& (*) - p(a,b) (**) \\
                           &=& \Res{z_1} \Res{z_2} \lbrace \left( [a(z_1),b(z_2)]c(w) i_{z_1,w} i_{z_2,w} \right. \\
                           &&  \left. - p(a,c)p(b,c) c(w) [a(z_1),b(z_2)] i_{w,z_1} i_{w,z_2}\right) \delta(z_1-w,x)\delta(z_2-w,y)\rbrace
\eeqns
The last step consists to showing that
\beqns
  (x-y)^n [Y(a(z),x),Y(b(w),y)]c(w)=0 \quad n>>0
\eeqns
But by assumption, $(a(z_1),b(z_2))$ is a local pair, i.e.
\beqns
(z_1 - z_2)^n [a(z_1),b(z_2)]=0 \quad n>>0
\eeqns
Now the trick is to use the identity
$$ x-y = (z_1 - z_2) - ((z_1-w) - x) + ((z_2-w) - y) \ .$$
Then we can expand $(x-y)^n$ in powers of $(z_1 - z_2)$, $((z_1-w) - x)$ and of $((z_2-w) - y)$ (finite sum).  Inserting this expansion in
\beqns
  (x-y)^n [Y(a(z),x),Y(b(w),y)]c(w)=0 
\eeqns
and using the computation of the commutator from the previous steps, we see that the terms in $(z_1 - z_2)^n$ will vanish by assumption, and the others, in non vanishing powers of $((z_1-w) - x)$ and/or $((z_2-w) - y)$, will also vanish by the presence of the distributions $\delta(z_1-w,x)$ and $\delta(z_2-w,y)$.
\end{proof}
\end{exe}

\begin{remarque}[Vertex Operator Algebra]
A distinguished vector $\nu$ arises in vertex algebras ($V$, say) studied in conformal theories. This (even) vector $\nu$ is such that the associated field $Y(\nu,z)=\sum_{n\in\corps{Z}}L_n z^{-n-2}$ is a Virasoro field with central charge $c$ (i.e. its Fourier coefficients span a Virasoro algebra with the central element acting as $cI_V$ ($c\in \corps{C}$) on $V$), for which $L_{-1}=T$ and $L_0$ is diagonalizable on $V$.  Such a vector is a called a \textit{conformal vector}\index{conformal vector}.  The vertex algebra is in turn called a \textit{conformal vertex algebra of rank $c$} and is often (depending on the literature) referred to as a \textbf{vertex operator algebra}\index{vertex operator algebra}.  The field $Y(\nu,z)$ is called an \textit{energy-momentum field}\index{energy-momentum field} of the vertex algebra $V$.  
\end{remarque}

\section{Unicity, $n$-products theorems}
In order to state those theorems, we will need the following results. 
\subsection{Preliminary results}
We need the following lemma to prove in an easier way that certain (possibly elaborate) expansions are in fact expansions in positive powers with respect to the variable considered and also to determine their \textit{initial} value.
\begin{lemme}\label{autrelemmesuper}
Let $V$ be a linear space and $\ket{0}$ a vector in $V$.  Let $a(z)$ and $b(z)$ two $End(V)$-valued fields such that 
\beqns
a_{(n)} \ket{0} &=& 0 \\
b_{(n)} \ket{0} &=& 0
\eeqns
for $n\in \corps{Z}^+$.  Then, $\forall N\in \corps{Z}$, $a_{(N)}b \ket{0}$ is an expansion in non negative powers of $z$, whose constant term is given by
$$ a(z)_{(N)}b(z) \ket{0} |_{z=0} =  a_{(N)}b_{(-1)} \ket{0} \ . $$
\end{lemme}
\begin{proof}
Let $N\geq 0$.  In that case,
\beqns
  a(z)_{(N)}b(z)\ket{0} &=& \Res{x} [a(x),b(z)](x-z)^N \ket{0} \\
                        &=& \Res{x} \sum_{m,n\in\corps{Z}} [a_{(m)},b_{(n)}] x^{-1-m} z^{-1-n} \sum_{j=0}^N \binom{N}{j} x^j (-z)^{N-j} \ket{0} \\
                        &=& \Res{x} \sum_{m,n,j} (-1)^{N-j} \binom{N}{j}[a_{(m)},b_{(n)}] x^{j-m-1} z^{N-j-n-1} \ket{0}\\
                        &=& \sum_{n,j} (-1)^{N-j} \binom{N}{j} [a_{(j)},b_{(n)}] z^{N-j-n-1} \ket{0}\\
                        &=& \sum_{n<0} \sum_{j=0}^N \binom{N}{j} a_{(j)} b_{(n)} z^{N-j-n-1} \ket{0}
\eeqns
Since $n<0$ and $0<j<N$, $N-j-n-1 \geq 0$.  The constant term is obtained for $j=N$ and $n=-1$, which proves the lemma for $N\geq 0$.

To prove the other case, let's compute, with $N\geq 0$ :
\beqns
  a(z)_{(-1-N)}b(z) \ket{0} &=& \frac{\partial_z^N a(z)_+}{N!} b(z) \ket{0} + p(a,b) b(z) \frac{\partial_z^N a(z)_-}{N!} \\
                            &=& \frac{\partial_z^N a(z)_+}{N!} b(z)_+ \ket{0}
\eeqns
As expected, this expression contains positive powers in $z$ alone.  In addition, since 
$$ \partial^N a(z)_+ = \sum_{m\geq N} m(m-1)\cdots(m-N+1) a_{(-1-m)} z^{m-N} \ ,$$
we can write
\beqns
  \frac{\partial_z^N a(z)_+}{N!} b(z)_+ \ket{0} &=& \sum_{j\geq 0} \sum_{m\geq N} \frac{m(m-1)\cdots(m-N+1)}{N!} a_{(-1-m)} b_{(-1-j)} z^{m-N+j}
\eeqns
The constant term is obtained for $j=0$ and $m=N$, which concludes the proof.
\end{proof}

The following simple lemma turns out to be very convenient to prove equalities involving possibly elaborate  infinite expansions (in positive powers though) at both sides.
\begin{lemme}\label{superrrr}
Let $A$ a linear operator on a linear space $V$.  The differential equation
$$ \frac{df(z)}{dz} = A f(z) \ , \quad f(z)\in V[[z]] $$
admits a unique solution, given an initial condition $f(0)=f_0$.
\end{lemme}
\begin{proof}
Setting
$$f(z)=\sum_{i\in\corps{Z}^+} f_i z^i \quad , \ f_i \in V$$
and inserting this equality in the differential equation, we easily obtain, after identification of the coefficients :
$$ f_{j+1} = \frac{A f_j}{j+1} \quad \text{where } j=0,1,\ldots $$
Recursively determined starting from the initial condition $f(0)=f_0$, the solution is unique.
\end{proof}

As already mentioned, the previous lemma is interesting in proving equalities involving possibly elaborate expansions (in positive powers though) at both sides.  The difficulty is divided in two steps, the first being to prove that both members are expansions in positive powers in a suitably chosen variable (making use of lemma \ref{autrelemmesuper} if needed), and the second being to prove that both sides obey the same initial condition (which is far easier) and the same differential equation (with respect to the chosen variable from the first step), which can often be guessed, at least for one of the members.  In fact, the real difficulty actually stems from \textit{finding} the candidate identity itself! We illustrate this in the following proposition.
\begin{prop}\label{propexpvertexvide}
Let $V$ be a vertex algebra.  Then the following equalities are satisfied :
\begin{enumerate}
\renewcommand{\labelenumi}{\emph{(\theenumi)}}
\item $Y(a,z) \ket{0} = e^{zT}a$
\item $e^{zT} Y(a,w) e^{-zT}=i_{w,z} Y(a,z+w)$
\item $\left( Y(a,z)_{(n)} Y(b,z) \right) \ket{0} = Y\left( a_{(n)}b,z \right) \ket{0}$
\end{enumerate}
\end{prop}
\begin{proof}
Members of equalities (1), (2) and (3) are formal expansions in positive powers in $z$ with coefficients belonging to $V$, to $End(V)[[w,w^{-1}]]$ and to $V$ respectively, which is trivial to show for the first case.  For the second one as well : recall that $i_{w,z}$ leads to an expansion in positive powers of $z/w$.  For the third case, we just make use of lemma \ref{autrelemmesuper}. 

Therefore, lemma \ref{superrrr} can be applied.  Those equalities are then proved as explained above, choosing $z$ as the suitable variable with respect to which the differential equation is written and the initial condition evaluated.

Let's prove (1).  Each member, denoted by $X(z)$, satisfies
$$ \frac{dX(z)}{dz} = TX(z) \ .$$
Indeed, for the member containing the exponential, this is trivial. For the other one, this comes from the translation covariant axiom and from $T\ket{0}=0$.  Moreover, the initial conditions are identical, for the term containing the exponential as for the other one, for which we apply the vacuum axiom. We indeed obtain $X(0)=a$ in both cases.

Let's prove (2).  Consider the RHS.  Let $X(z)=i_{w,z} Y(a,z+w)$.  We have
\beqns
  \frac{dX(z)}{dz} &=& \frac{d}{dz} \left( i_{w,z} Y(a,z+w) \right) \\
                   &=& i_{w,z} \frac{d}{dz}Y(a,z+w) \quad \text{since $[i_{w,z},\partial_z]=0$} \\
                   &=& i_{w,z} [T,Y(a,z+w)] \\
                   &=& [T,X(z)]
\eeqns
\beqns
  X(0) &=& \left.\left( i_{w,z} Y(a,z+w) \right)\right|_{z=0} \\
       &=& i_{w,0} Y(a,w) \\
       &=& Y(a,w)
\eeqns
For the LHS, let $X(z)=e^{zT} Y(a,w) e^{-zT}$.  We have
\beqns
  \frac{dX(z)}{dz} &=& \frac{d}{dz} \left( e^{zT} Y(a,w) e^{-zT} \right) \\
                   &=& \left( T e^{zT} \right) Y(a,w) e^{-zT} + e^{zT} Y(a,w) \left(-T \right) \\
                   &=& [T,X(z)]
\eeqns
\beqns
  X(0) &=& \left.\left( e^{zT} Y(a,w) e^{-zT} \right)\right|_{z=0} \\
       &=& Y(a,w)
\eeqns

Let's prove (3). The LHS is a series in positive powers of $z$, as it can be shown by applying lemma \ref{autrelemmesuper}.
Let's check the initial condition.  Let $X(z)=Y\left( a_{(n)}b,z \right)\ket{0}$.
\beqns
  X(0) &=& \sum_n \left( a_{(n)}b \right)_{(m)} \ket{0} \left.z^{-1-m}\right|_{z=0} \\
       &=& \left(a_{(n)}b\right)_{(-1)}\ket{0} \\
       &=& a_{(n)}b \quad \text{since $d_{(-1)}\ket{0}=d$ where $d=a_{(n)}b$}
\eeqns
For the LHS, setting $X(z)=\left( Y(a,z)_{(n)} Y(b,z) \right) \ket{0}$, we have
\beqns
 X(0) &=& \left.\left( Y(a,z)_{(n)} Y(b,z) \right) \ket{0}\right|_{z=0} \\
      &=& a_{(n)}\left(b_{(-1)}\ket{0}\right) \quad \text{cfr. lemme \ref{autrelemmesuper}} \\
      &=& a_{(n)}b \quad \text{since $b_{(-1)}\ket{0}=b$}
\eeqns
Both members satisfy the same differential equation. By setting $X(z)=Y\left( a_{(n)}b,z \right)\ket{0}$, we have
\beqns
 \frac{dX(z)}{dz} &=& \frac{d}{dz} Y\left( a_{(n)}b,z \right)\ket{0} \\
                  &=& [T,Y\left( a_{(n)}b,z \right)]\ket{0}] \\
                  &=& TX(z) \quad \text{since $T\ket{0}=0$}
\eeqns
Similarly, setting $X(z)=\left( Y(a,z)_{(n)} Y(b,z) \right) \ket{0}$, we find
\beqns
 \frac{dX(z)}{dz} &=& \frac{d}{dz} \left( \left( Y(a,z)_{(n)} Y(b,z) \right) \ket{0} \right) \\
                  &=& \partial_z Y(a,z)_{(n)} Y(b,z)\ket{0} + Y(a,z)_{(n)}\partial_z Y(b,z)\ket{0} \\
                  &=& [T,Y(a,z)]_{(n)}Y(b,z)\ket{0} + Y(a,z)_{(n)}[T,Y(b,z)]\ket{0}
\eeqns
By the generalized $n$-product formula, we can write, using $T\ket{0}=0$ :
\beqns
  \frac{dX(z)}{dz} &=&  [T,Y(a,z)]_{(n)}Y(b,z)\ket{0} + Y(a,z)_{(n)}[T,Y(b,z)]\ket{0} \\
                   &=&  \Res{z} \left( i_{x,z} (x-z)^n [T,Y(a,x)] Y(b,z) \ket{0} - i_{z,x} (x-z)^n p(a,b) Y(b,z) [T,Y(a,x)]\ket{0} \right. \\
                   &&  \left. + i_{x,z} (x-z)^n Y(a,x) [T,Y(b,z)]\ket{0} - i_{z,x} (x-z)^n p(a,b) [T,Y(b,z)]Y(a,x)\ket{0} \right) \\
                   &=& \Res{z} \left( i_{x,z} (x-z)^n T Y(a,x)Y(b,z) \ket{0} - i_{x,z} (x-z)^n Y(a,x) T Y(b,z) \ket{0} \right. \\
                   && - i_{z,x}(x-z)^n p(a,b) Y(b,z) T Y(a,x) \ket{0} +  i_{x,z} (x-z)^n Y(a,x) T Y(b,z) \ket{0} \\ 
                   && \left. - i_{z,x} (x-z)^n p(a,b) TY(b,z) Y(a,x) \ket{0} + i_{z,x}(x-z)^n p(a,b) Y(b,z) T Y(a,x) \ket{0} \right) \\
                   &=& \Res{z} \left( i_{x,z} (x-z)^n T Y(a,x)Y(b,z) \ket{0} - i_{z,x} (x-z)^n p(a,b) TY(b,z) Y(a,x) \ket{0} \right) \\
                   &=& T \Res{z}\left(i_{x,z} (x-z)^n  Y(a,x)Y(b,z) - i_{z,x} (x-z)^n p(a,b) Y(b,z) Y(a,x) \right) \ket{0} \\
                   &=& T Y(a,z)_{(n)}Y(b,z)\ket{0} \\
                   &=& TX(z)
\eeqns
which concludes the proof.

\end{proof}

\begin{remarque}
Note that the second equality of the previous proposition justifies, in the formal setting, the name ``infinitesimal generator of translation'' given to $T$.
\end{remarque}

\subsection{Unicity theorem}\index{unicity theorem}
\begin{theo}
Let $V$ be a vertex algebra.  The set of $\{ Y(a,z) \}$ being local, each field $Y(a,z)$ is uniquely determined by its ``initial'' value on $\ket{0}$.
\end{theo}
\begin{proof}
Let $B(z)$ an $End(V)$-valued field such that
\begin{enumerate}
\renewcommand{\labelenumi}{(\theenumi)}
\item $\left( B(z),Y(a,z) \right)$ is local pair for all $a\in V$,
\item $B(z)\ket{0} = Y(b,z) \ket{0}$.
\end{enumerate}
We have to prove that then $B(z)=Y(b,z)$. As stressed in the theorem, the locality assumption is determinant here.
Let's set $B_1(z)=B(z)-Y(b,z)$.  Then $(B_1(z),Y(a,z))$ is a local pair and $B_1(z)\ket{0}=0$. But locality means
\beqns
  (z-w)^N [B_1(z),Y(a,w)] = 0 \quad \text{for } N>>0
\eeqns
Being applied on $\ket{0}$, we obtain
\beqns
  (z-w)^N B_1(z)Y(a,w)\ket{0} &=& (\pm 1) (z-w)^N Y(a,w)B_1(z)\ket{0} \quad \text{for } N>>0 \\
                        &=& 0 \quad \text{since } B_1(z)\ket{0}=0
\eeqns
By the first property of proposition \ref{propexpvertexvide}, we have
\beqns
  Y(a,w)\ket{0} = e^{wT}a
\eeqns
The previous equality becomes, setting $w=0$ :
\beqns
  z^N B_1(z) a=0 &\Rightarrow& B_1(z)a = 0 \quad \forall a\in V \\
                 &\Rightarrow& B_1(z) = 0
\eeqns
hence the statement.
\end{proof}

\subsection{$n$-products theorem}\index{$n$-products theorem}
The following theorem is important for its corollaries.  We use the previous one to prove it.
\begin{theo}
Let $V$ be a vertex algebra. We have the following identity :
$$ Y\left( a_{(n)}b,z \right) = Y(a,z)_{(n)} Y(b,z) \ .$$
\end{theo}
\begin{proof}
Let's set $B(z)=Y(a,z)_{(n)}Y(b,z)$ and apply the unicity theorem.  The assumptions of the latter are clearly satisfied, since $(B(z),Y(a_{(n)}b,z))$ is local pair $\forall a,b\in V$ by Dong's lemma and $B(z)\ket{0}=Y(a_{(n)}b,z)\ket{0}$ by the third equality of proposition~\ref{propexpvertexvide}.
\end{proof}

\begin{coroll}
The following equalities are satisfied :
\begin{enumerate}
\renewcommand{\labelenumi}{\emph{(\theenumi)}}
\item $Y\left( a_{(-1)}b,z \right) = :Y(a,z)Y(b,z): $
\item $Y(Ta,z)= \partial_z Y(a,z)$
\item $\lbrack Y(a,z),Y(b,w)\rbrack =\sum_{j\geq 0} Y\left( a_{(j)}b,w \right) \frac{\partial^j_w \delta(z,w)}{j!}$
\end{enumerate}
\end{coroll}
\begin{proof}
The property (1) is a direct consequence of the previous theorem in the case $n=-1$.

To prove (2), we apply this theorem for $b=\ket{0}$ and $n=-2$ :
\beqns
 Y(a_{(-2)}\ket{0},z) &=& Y(a,z)_{(-2)} Y(\ket{0},z) \\
                      &=& : \left( \partial_z Y(a,z) \right) I_V : \\
                      &=& \partial_z Y(a,z)
\eeqns
But, $a_{(-2)}\ket{0}=T(a)$, according to remark~\ref{actionTa}, hence the equality.

The third equality comes from the theorem and the locality (decomposition theorem).
\end{proof}

\begin{remarque}
The second property of the previous corollary is equivalent to : $(Ta)_{(n)}b=-na_{(n-1)}b$ $\forall n\in \corps{Z}$.
\begin{proof}
Indeed, $Y(Ta,z)=\sum_n (Ta)_{(n)}b z^{-1-n}$ and $\partial_z Y(a,z) = \sum_n (-n)a_{(n-1)} z^{-1-n}$, hence the statement.
\end{proof}
\end{remarque}

\begin{remarque}
The third property shows that a vertex algebra can be seen as a field representation acting on itself.  Indeed, we see that by the state-field correspondence and by the fact that the sum is finite, $a_{(j)}b=0$ for $j>>0$ and $b\in V$, so that $Y(a,z)$ acts on $V$ as a End(V)-valued field. [Recall $a_{(j)}b \doteq a_{(j)}(b)$, i.e. $a_{(j)} \in$ End(V).]
\end{remarque}

\subsection{Borcherds' identity}
Let us have a look back at the Borcherds' identity\index{Borcherds' identity}.  We now show that it is merely a consequence of the $n$-product theorem, together with the locality axiom. The Borcherds' identity we saw at the very beginning translates in terms of formal distributions as follows.
\begin{eqnarray} Y(a,z)Y(b,w) i_{z,w} (z-w)^n - p(a,b) Y(b,w)Y(a,z)i_{w,z} (z-w)^n \nonumber \\
\quad\quad\quad\quad = \sum_{j\in\mathbb{Z}^+} Y(a_{(n+j)}b,w) \frac{\partial_w^j \delta(z,w)}{j!}  \label{borchccc}
\end{eqnarray}
for all $n\in\mathbb{Z}$ and $a,b\in V$.

We first prove \cite{DsK2} this identity and then that it is equivalent to Borcherds' identity.

\begin{prop}
Let $V$ be a vertex algebra. Then the identity \eqref{borchccc} is verified.
\end{prop}
\begin{proof}
Multiplying the LHS by $(z-w)^N$ for $N\geq -n$, we have
$$ (z-w)^N \ \cdot \ \text{LHS} = (z-w)^{n+N} [Y(a,z),Y(b,w)] \ . $$
Now since $\{Y(a,z)|a\in V\}$ is a local family, we can write
$$ (z-w)^N \big( (z-w)^n [Y(a,z),Y(b,w)] \big) = 0 $$
taking $N$ large enough (such that we also have $N\geq -n$).
By the decomposition theorem, we obtain 
$$ \text{LHS} = (z-w)^n [Y(a,z),Y(b,w)] = \sum_{j\geq 0} c^j (w) \frac{\partial_w^j \delta(z,w)}{j!} \ , $$
where 
\begin{eqnarray*}
c^j (w) &=& \text{Res}_z (z-w)^{n+j} [Y(a,z),Y(b,w)] \\
&=& Y(a,w)_{(n+j)}Y(b,w) 
\end{eqnarray*}
We conclude the proof by using $n$-product theorem according to which
$$ Y(a,w)_{(n+j)}Y(b,w) = Y(a_{(n+j)}b,w) \ .$$
\end{proof}

\begin{prop}
The identity \eqref{borchccc}, translated in terms of the Fourier modes, leads to the original Borcherds' identity.
\end{prop}
\begin{proof}
Making the following substitution
\begin{eqnarray*}
Y(a,z) &=& \sum_{k\in\mathbb{Z}} a_{(k)} z^{-k-1} \\
Y(b,w) &=& \sum_{l\in\mathbb{Z}} b_{(l)} w^{-l-1} \\
Y(a_{(n+j)}b,w) &=& \sum_{r\in\mathbb{Z}} \left( a_{(n+j)}b\right)_{(r)} w^{-r-1}
\end{eqnarray*}
and recalling that
\begin{eqnarray*}
i_{z,w} (z-w)^n &=& \sum_{j\geq 0} \binom{n}{j} (-w)^j z^{n-j} \\
i_{w,z} (z-w)^n &=& \sum_{j\geq 0} \binom{n}{j} z^j (-w)^{n-j}
\end{eqnarray*}
and that
\begin{eqnarray*}
\frac{\partial_w^j \delta(z,w)}{j!} &=& \sum_{s\in\mathbb{Z}} \binom{s}{j}z^{-s-1}w^{s-j} \ ,
\end{eqnarray*}
we apply both sides of \eqref{borchccc} to any element $c\in V$ and obtain
\begin{eqnarray*}
\sum_{k,l \in \mathbb{Z}} \sum_{j\geq 0} \binom{n}{j} (-1)^j a_{(k)}\left( b_{(l)}c \right) z^{n-j-k-1} w^{j-l-1}\quad\quad\quad\quad\quad\quad\quad\quad\quad\quad\quad\quad\quad\quad\quad\quad\quad\quad\quad\quad\quad\quad\quad\quad \\
- p(a,b) \sum_{k,l \in \mathbb{Z}} \sum_{j\geq 0} \binom{n}{j} (-1)^{n-j} b_{(l)}\left(a_{(k)}c \right) z^{j-k-1} w^{n-j-l-1} \quad\quad\quad\quad\quad\quad\quad\quad\quad\quad\quad \\
= \sum_{r,s \in \mathbb{Z}} \sum_{j\geq 0} \binom{s}{j} \left( a_{(n+j)}b\right)_{(r)}c z^{-1-s} w^{s-j-r-1}\quad\quad\quad\quad\quad\quad\quad\quad\quad
\end{eqnarray*}
Now we identify the coefficients of $z^{-1-s}w^{-1-t}$ by relabelling the indices in the following way.  We set $k=n+s-j$ and $l=j+t$ for the first term of the LHS, $k=j+s$ and $l=n-j+t$ for the second term of the LHS, and finally $r=s+t-j$ for the RHS.  We then obtain
\begin{eqnarray*}
\sum_{j\in \mathbb{Z}^+} \binom{s}{j} \left(a_{(n+j)}b \right)_{(s+t-j)}c \quad\quad\quad\quad\quad\quad\quad\quad\quad\quad\quad\quad\quad\quad\quad\quad\quad\quad\\ 
=  \sum_{j\in \mathbb{Z}^+} (-1)^j \binom{n}{j} \Big( a_{(n+s-j)}\left(b_{(j+t)}c\right) - (-1)^n b_{(n+t-j)}\left(a_{(j+s)}c\right)\Big) 
\end{eqnarray*}

\end{proof}

\section{Quasi-symmetry}
Recall the property of quasi-commutativity of the normal product.  Actually, this property comes from a more general result, namely the property of \textit{quasi-symmetry}\index{quasi-symmetry} in the framework of a vertex algebra.
\begin{theo}
Let $V$ be a vertex algebra. Then we have the property of \emph{quasi-symmetry} :
$$ Y(a,z)b=p(a,b) e^{zT} Y(b,-z) a \ . $$
\end{theo}
\begin{proof}
By the locality axiom,
\beqns
 (z-w)^N Y(a,z) Y(b,w) \ket{0} = (z-w)^N p(a,b) Y(b,w) Y(a,z) \ket{0} \quad \text{ for } N>>0
\eeqns
But $Y(b,w)\ket{0}=e^{wT}b$ and $Y(a,z)\ket{0}=e^{zT}a$.  So that :
\beqns
 (z-w)^N Y(a,z)e^{wT}b &=& (z-w)^N p(a,b) Y(b,w) e^{zT}a \\
                       &=& (z-w)^N p(a,b) e^{zT}e^{-zT}Y(b,w)e^{zT}a \\
                       &=& (z-w)^N p(a,b) e^{zT} i_{w,z}Y(b,w-z)a \quad \text{cfr. \ref{propexpvertexvide}(3)}
\eeqns
 By assumption, $Y(b,w-z)$ is a field, hence
 $$ Y(b,w-z)a \in V[[w-z]][(w-z)^{-1}] $$
 Choose $N$ large enough to have
 $$ (z-w)^N i_{w,z} Y(b,w-z)a \in V[[w-z]] $$
 Now, $i_{w,z}$ commutes with $(z-w)^N$.  Interchanging them and for the choice $N$ made above, we obtain
  $$ i_{w,z} \left( (z-w)^N  Y(b,w-z)a \right) = (z-w)^N  Y(b,w-z)a$$
  Inserting this expression in the previous equality, we have
  \beqns
    (z-w)^N Y(a,z)e^{wT}b = (z-w)^N p(a,b) e^{zT}Y(b,w-z)a
  \eeqns
  Setting $w=0$ and dividing by $z^N$, we get the expected result.
\end{proof}

\begin{remarque}
The property of quasi-symmetry can be written differently :
$$ a_{(n)}b = -p(a,b) (-1)^n \sum_{j\geq 0} \frac{(-T)^j}{j!} \left( b_{(n+j)}a \right) \quad \forall n\in\corps{Z} $$
This expression defines the so-called ``Borcherds' $n$-products''\index{Borcherds' $n$-products}.
\begin{proof}
  Starting from $Y(a,z)b=p(a,b)e^{zT}Y(b,-z)a$, we expand each member.  On the one hand,
  \beqns
    Y(a,z)b = \sum_{n} a_{(n)}b z^{-1-n} \ .
  \eeqns
  On the other hand,
  \beqns
    p(a,b)e^{zT}Y(b,-z)a &=& p(a,b) \sum_{j\geq 0} \frac{T^j}{j!} z^j \sum_{n} b_{(n)}a (-z)^{-1-n} \\
                         &=& \sum_n (-1)^{-1-n} p(a,b) \sum_{j\geq 0} z^{-1-n+j} \frac{T^j}{j!} \left( b_{(n)}a \right) \\
                         &=& \sum_{m,j}(-1)^{-1-m-j} p(a,b) z^{-1-m} \frac{T^j}{j!} \left(b_{(m+j)}a \right)\\
                         &=& \sum_m \left( -p(a,b) (-1)^m \sum_{j\geq 0} \frac{-T^j}{j!} \left(b_{(m+j)}a\right) \right) z^{-1-m}
  \eeqns
  where $n=m+j$.  Setting $m=n$ in the last line and equating both expanded members, we get the result by identifying the coefficients of $z^{-1-n}$.
\end{proof}
\end{remarque}

\begin{remarque}
Setting $n=-1$ in the Borcherds' $n$-products formula, we get the property of quasi-commutativity\index{normal product!quasi-commutativity} of the normal product, i.e.
$$ a_{(-1)} b - p(a,b) b_{(-1)} a = \int_{-T}^0 \lbrack a_\lambda b\rbrack  d\lambda \ ,$$
where $\lbrack a_\lambda b\rbrack =\sum_{j\geq 0} \frac{\lambda^j}{j!} a_{(j)}b$.  The normal product fails to be commutative (hence name of the property), due to the presence of the RHS, which is sometimes referred to as the ``quantum correction''.
\end{remarque}

\subsection{Quasi-associativity of the normal product}
\begin{prop}\label{quasiassoc}
The normal product is \textit{quasi-associative}\index{normal product!quasi-associativity} in the following sense :
\beqns
\left(a_{(-1)}b \right)_{(-1)}c - a_{(-1)} \left( b_{(-1)}c \right) &=& \sum_{j\geq 0} a_{(-j-2)}\left( b_{(j)}c \right) + p(a,b) \sum_{j\geq 0} b_{(-j-2)} \left( a_{(j)}c \right) \\
&=& \left( \int_0^T d\lambda \, a \right)_{(-1)} \lbrack b_\lambda c\rbrack  + p(a,b) \left( \int_0^T d\lambda \, b \right)_{(-1)} \lbrack a_\lambda c\rbrack  \ ,
\eeqns
\end{prop}
\begin{remarque}\label{interintegr}
The integrals above have to be interpreted as follows.  In the first place, we expand the $\lambda$-bracket in positive powers of $\lambda$.  The monomials in $\lambda$ are then placed at the left under the integral sign, before integrating inside the parentheses following the usual rules of integration.
\end{remarque}
\begin{proof}[Proof of proposition \ref{quasiassoc}] 
Let's prove the equality
\beqns
\left(a_{(-1)}b \right)_{(-1)}c - a_{(-1)} \left( b_{(-1)}c \right) &=& \sum_{j\geq 0} a_{(-j-2)}\left( b_{(j)}c \right) + p(a,b) \sum_{j\geq 0} b_{(-j-2)} \left( a_{(j)}c \right)
\eeqns
In the LHS, $\left(a_{(-1)}b \right)_{(-1)}$ is an endomorphism of $V$. More precisely, it shows up as the coefficient for $n=-1$ of the $End(V)$-valued field $Y(a_{(-1)}b,z)=\sum_{n}\left(a_{(-1)}b \right)_{(n)} z^{-n-1}$.  By the $n$-product theorem, we have $Y(a_{(-1)}b,z)=Y(a,z)_{(-1)}Y(b,z)=:Y(a,z)Y(b,z):$, so that $\left(a_{(-1)}b \right)_{(-1)}$ appears as the coefficient for $n=-1$ of the $End(V)$-valued field 
$:Y(a,z)Y(b,z):$, which can be written, setting $A(z)=Y(a,z)$ and $B(z)=Y(b,z)$ for clarity :
\beqns
 :A(z)B(z): &=& A(z)_+ B(z) + p(A,B) B(z) A(z)_- \\
            &=& \sum_n :AB:_{(n)} z^{z-n-1}
\eeqns
with (by proposition \ref{prodordcomp})
\beqns
 :AB:_{(n)} = \sum_{j=-1}^{-\infty} a_{(j)}b_{(n-j-1)} + p(a,b) \sum_{j=0}^\infty b_{(n-j-1)}a_{(j)} \ .
\eeqns
Therefore, we obtain for $n=-1$ :
\beqns
 \left(a_{(-1)}b \right)_{(-1)} = \sum_{j=-1}^{-\infty} a_{(j)}b_{(-2-j)} + p(a,b) \sum_{j=0}^\infty b_{(-2-j)}a_{(j)}
\eeqns
Hence
\beqns
 \left(a_{(-1)}b \right)_{(-1)} - a_{(-1)} b_{(-1)}  &=& \sum_{j=-1}^{-\infty} a_{(j)}b_{(-2-j)} + p(a,b) \sum_{j=0}^\infty b_{(-2-j)}a_{(j)} - a_{(-1)} b_{(-1)} \\
 &=& a_{(-1)}b_{(-1)} + \sum_{j=-2}^{-\infty} a_{(j)}b_{(-2-j)} + p(a,b) \sum_{j=0}^\infty b_{(-2-j)} a_{(j)} - a_{(-1)}b_{(-1)}
\eeqns
With $k=-j-2$, the first sum runs over non negative values of $k$ and we obtain
\beqns
\left(a_{(-1)}b \right)_{(-1)} - a_{(-1)} b_{(-1)} = \sum_{k=0}^\infty a_{(-k-2)} b_{(k)} + p(a,b) \sum_{k=0}^\infty b_{(-2-k)}a_{(k)} \ .
\eeqns
Applying both members to en element $c\in V$, we get the expected result.

Finally, let's prove
\beqn
\left(a_{(-1)}b \right)_{(-1)}c - a_{(-1)} \left( b_{(-1)}c \right) &=& \left( \int_0^T d\lambda \, a \right)_{(-1)} \lbrack b_\lambda c\rbrack  + p(a,b) \left( \int_0^T d\lambda \, b \right)_{(-1)} \lbrack a_\lambda c\rbrack  \label{azf68}
\eeqn
Start from the RHS.  Recalling remark \ref{interintegr}, the first term can be written
\beqn
 \left( \int_0^T d\lambda \, a \right)_{(-1)} \lbrack b_\lambda c\rbrack 
         &=& \sum_{j\geq 0} \left( \int_0^T \frac{\lambda^j}{j!} d\lambda a \right)_{(-1)}\left( b_{(j)}c \right) \nonumber \\
         &=& \sum_{j\geq 0} \left( \left. \frac{\lambda^{j+1}}{(j+1)!} a \right|_0^T \right)_{(-1)}\left( b_{(j)}c \right) \nonumber \\
         &=& \sum_{j\geq 0} \left( \frac{T^{j+1}}{(j+1)!} a \right)_{(-1)}\left( b_{(j)}c \right) \label{eeq4978}
\eeqn
We will show that
\beqns
 \sum_{j\geq 0} \left( \frac{T^{j+1}}{(j+1)!} a \right)_{(-1)}\left( b_{(j)}c \right) = \sum_{j=0}^\infty a_{(-j-2)}\left( b_{(j)}c \right)
\eeqns
But by the translation covariant axiom, 
 $$(Ta)_{(n)}=-na_{(n-1)} \ .$$  
 By induction, we obtain
\beqns
  (T^{j+1}a)_{(n)} = (-1)^{j+1} n(n-1)\cdots(n-j) a_{(n-j-1)}
\eeqns
For $n=-1$, we then have
\beqns
 (T^{j+1}a)_{(-1)} &=& (-1)^{j+1} (-1)(-2)\cdots(-j-1) a_{(-j-2)} \\
                   &=& (-1)^{j+1}(-1)^{j+1}(j+1)! a_{(-j-2)}
\eeqns
Hence
\beqns
 \left( \frac{T^{j+1}}{(j+1)!}a \right)_{(-1)} = a_{(-j-2)}
\eeqns
Applying both members to $\left(b_{(j)}c\right)$, summing over $j$ and using (\ref{eeq4978}), we obtain
\beqns
 \left( \int_0^T d\lambda \, a \right)_{(-1)} \lbrack b_\lambda c\rbrack = \sum_{j\geq 0} a_{(-j-2)} \left( b_{(j)}c \right)
\eeqns
We proceed in a similar way for the second term of the RHS of (\ref{azf68}).
\end{proof}
\begin{remarque}
The extra work done by introducing formal integrals might \textit{a priori} appear as motivated by a matter of taste.   But in fact, it permits us to reveal a crucial symmetry property in a straightforward way, as we will see it in the next section. The normal product fails to be associative, due to the presence of the RHS, which is referred to as the ``quantum correction'', as in the case of the quasi-commutativity.
\end{remarque}

\subsection{Normal product and Lie superbracket}\label{prodordcr}
Even though the normal product is quasi-associative, a Lie superbracket\index{normal product!Lie superbracket} can nevertheless be defined by the following :
$$ \lbrack a,b\rbrack  \doteq :ba: - p(a,b) :ba: \ . $$
\begin{remarque}
Even though the normal product fails to be associative, the bracket defined above is Lie thanks to the fact that the \textit{associator}\index{associator} is supersymmetric. 
\end{remarque}
\begin{prop}
The bracket defined by $ \lbrack a,b\rbrack  \doteq :ab: - p(a,b) :ba: $ is a Lie bracket.
\end{prop}
\begin{proof}
The property  $[b,a]=-p(a,b)[a,b]$ is clearly satisfied.  The only non trivial identity is Jacobi.
By definition,
\beqn
 \lbrack a,\lbrack b,c\rbrack \rbrack  &=& : a(:bc:): - p(b,c):a(:cb:): - p(a,b)p(a,c):(:bc:)a: \nonumber \\
                                       && \quad\quad\quad\quad\quad\quad\quad\quad\quad\quad\quad\quad + p(a,b)p(a,c)p(b,c) :(:cb:)a: \label{jacobprod1}\\
 \lbrack \lbrack a,b\rbrack ,c\rbrack  &=& :(:ab:)c: - p(a,b) :(:ba:)c: - p(a,c)p(b,c):c(:ab:): \nonumber\\
                                       && \quad\quad\quad\quad\quad\quad\quad\quad\quad\quad\quad\quad + p(a,c)p(b,c)p(a,b) :c(:ba:): \label{jacobprod2}\\
 p(a,b)\lbrack b,\lbrack a,c\rbrack \rbrack  &=& p(a,b) :b(:ac:): - p(a,c)p(a,b) :b(:ca:): - p(b,c):(:ac:)b: \nonumber\\
                                       && \quad\quad\quad\quad\quad\quad\quad\quad\quad\quad\quad\quad+ p(b,c)p(a,c):(:ca:)b: \label{jacobprod3}
\eeqn
To prove Jacobi, we have to show that if we subtract the first equation by the other two, we get zero.  We had indeed previously proved that, together with the property $[b,a]=-p(a,b)[a,b]$, the equivalence with Jacobi identity was assured.
Let's compute explicitly the quantity
\beqn
  \text{(\ref{jacobprod1})} - \text{(\ref{jacobprod2})} - \text{(\ref{jacobprod3})} 
           &=&  :a(:bc:): - :(:ab:)c: \label{jacobprod-1}\\
           && - p(b,c) \left( :a(:cb:): - :(:ac:)b: \right) \label{jacobprod-2}\\
           && - p(a,b)p(a,c) \left( :(:bc:)a: - :b(:ca:): \right) \label{jacobprod-3} \\
           && + p(a,b)p(a,c)p(b,c) \left( :(:cb:)a: - :c(:ba:): \right) \label{jacobprod-4}\\
           && + p(a,b) \left( :(:ba:)c: - :b(:ac:):  \right) \label{jacobprod-5}\\
           && + p(a,c)p(b,c) \left( :c(:ab:): - :(:ca:)b: \right) \label{jacobprod-6}
\eeqn
For a superbracket defined on an associative algebra, we would get the result directly.  Here, the normal product fails to be associative, but another property will do the trick : the \textit{associator} is supersymmetric.  More precisely, by proposition \ref{quasiassoc}, we see that
\beqns
 :a(:bc:): - :(:ab:)c: = p(a,b) \left( :b(:ac:): - :(:ba:)c: \right)
\eeqns
The \textit{associator} is then said supersymmetric under the interchange $a\leftrightarrow b$. 

Therefore,
$ \text{(\ref{jacobprod-1})} + \text{(\ref{jacobprod-5})} = 0 \ .$
Similarly, interchanging $a\leftrightarrow c$ :
\beqns
  \text{(\ref{jacobprod-2})} &=& -p(b,c)p(a,c) \left( :c(:ab:): - :(:ca:)b: \right) \ ,
\eeqns
hence
$ \text{(\ref{jacobprod-2})} + \text{(\ref{jacobprod-6})} = 0 \ .$
Finally, under $b\leftrightarrow c$ :
\beqns
 \text{(\ref{jacobprod-3})} -p(a,b)p(a,c)p(b,c) \left( :(:cb:)a: - :c(:ba:): \right) \ ,
\eeqns
hence
$ \text{(\ref{jacobprod-3})} + \text{(\ref{jacobprod-4})} = 0 \ ,$
which concludes the proof.
\end{proof}

\begin{remarque}
Recalling the quasi-commutativity property, we have 
$$ : a b : - p(a,b) : b a : = \int_{-T}^0 \lbrack a_\lambda b\rbrack  d\lambda \ .$$
Since the LHS was just shown to be a Lie bracket and the RHS is well defined for a Lie conformal superalgebra $\cal R$, we can now give the latter the structure of Lie superalgebra\index{Lie conformal algebra!Lie superalgebra}, by endowing it with the superbracket defined by $[a,b]\doteq\int_{-T}^0 \lbrack a_\lambda b\rbrack  d\lambda$, for $a,b\in\cal R$.
\end{remarque}

\section{Vertex Algebra and Lie Conformal Algebra}
It is the Fourier transform that permits us to go from a vertex algebra to the structure of Lie conformal algebra.  In other words, any vertex algebra $V$ can be endowed with a $\lambda$-bracket between elements of $V$, giving the latter the structure of Lie conformal algebra\index{Lie conformal algebra}.  It suffices to define the $\lambda$-bracket by :
$$ \Fzz \left( Y(a,z)b \right) = \sum_{j\geq 0} \frac{\lambda^j}{j!} a_{(j)}b \doteq \lbrack a_\lambda b\rbrack  $$
\begin{proof}
Indeed, 
\beqns
  \Fzz \left( Y(a,z)b \right) &=& \Res{z} e^{\lambda z} \left( \sum_{n\in\corps{Z}} a_{(n)}b z^{-n-1} \right) \\
                              &=& \Res{z} \sum_{j\geq 0} \frac{\lambda^j}{j!} z^j \sum_n a_{(n)}b z^{-n-1} \\
                              &=& \sum_{n,j} \frac{\lambda^j}{j!} a_{(n)}b \Res{z}z^{j-n-1} \\
                              &=& \sum_{j\geq 0} \frac{\lambda^j}{j!} a_{(j)}b \\
                              &=& [a_\lambda b]
\eeqns
\end{proof}
\begin{remarque}
Let's stress the fact that $a_{(j)}b$ should not be considered as the particular $j$-product between formal distributions.  In this context, $a_{(j)}b$ is identified to the $j$-product between elements $a,b$ of $V$ by the relation $a_{(j)}b = a_{(j)}(b)$, where $a_{(j)} \in End(V)$ is the $j^{\text{th}}$ coefficient of the $End(V)$-valued field $Y(a,z)=\sum_j a_{(j)} z^{-j-1}$.  To avoid confusion, we could write $a_{(j)}b=\left( Y(a,z) \right)_{(j)}b$.  Therefore, the $\lambda$-bracket defined as such should not be seen as the particular one defined between formal distributions, but as a bracket satisfying the properties \eqref{cmod01}-\eqref{jacol}, with $\partial=T$, which is verified in the following proposition.
\end{remarque}
\begin{prop}
The product defined above is indeed a $\lambda$-bracket.
\end{prop}
\begin{proof}
Consider the following equality :
\beqns
 \Res{z}e^{\lambda z}Y(Ta,z)b &=& \Res{z} e^{\lambda z} \sum_n (Ta)_{(n)}b z^{-n-1} \\
                              &=& \sum_{j\geq 0} \frac{\lambda^j}{j!} (Ta)_{(j)}b \\
                              &=& [Ta_\lambda b]
\eeqns
Using $Y(Ta,z)=\partial_z Y(a,z)$, we can write :
\beqns
  [Ta_\lambda b] &=& \Res{z}e^{\lambda z}Y(Ta,z)b \\
                 &=& \Res{z}e^{\lambda z} \partial_z Y(a,z) b \\
                 &=& -\lambda \Res{z} e^{\lambda z} Y(a,z) b \\
                 &=& - \lambda [a_\lambda b]
\eeqns
Hence we get \eqref{cmod01}.

In addition, as $a_{(j)}(Tb)=T(a_{(j)}b) - (Ta)_{(j)}b$, we can write :
\beqns
[a_\lambda Tb] &=& \sum_{j\geq 0} \frac{\lambda^j}{j!} (a)_{(j)}(Tb) \\ 
               &=& \sum_{j\geq 0} \frac{\lambda^j}{j!} T\left( (a)_{(j)}b \right) - \sum_{j\geq 0} \frac{\lambda^j}{j!} (Ta)_{(j)}b \\
               &=& T [a_\lambda b] - [Ta_\lambda b] \\
               &=& (T+\lambda) [a_\lambda b]
\eeqns
and we obtain \eqref{cmod02}.

Using the quasi-symmetry property, we obtain
\beqns
  [a_\lambda b] &=& \Fzz Y(a,z)b \\
                &=& p(a,b) \Fzz e^{zT}Y(b,-z)a \\
                &=& p(a,b) \Res{z} e^{(\lambda + T)z} Y(b,-z)a \\
                &=& p(a,b) \Fzd{z}{\lambda + T} Y(b,-z)a \\
                &=& p(a,b) \left( -\Fzd{z}{-\lambda-T}Y(b,z)a \right)\\
                &=& -p(a,b) [b_{-\lambda - T}a] \quad \text{by \ref{trfour1d}(3)}
\eeqns
hence \eqref{skcomml}.

Let's prove Jacobi \eqref{jacol}.  By the $n$-products theorem,
\beqns
  Y([a_\lambda b],z) &=& \sum_{j\geq 0} \frac{\lambda^j}{j!} Y(a_{(j)}b,z) \\
                     &=& \sum_{j\geq 0} \frac{\lambda^j}{j!} Y(a,z)_{(j)}Y(b,z) \\
                     &=& [Y(a,z)_\lambda Y(b,z)] 
\eeqns
Apply the LHS on $c\in V$, before applying $\Fzd{z}{\lambda+\mu}$ :
\beqns
  \Fzd{z}{\lambda+\mu} Y([a_\lambda b],z)c = [[a_\lambda b]_{\lambda + \mu}c] 
\eeqns
Doing the same on the RHS, we have
\beqns
  \Fzd{z}{\lambda+\mu} [Y(a,z)_\lambda Y(b,z)]c &=& [\Fzz Y(a,z),\Fzd{z}{\mu} Y(b,z)]c \quad \text{by \ref{trfour1d}(2)} \\
                                                &=& \Fzz Y(a,z) [b_\mu c] - p(a,b) \Fzd{z}{\mu} Y(b,c) [a_\lambda c] \\
                                                &=& [a_\lambda [b_\mu c]] - p(a,b) [b_\mu[a_\lambda c]]
\eeqns
Hence
\beqns
 [[a_\lambda b]_{\lambda + \mu}c] = [a_\lambda [b_\mu c]] - p(a,b) [b_\mu[a_\lambda c]]
\eeqns
or
\beqns
 [a_\lambda [b_\mu c]] = [[a_\lambda b]_{\lambda + \mu}c] + p(a,b) [b_\mu[a_\lambda c]]
\eeqns
and we indeed obtain \eqref{jacol}.
\end{proof}

\section{Third definition}
The notions of Lie conformal superalgebra, normal product and their respective properties were stated and proved in the preceding sections.   Those can all be gathered and permit to build a definition of a vertex algebra equivalent to the first two we have seen so far~\cite{DsK2}.
\begin{defin}\label{algvert3}
A vertex algebra is the data of five elements $(V,\ket{0},T,[\cdot_\lambda \cdot],:\cdot \, \cdot :)$, where
\begin{enumerate}
\renewcommand{\labelenumi}{(\theenumi)}
\item $(V,T,[\cdot_\lambda \cdot])$ is a Lie conformal algebra,
\item $(V,\ket{0},T,:\cdot \, \cdot :)$ is a unital and differential superalgebra, \textit{i.e.} $T$ is a superderivation, where the normal product $:\cdot \, \cdot:$ is quasi-commutative in the following sense :
\beqns
:ab: - p(a,b) :ba: = \int_{-\partial}^0 [a_\lambda b] d\lambda
\eeqns
and quasi-associative in the following sense :
\beqns
\left(a_{(-1)}b \right)_{(-1)}c - a_{(-1)} \left( b_{(-1)}c \right) = \left( \int_0^T d\lambda \, a \right)_{(-1)} \lbrack b_\lambda c\rbrack  + p(a,b)
\left( \int_0^T d\lambda \, b \right)_{(-1)} \lbrack a_\lambda c\rbrack  \ ,
\eeqns
\item the $\lambda$-bracket $[\cdot_\lambda \cdot]$ and the normal product $:\cdot \, \cdot :$ are related by the non-abelian Wick's formula, \textit{i.e.}
\beqns
    [a_\lambda : bc :] = : [a_\lambda b] c : + p(a,b) : b[a_\lambda c] : + \int_0^\lambda [[a_\lambda b]_\mu c] d\mu 
\eeqns
\end{enumerate}
\end{defin}

\begin{remarque}
As seen before, the associator $::ab:c:-p(a,b):b:ac::$ is supersymmetric under the interchange $a\leftrightarrow b$. Actually, the latter can be considered as a condition and replace the condition of quasi-associativity in definition \ref{algvert3} \cite{BK}.
\end{remarque}

We see that the notion of \textit{vertex algebra} can be defined with the help of the algebraic structures introduced in previous chapters. This chapter thus concludes by illustrating the way that all those notions get tangled with each other.

%
%
%

\printindex

\end{document}